\documentclass[11pt]{amsart}

\usepackage{amssymb} \usepackage{amsmath} \usepackage{amsthm}
\usepackage{latexsym} \usepackage{exscale}
\usepackage{epsf}
\makeindex
\newtheorem{thm}{Theorem}[section]
\newtheorem{lemma}[thm]{Lemma}
\newtheorem{prop}[thm]{Proposition}
\newtheorem{cor}[thm]{Corollary}
\newtheorem{defn}{Definition}[section]
\newtheorem{claim}{Claim}[section]
\DeclareMathOperator{\sgn}{sgn}
\DeclareMathOperator{\spn}{span}
\DeclareMathOperator{\re}{Re}
\numberwithin{equation}{section}
\numberwithin{figure}{section}

\newcommand{\lam}{\ensuremath{\lambda}}

\newcommand{\evans}{\ensuremath{D(\lambda)}}

\newcommand{\C}{\ensuremath{\mathbf{C}}}

\newcommand{\R}{\ensuremath{\mathbf{R}}}

\newcommand{\E}{\ensuremath{\mathcal{E}}}

\newcommand{\z}{\ensuremath{\zeta}}

\newcommand{\pf}{\begin{proof}}

\newcommand{\foorp}{\end{proof}}

\newcommand{\ip}[2]{\ensuremath{\langle #1,#2\rangle}}

\begin{document}
\title[One-dimensional Stability]{On the One-dimensional Stability of Viscous Strong Detonation Waves}
\author{Gregory Lyng}
\address{Department of Mathematics, University of Michigan, Ann Arbor, MI 48109}
\email{glyng@umich.edu}
\author{Kevin Zumbrun}
\address{Department of Mathematics, Indiana University, Bloomington, IN 47405}
\email{kzumbrun@indiana.edu}
\date{February 12, 2003 \\ Revised: \today}

\begin{abstract} Building on  Evans function techniques developed to study the stability of viscous shocks, we examine the stability of viscous strong detonation wave solutions of the reacting Navier-Stokes equations. The primary result, following \cite{AGJ,GZ}, is the calculation of a \emph{stability index} whose sign determines a necessary condition for spectral stability. We show that for an ideal gas this index can be evaluated in the ZND limit of vanishing dissipative effects. Moreover, when the heat of reaction is sufficiently small, we prove that strong detonations are spectrally stable provided the underlying shock is stable. Finally, for completeness, the stability index calculations for the nonreacting Navier-Stokes equations are included. 
\end{abstract}
\maketitle
\tableofcontents
%
\section{Introduction and Preliminaries}\label{Intro}
%
%
%
\subsection{Introduction}
Laboratory and numerical experiments indicate that detonations have quite sensitive stability properties. Indeed steady planar detonations subjected to one-dimensional longitudinal perturbations may change form to ``galloping'' detonations in which the velocity fluctuates periodically in time. Such detonations have been been predicted numerically \cite{FickWood} and observed experimentally in various settings by \cite{GordMooHarp}, \cite{Manson}, and \cite{Mundy}. Another instability, this one with 3-dimensional structure, is the ``spinning detonation'' long-known in lab experiments \cite{CampWood1,CampWood2} and more recently captured numerically for the ZND model in both \cite{BMR} and \cite{KasStew}. A three-dimensional perturbation of a steady detonation wave propagating down a tube with a circular cross-section may bifurcate to a wave with a complex rotating structure which traces a helical path along the boundary of the tube. This structure is typically followed by localized regions of extremely high pressure. 

Due to this sensitivity and the complicated, nonlinear nature of reacting gas dynamics (which includes such difficulties as turbulence, boundary layers, and complex chemical interactions), stability analyses of detonation waves are largely numerical studies of the ZND model (which neglects dissipative effects) as in e.g. \cite{BMR},\cite{LeeStew},\cite{SS1,SS2,SS3}, or are  restricted to various incarnations of the Majda or Majda-Rosales\footnote{For the remainder of the paper, we refer to all models with the simplifying feature of scalar kinetics as the Majda model.} models (unphysical analogues of Burgers equation) as in \cite{LLT},\cite{LY},\cite{LYi}, \cite{Li1,Li2,Li3,Li4,Li5},\cite{RV},and \cite{Sz}. Our approach, utilizing the Evans function, allows the treatment of the reacting Navier-Stokes equations and yields an explicitly computable quantity known as the \emph{stability index}.

In Eulerian coordinates the reacting Navier-Stokes equations modeling the simplest possible one-step chemical reaction can be written
\begin{eqnarray}
\rho_t + (\rho u)_x & = & 0, \label{eq:mass} \\
(\rho u)_t + (\rho u^2+p)_x & = & (\nu u_x)_x, \label{eq:momentum} \\
\tilde{\E}_t+[\rho u\tilde{\E}+up]_x & = & (\theta T_x)_x+(q\rho d Y_x)_x+(\nu uu_x)_x,\label{eq:energy} \\
(\rho Y)_t+(\rho uY)_x & = & (\rho dY_x)_x -k\rho Y\varphi(T). \label{eq:progress}
\end{eqnarray}
In (\ref{eq:mass})-(\ref{eq:progress}) and below we use, unless stated otherwise, the notations: $\rho$, $u$, $p$, $\E$, $T$, and $Y$ represent respectively density, velocity, pressure, total energy, temperature, and mass fraction of reactant. The use of the tilde denotes that the energy $\tilde{\E}=\rho(u^2/2+\tilde{e})$ is modified from the standard gas-dynamical energy $\E=\rho(u^2/2+e)=\rho E$ due to heat produced in the chemical reaction by  $\tilde{e}=e+qY$. The positive constants $\nu$, $\theta$, and $d$ represent the effects of viscosity, heat conductivity, and species diffusion. The positive constants $k$ and $q$ measure the rate of reaction and the heat released in reaction, respectively, and the form of the so-called ignition function $\varphi$ is discussed in detail below. The system is closed by specifying equations of state, $p=p(\rho,e,Y)$ and $T=T(\rho,e,Y)$.  We begin by assuming only that $p$ and $T$ are independent of $Y$, but for some portions of the analysis we shall assume further an ideal, polytropic gas, i.e
$$p(\rho,e)=\Gamma\rho e,\quad T(\rho,e)=c_v^{-1}e,$$
where $c_v$, the specific heat at constant volume, and $\Gamma$, known as the Gruneisen coefficient, are constants. 
Equations \eqref{eq:mass}-\eqref{eq:progress} are standard; a derivation can be found in \cite{Wi}.

\begin{figure}[ht]
\centerline{\epsfxsize8truecm\epsfbox{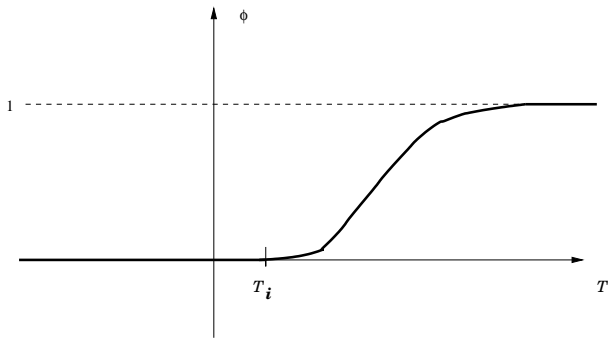}}
\caption{The Ignition Function $\varphi$}
\label{ignite}
\end{figure}

Often $\varphi$ is assumed to satisfy the Arrhenius law, so
$$\varphi(T)=\exp\left(-\frac{E_A}{RT}\right),$$
where $E_A$ is the activation energy and $R$ is the gas constant (assuming the ideal gas law, $R=c_v\Gamma$ ). However, nonvanishing of the exponential creates a problem known as the ``cold-boundary difficulty.'' Essentially, nonvanishing of $\varphi$ precludes the unburned state from being a rest point of the traveling wave ODE. In place of the Arrhenius kinetics we make the standard assumption that the smooth function $\varphi$ vanishes for temperatures below some ignition temperature, $T_i$, and is identically 1 for some larger value of $T$. This circumvents the cold boundary difficulty; see Figure~\ref{ignite}.
This model includes dissipative effects neglected by the ZND model, and allows for complete gas dynamical effects unlike the Majda model.  We remark that an artificially strictly parabolic multidimensional version is considered in the appendix of \cite{Z4}; here we include the additional difficulty of partial parabolicity. 
%
%
%
%
\subsection{Plan of the Paper}
In Section~\ref{Intro} we provide the relevant background material for our analysis. We first discuss the standard ZND model as a prelude to a discussion of the analysis of \cite{GS1}. ZND solutions are singular solutions in the context of geometric singular perturbation theory used therein. The structure of these singular solutions allows us to evaluate the stability index, which we do in Section~\ref{RNSsection}. The backgound material in Section~\ref{Intro} concludes with a description of the Evans function theory for the stability of viscous shock waves. Section~\ref{NSsection} contains the stability index calculations for the (nonreacting) Navier-Stokes equations. These prove useful in Section~\ref{RNSsection}, where we examine the reacting system \eqref{eq:mass}-\eqref{eq:progress}. In Appendix~\ref{appendix} we include a revised version of an appendix of \cite{Z4} which is used in the stability index calculations.
%
%
\subsection{The ZND Model}
Setting the constants $\nu,\theta,$ and $d$ equal to zero in \eqref{eq:mass}-\eqref{eq:progress} yields the ZND model introduced independently by Zeldovich, von Neumann, and D\"{o}ring. In their formulation dissipative effects are neglected and the reaction rate is assumed to be finite. This is a refinement of the early Chapman-Jouget model in which the reaction was assumed to take place instantaneously. The ZND model, in Eulerian coordinates, then has the form
\begin{eqnarray*}
\rho_t + (\rho u)_x & = & 0, \label{ZNDmass} \\
(\rho u)_t + (\rho u^2+p)_x & = & 0, \label{ZNDmomentum} \\
\tilde{\mathcal{E}}_t+((\tilde{\mathcal{E}}+p)u)_x & = & 0,\label{ZNDenergy} \\
(\rho Y)_t+(\rho uY)_x & = &  -k\rho Y\varphi(T). \label{ZNDprogress}
\end{eqnarray*}
As above the constant $k$ is the reaction rate, and $\varphi$ is the ignition function. 
Strong detonations in the ZND model are initiated by a (purely) gas dynamical shock, called the Neumann shock,  which heats the gas by compressing it. The increase in the temperature to a sufficiently high level ``turns on'' $\varphi$ and thus starts the reaction. Thus these waves have the structure of a gas dynamical shock followed by a reaction zone resolving to the final burned state.  This is seen in the characteristic ``detonation spikes'' in the temperature and pressure profiles in strong agreement with observed features in laboratory experiments. We'll see this structure in our discussion of \cite{GS1} below.
%
%
%
\subsection{Existence of Strong Detonations}
In \cite{GS1}, existence of traveling wave solutions of Equations~(\ref{eq:mass})-(\ref{eq:progress}) was studied using the techniques of geometric singular perturbation theory (GSPT). As the orientations of the singular manifolds constructed in that argument will play a role in our analysis, we recap the argument here. An interesting feature of the GSPT analysis of \eqref{eq:mass}-\eqref{eq:progress} is the recovery of the shock layer analysis of \cite{Gi}. 
\subsubsection{The Hugoniot Curve}
\begin{figure}[ht]
\centerline{\epsfxsize12truecm\epsfbox{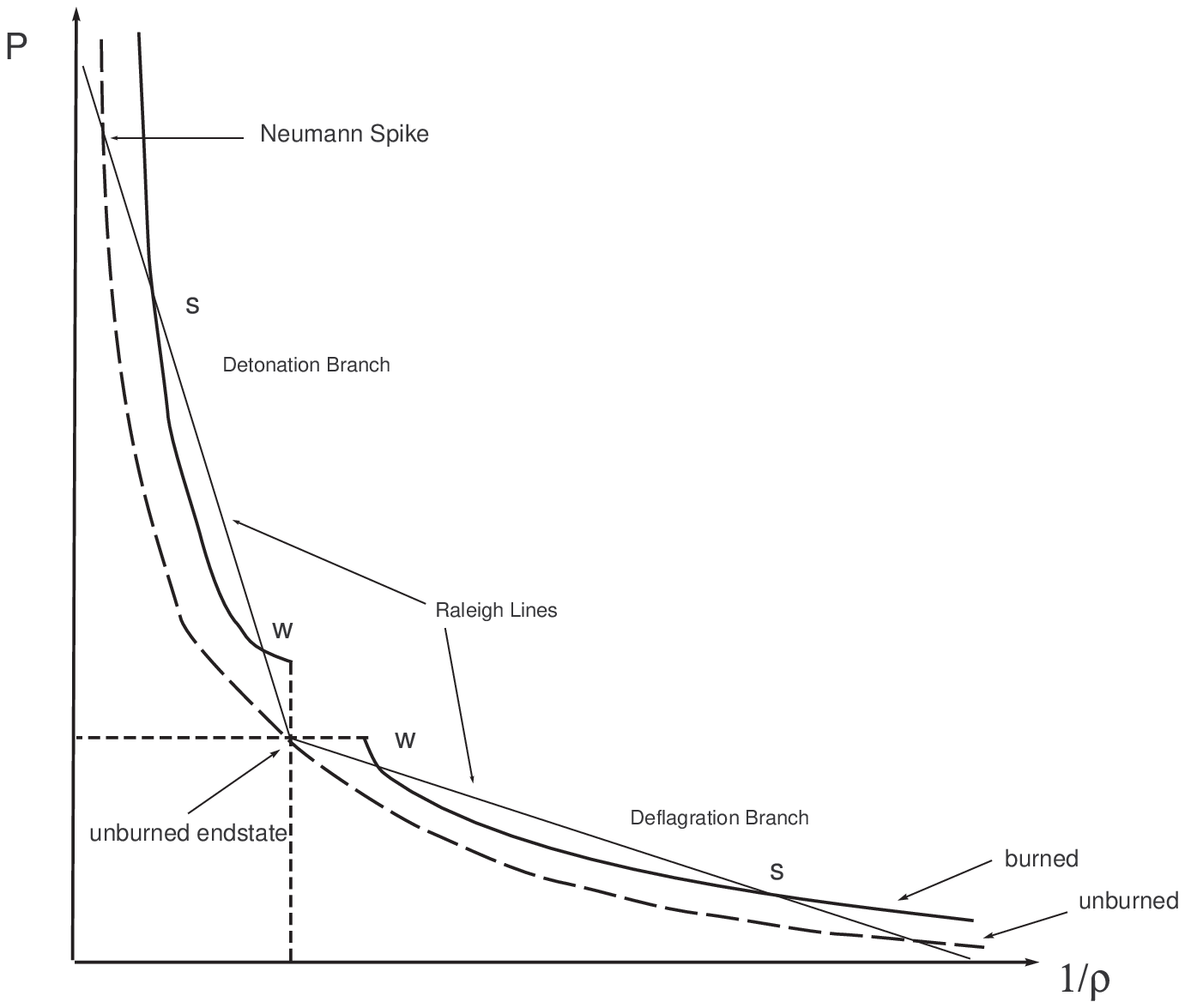}}
\caption{The Hugoniot Curve}
\label{hugoniot}
\end{figure}

Traveling wave solutions of Equations~(\ref{eq:mass})-(\ref{eq:progress}) are those which depend only on $\xi =x-st$. This ansatz reduces the system \eqref{eq:mass}-\eqref{eq:progress} to a system of ordinary differential equations~(ODEs). By Galilean invariance we may, without loss of generality, set $s=0$, so the system is 
\begin{eqnarray}
(\rho u)' & = & 0, \label{tw:mass} \\
(\rho u^2+p)' & = & (\nu u')', \label{tw:momentum} \\
\left(\rho u( \frac{u^2}{2}+\tilde{e})+up\right)' & = & (\theta T')'+(q\rho d Y')'+(\nu uu')',\label{tw:energy} \\
(\rho uY)' & = & (\rho dY')' -k\rho Y\varphi(T), \label{tw:progress}
\end{eqnarray}
where $'$ denotes differentiation with respect to $x$.
From the first equation (\ref{tw:mass}) it follows that the mass flux $m=\rho u$ has a constant value. Moreover each of (\ref{tw:momentum}) and (\ref{tw:energy}) can be integrated up once. We suppose that an unburned state $\rho_+,u_+,p_+,Y_+$ has been fixed at $+\infty$. Then the momentum equation integrates to
$$\rho u^2+p-(\rho_+u_+^2+p_+)=\nu u'.$$
For a possible connection to a burned state $\rho ,u,p,Y$ at $-\infty$, it is necessary that the state be a rest point of the ODE, or more precisely $\rho ,u,p,Y$ must satisfy
$$\rho u^2+p-(\rho_+u_+^2+p_+)= 0.$$
Searching for all such states leads to the expression
\begin{equation*}
p-p_+ = -m^2\left(\frac{1}{\rho}-\frac{1}{\rho_+}\right).
\end{equation*}    
This equation describes a line in the specific volume-pressure plane with slope $-m^2$. It is referred to as the \emph{Raleigh Line}. Similarly integrating the energy equation and searching for rest points again, leads to the equation for the \emph{Hugoniot Curve}
$$\tilde{e}-\tilde{e}_+=-\frac{1}{2}(p+p_+)\left(\frac{1}{\rho}-\frac{1}{\rho_+}\right).$$
The intersection of the Raleigh line and the Hugoniot curve in the specific volume-pressure plane determine the possible burned states corresponding to the fixed unburned state $\rho_+,u_+,p_+,Y_+$. An important distinction from the nonreacting case is the fact that there may be one, two, or no possible burned states. See Figure \ref{hugoniot}.
Another interesting feature of Figure~\ref{hugoniot} is the fact that the Hugoniot curve splits into two branches. This indicates that the conservation relations are compatible with two distinct types of processes, just as observed by early experimentalists. The compressive solutions are called \emph{detonations} while the expansive solutions are referred to as \emph{deflagrations}. Accordingly the two branches of the Hugoniot curve are called the detonation branch and deflagration branch. Each of these branches is further subdivided into two sections. Here we focus on the detonation branch, but similar characterizations apply to the deflagration branch. From the diagram it's clear that there is a unique value of $-m^2$ so that the Raleigh line is tangent to the Hugoniot curve. The point of tangency is called the Chapman-Jouget point, and a detonation connecting to the burned state identified by that point is called a Chapman-Jouget detonation. For values of $-m^2$ which are smaller than this unique value, the Raleigh line intersects the Hugoniot curve twice. The larger of these two intersections is the burned state corresponding to a \emph{strong detonation} while the smaller corresponds to the burned state of a \emph{weak detonation}. One distinction between these waves is the following: strong detonations satisfy the Lax characteristic condition while weak detonations are undercompressive. Here compressivity refers to the number of incoming characteristics. Finally, if $-m^2$ is too large, the Raleigh line and the Hugoniot curve do not intersect and there are no possible burned states which are compatible with the unburned state.


%
%
%
\subsubsection{GSPT Analysis of Detonation Waves}
We begin with the briefest of introductions to GSPT. Consider a system of singularly perturbed ODE
\begin{equation}
\begin{cases} \frac{dx}{dt}  =  f(x,y), \\ \epsilon\frac{dy}{dt}  =  g(x,y), \end{cases}
\tag{$\text{Slow}_{\epsilon}$} \label{slow}
\end{equation}
where $\epsilon$ is small. We call such a system, the \emph{slow} system. By rescaling the independent variable by $\tau =t/\epsilon$, we obtain the equivalent (when $\epsilon \ne 0$) \emph{fast} system
\begin{equation}
\begin{cases} \frac{dx}{d\tau}  =  \epsilon f(x,y), \\ \frac{dy}{d\tau}  =  g(x,y). \end{cases}
\tag{$\text{Fast}_{\epsilon}$}\label{fast}
\end{equation}
Looking at \eqref{slow} and \eqref{fast}, it's clear that there are then two distinguished limiting systems when $\epsilon =0$. They are the \emph{reduced} problem
\begin{equation}
\begin{cases} \frac{dx}{dt}  =  f(x,y), \\ 0  =  g(x,y), \end{cases}
\tag{$\text{Slow}_{0}$} \label{slow0}
\end{equation}
and the \emph{layer} problem
\begin{equation}
\begin{cases} \frac{dx}{d\tau}  =  0, \\ \frac{dy}{d\tau}  =  g(x,y). \end{cases}
\tag{$\text{Fast}_{0}$}\label{fast0}
\end{equation}
The basic idea then is to construct solutions of the original system as smooth perturbations of the composite orbits of the decoupled limiting systems \eqref{slow0} and \eqref{fast0}. For more details consult \cite{Fen} and \cite{Szm}.

With this framework in mind, we take up the analysis of \cite{GS1}. We note that here, following \cite{GS1}, we are assuming an ideal, polytropic gas, so that $p=R\rho T$, $e=c_vT$, and $\Gamma=\gamma-1=R/c_v$. We remark that existence for more general equations of state has been shown by different methods in \cite{HR}, but the singular manifolds constructed in the GSPT analysis \cite{GS1} play a key role in the evaluation of our stability condition. Namely, they contain the necessary geometric information about the profile. 
Using $m$ and the integrated versions of (\ref{tw:momentum}) and (\ref{tw:energy}), one finds by some further manipulation that
\begin{eqnarray*}
\nu u_x&  =& m(u-u_{\pm})+mR\left(\frac{T}{u}-\frac{T_{\pm}}{u_{\pm}}\right),\\
\theta T_x +\nu uu_x+q\rho dY_x & =& m\left((R+c_v)(T-T_{\pm})+q(Y-Y_{\pm})+\frac{1}{2}(u^2-u^2_{\pm})\right). 
\end{eqnarray*}
Next define the variable $Z$ by the relationship
$$Z=Y-\rho d\frac{Y_x}{m}=Y-d\frac{Y_x}{u},$$
and note that $Y_x$ vanishes at $\pm\infty$ to obtain
$$Z_-=Y_-,\qquad Z_+=Y_+.$$ 
The equation for $Y$ can thus be rewritten as 
$$Z_x=-k\frac{Y}{u}\varphi(T). \quad(\mbox{Note}\; u\ne 0)$$ 
Finally, rescaling to make the equations dimensionless, one arrives at the system 
\begin{eqnarray}
\nu u_x & = & u-1+\frac{1}{\gamma M^2}\left(\frac{T}{u}-1\right), \label{singmom}\\
\theta T_x & = & T-1 -\frac{\gamma -1}{\gamma}(T-u)+qZ-\frac{(\gamma-1)M^2}{2}(u-1)^2, \label{singeng} \\
dY_x & = & u(Y-Z), \label{singprog}\\
Z_x & = & -\frac{Y}{u}\varphi(T). \label{singextra} 
\end{eqnarray}
Here all quantities have been rescaled, $M$, defined by $M^2=u^2/(\gamma RT)$, is the Mach number. 

\begin{figure}[ht]
\centerline{\epsfxsize12truecm\epsfbox{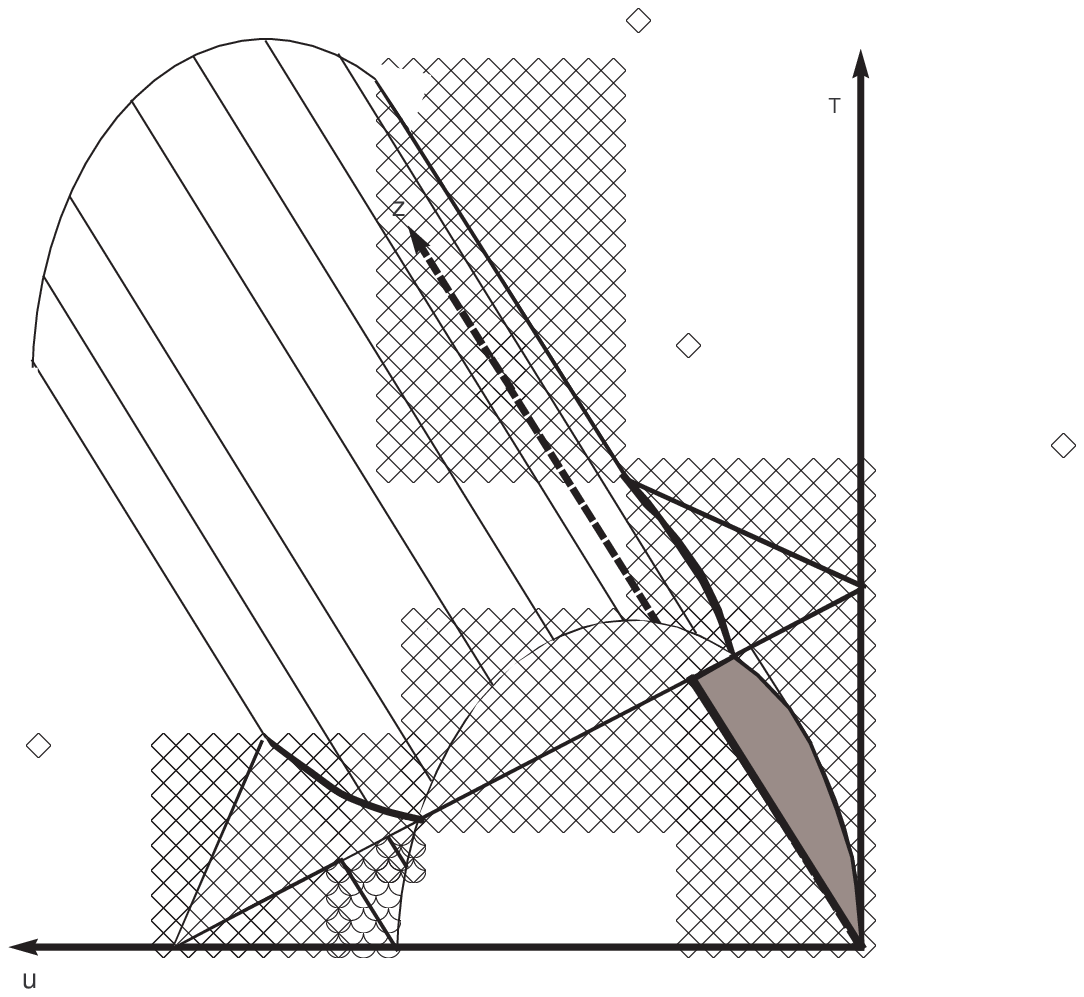}}
\caption{The Intersection of $\mathcal{K}$ and the Trough}
\label{manifold}
\end{figure}

\begin{figure}[ht]
\centerline{\epsfxsize10truecm\epsfbox{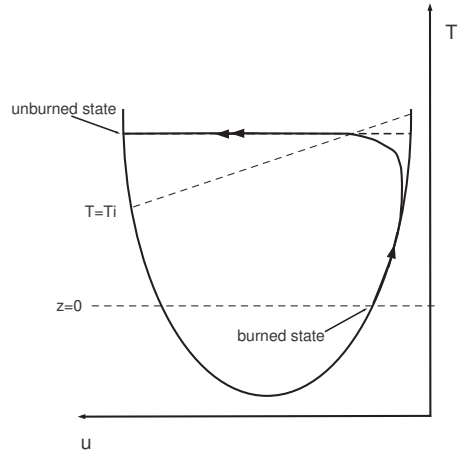}}
\caption{Detonation}
\label{detonation}
\end{figure}

The values of the dissipative coefficients ($\nu, \theta, \mbox{and}\: d$) are typically quite small. Taking advantage of this smallness, the next step is to fix small values $\hat{\nu},\hat{\theta},$ and $\hat{d}$, and then to set $\nu=\epsilon \hat{\nu}$, $\theta = \epsilon\hat{\theta}$ and $d=\epsilon\hat{d}$ so that the system (\ref{singmom})-(\ref{singextra}) takes the form  
\begin{eqnarray}
\epsilon\hat{\nu} u_x & = & u-1+\frac{1}{\gamma M^2}\left(\frac{T}{u}-1\right), \label{epssingmom}\\
\epsilon\hat{\theta} T_x & = & T-1 -\frac{\gamma -1}{\gamma}(T-u)+qZ-\frac{(\gamma-1)M^2}{2}(u-1)^2, \label{epssingeng} \\
\epsilon\hat{d}Y_x & = & u(Y-Z), \label{epssingprog}\\
Z_x & = & -\frac{Y}{u}\varphi(T). \label{epssingextra} 
\end{eqnarray}
Here $\epsilon$ is supposed to be small, and this system is singularly perturbed. Setting $\epsilon =0$ yields the reduced (slow flow) system
\begin{eqnarray}
0 & = & u-1+\frac{1}{\gamma M^2}\left(\frac{T}{u}-1\right), \label{redmom}\\
0 & = & T-1 -\frac{\gamma -1}{\gamma}(T-u)+qZ-\frac{(\gamma-1)M^2}{2}(u-1)^2, \label{redeng} \\
0 & = & u(Y-Z), \label{redprog}\\
Z_x & = & -\frac{Y}{u}\varphi(T). \label{redextra} 
\end{eqnarray}
Equations (\ref{redmom})-(\ref{redprog}) define a one-dimensional manifold $\mathcal{C}$ upon which equation (\ref{redextra}) describes a flow.  Upon noting that (\ref{redmom}) is independent of $Y$ and $Z$; and (\ref{redeng}) is independent of $Y$; and when $u\ne 0$, (\ref{redprog}) implies that $Y=Z$; $\mathcal{C}$ can be visualized in three-dimensional $uTZ$-space. The equation (\ref{redmom}) describes a parabolic trough. Using (\ref{redmom}) in (\ref{redeng}) yields
$$0=T\frac{\gamma +1}{2\gamma}+qZ+u\frac{\gamma -1}{2\gamma}(1+\gamma M^2)-1-\frac{(\gamma -1)M^2}{2},$$
which describes a plane $\mathcal{K}$ in $uTZ$-space. The manifold $\mathcal{C}$ is exactly the intersection of this plane and the parabolic trough. This intersection is pictured in Figure~\ref{manifold}. Note that $\mathcal{C}$ splits into two branches, as the requirement that there be two burned end states corresponding to the fixed unburned state forces the vertex of the intersection to have a negative $Z$ coordinate. In fact the vertex is exactly at $Z=0$ in the boundary Chapman-Jouget case.
Note also that all the rest points of (\ref{epssingmom})-(\ref{epssingextra}) are contained in $\mathcal{C}$. Rescaling the independent variable by $\xi=x/\epsilon$ in (\ref{epssingmom})-(\ref{epssingextra}) yields the equivalent (when $\epsilon \neq 0$) system
\begin{eqnarray}
\hat{\nu} u_{\xi} & = & u-1+\frac{1}{\gamma M^2}\left(\frac{T}{u}-1\right), \label{fastsingmom}\\
\hat{\theta} T_{\xi} & = & T-1 -\frac{\gamma -1}{\gamma}(T-u)+qZ-\frac{(\gamma-1)M^2}{2}(u-1)^2, \label{fastsingeng} \\
\hat{d}Y_{\xi} & = & u(Y-Z), \label{fastsingprog}\\
Z_{\xi} & = & \epsilon\left(-\frac{Y}{u}\varphi(T)\right). \label{fastsingextra} 
\end{eqnarray}
Setting $\epsilon =0$ in (\ref{fastsingmom})-(\ref{fastsingextra}) yields the layer (fast flow) system
\begin{eqnarray}
\hat{\nu} u_{\xi} & = & u-1+\frac{1}{\gamma M^2}\left(\frac{T}{u}-1\right), \label{laysingmom}\\
\hat{\theta} T_{\xi} & = & T-1 -\frac{\gamma -1}{\gamma}(T-u)+qZ-\frac{(\gamma-1)M^2}{2}(u-1)^2, \label{laysingeng} \\
\hat{d}Y_{\xi} & = & u(Y-Z), \label{laysingprog}\\
Z_{\xi} & = & 0, \label{laysingextra} 
\end{eqnarray}
%
\begin{figure}[ht]
\centerline{\epsfxsize10truecm\epsfbox{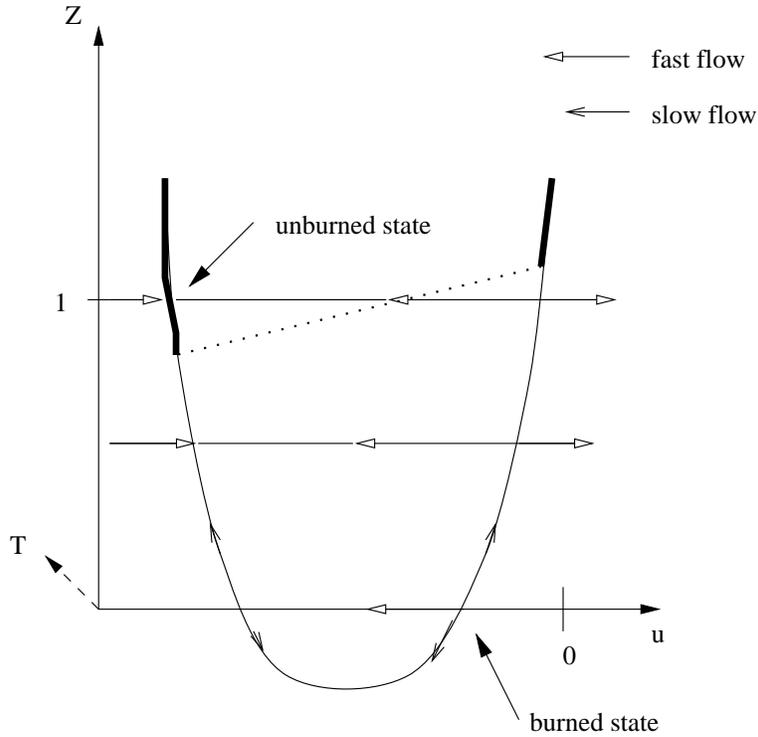}}
\caption{Singular Flow}
\label{singularflow}
\end{figure}
which is the gas dynamic shock problem of \cite{Gi} coupled to (\ref{laysingprog}). Thus for each constant $Z$ slice, the fast flow is described by the shock layer analysis of \cite{Gi}. Figure \ref{singularflow} represents the structure of the singular ($\epsilon =0$) flow looking down on $\mathcal{C}$ from a vantage point perpendicular to the plane $\mathcal{K}$. Hollow arrows represent fast flow, while single arrows represent slow flow. The fast flow in the plane $Z=\mbox{constant}$ corresponds to a nonreacting gas dynamical shock, while the slow flow proceeds on each branch of $\mathcal{C}$ representing the progress of the reaction. Combining, we see in Figure~\ref{detonation} the perturbed composite orbit of a strong detonation. Note the presence of the ZND structure, namely the gas dynamical shock to the Neumann spike which raises the temperature above ignition followed by a reaction resolving to the final totally burned state. The same structure can be seen in the analysis of the scalar Majda model \cite{M,RV,LyZ2}.
%
Having dealt with the question of existence (at least for small $\nu$, $\theta$, and $d$), we turn our attention to the tools of our stability analysis.
%
%
\subsection{The Evans Function and The Gap and Tracking Lemmas}
The proper notion of stability for traveling waves connecting constant end states, as are the detonations and shocks we consider, is that of \emph{orbital stability}, that is the convergence of the perturbed solution to the manifold of solutions which connect the same two end states. We refer to this as nonlinear stability. A weaker notion is that of \emph{linearized orbital stability} defined as convergence of the perturbation solving the linearized equations to the tangent manifold of solutions connecting the end states. Closely related to this concept of linear stability is \emph{spectral stability}; a profile is spectrally stable if the linearized operator about the wave, $L$, has no spectrum in the set $\{\lam\in \C\; |\;\re\lam\geq 0\}$ except $\lam =0$. We note that translation invariance implies that $0$ is an eigenvalue of the linearized operator, so necessarily $0\in\sigma(L)$. Spectral stability is clearly necessary for linear stability, itself necessary for nonlinear stability. Recent work by Mascia and Zumbrun \cite{MZ1,MZ2,MZ3} extending earlier work of Zumbrun and Howard \cite{ZH} indicates that spectral stability implies nonlinear orbital stability in the settings of viscous conservation laws and relaxation systems which are closely related to the combustion systems we consider. In light of these results, the determination of spectral stability can be regarded as the essential initial step in determining the stability of detonation waves. We also note that \cite{Z5,Zhandbook} contain large-amplitude nonlinear stability results, relevant to the detonation problem. 

\subsubsection{Background}
As noted above, a vital step in determining stability amounts to locating the spectrum of a linear operator. The search for spectrum is facilitated by the Evans function, $D(\lambda)$, a potent tool in the investigation of stability of traveling waves. This function, analytic on the unstable half-plane, is an infinite-dimensional analogue of the characteristic polynomial. Zeros of $D(\cdot)$ correspond to eigenvalues of the linearized operator about the wave. The Evans function was introduced in \cite{Ev1,Ev2,Ev3,Ev4} specifically to study nerve axon equations and further developed in \cite{AGJ} to the case of semilinear parabolic systems. In \cite{GZ} the use of the Evans function was extended to the case of $2\times 2$ viscous conservation laws. The extension to $n\times n$ systems has been completed in \cite{BSZ}.  

The construction of the Evans function is accomplished by analytically parametrizing the unstable/stable manifolds of the variable coefficient eigenvalue equations. This is done by comparing these objects to the corresponding unstable/stable manifolds for the limiting constant coefficient systems at $\pm \infty$. We now give an abbreviated description of the construction for the case of a system of viscous conservation laws
$$U_t+f(U)_x=(B(U)U_x)_x,$$
where $U,f\in \mathbf{R}^n$ and $B$ is an $n\times n$ matrix. As we shall be concerned with the cases of gas dynamics and combustion, it will be the case that the matrix $B$ is incompletely parabolic. However, the equations of compressible gas dynamics satisfy the symmetrizability, dissipativity, and block structure conditions of Kawashima \cite{Kaw}:  

\noindent
\emph{symmetrizability}
\begin{align} &\text{There exists a symmetrizer $A^0(U)$, symmetric and positive definite,} \nonumber\\ &\text{such that $A^0(U)A(U)$ is symmetric and $A^0(U)B(U)$}\tag{+}\label{symmetrizability} \\ &\text{is symmetric and positive semidefinite.}\nonumber
\end{align}
\emph{dissipativity}
\begin{equation}
\label{dissipativity}\text{There is no eigenvector of $A(U)$ lying in the kernel of $B(U)$.} \tag{++}
\end{equation}
\emph{block structure}  
\begin{equation}
\label{blockstructure}\text{The right kernel of $B(U)$ is independent of $U$.} \tag{+++}
\end{equation}
In the above $A(U)$ denotes the Jacobian matrix of the flux $f$.
We discuss the construction in the case of combustion in more detail below. A viscous profile,
$$U(x,t)=\bar{U}(x-st),\qquad \bar{U}(\pm\infty)=U_{\pm},$$
is a solution of the (integrated) traveling wave ODE
$$-s(U-U_-)+f(U)-f(U_-)=(B(U)U')'.$$
Taking, without loss of generality, $s=0$ and linearizing about $\bar{U}(x)$, one obtains an equation modeling the approximate evolution of a small disturbance, $v$. This equation has the form 
$$v_t=Lv:=(Bv_x)_x-(Av)_x.$$  
Here $B(x)=B(\bar{U}(x))$ and $A(x)$ is determined by the relationship $Av=f'(\bar{U}(x))v-B'(\bar{U}(x))v\bar{U}_x$. The operator $L$ is the linearized operator about the wave $\bar{U}$. 
\begin{defn}The profile $\bar{U}(\cdot)$ is \emph{spectrally stable} if the linearized operator $L$ about the wave has no spectrum in the closed unstable complex half plane $\{\lam \in\C \;|\; \re\lam\geq 0\}$ except at $\lam =0$.
\end{defn} 
The next lemma allows us to narrow our search for spectrum.
\begin{lemma} \label{noess} Assuming \eqref{symmetrizability}-\eqref{blockstructure}, the operator $L$ has no essential spectrum in $\{\lam\in\C\;|\; \re\lam\ge 0\}\backslash 0$.
\end{lemma}
The lemma follows by a standard argument of \cite{He} provided that the constant solutions $U\equiv U_{\pm}$ are linearly stable. Such stability follows at once from the condition
\begin{equation*}
\re \sigma (i\xi A(U)-|\xi|^2B(U))\le \frac{-\theta |\xi|^2}{1+|\xi|^2},\quad \theta >0,  
\end{equation*}
which is equivalent to \eqref{dissipativity} in the presence of \eqref{symmetrizability} by an argument of \cite{SK}. Details can be found in any of \cite{GZ,ZH,ZS,Z1,Z4}. Thus for the systems of our interest $\sigma_{ess}(L)$ is confined to the left complex half plane except the origin, thus the only possible unstable spectrum consists of isolated eigenvalues of finite multiplicity.  Determination of spectral stability is then reduced to checking that the operator $L$ has no unstable point spectrum. The focus then is on the eigenvalue equation for this operator,
$$Lw=\lam w,$$
and solutions $w\in L^2$ with corresponding eigenvalue $\lam$ such that $\re\lam >0$. 

The eigenvalue equation can be recast as a system of first order ODE
\begin{equation}W'=\mathbb{A}(x,\lam)W,\quad W\in\C^N .\label{firstODE} \end{equation}
Because the wave $\bar{U}$ connects constant states $U_{\pm}$, the matrix $\mathbb{A}$ has limits as $x\rightarrow \pm\infty$. Thus
$$\mathbb{A}(x,\lam)\longrightarrow \mathbb{A}_{\pm}(\lam)\;\;\;\mbox{as}\;\;x\rightarrow \pm\infty .$$
The idea then, is to connect the stable (resp. unstable) manifolds of the system $W'=\mathbb{A}(x,\lam)W$ at $\pm\infty$ to the stable (resp. unstable) subspaces $\mathcal{S}^{\pm}(\lam)$(resp. $\mathcal{U}^{\pm}(\lam)$) of the constant-coefficient systems at each of $\pm\infty$.
%
The procedure we are outlining requires that the system have \emph{consistent splitting of the asymptotic systems}.
\begin{defn} We say that the system (\ref{firstODE}) has \emph{consistent splitting on $\Omega\subset \C$} if the matrices $\mathbb{A}_{\pm}(\lam)$ are both hyperbolic for all \lam\ in a region $\Omega\subset\C$, and there is an integer $k$ such that the stable (resp. unstable) subspaces of $\mathbb{A}_+$ and $\mathbb{A}_-$ are both $k$-dimensional (resp. $(N-k)$-dimensional). \end{defn}  
The set $\Omega$ is called the \emph{region of consistent splitting}. In \cite{GZ} it is shown that the linear stability of the constant solutions is equivalent to: 
\begin{align} &\text{Eigenvalues $\mu^{\pm}(\lam)$ of $\mathbb{A}_{\pm}(\lam)$ have nonvanishing real part}\nonumber\\
	      &\text{for all \lam\ with $\re\lam >0$}. \tag{**} \label{consplit}
\end{align}
It follows that $\Omega$, the region of consistent splitting, contains at least the unstable complex half plane. One important consequence  of \eqref{consplit} is that the the number of stable/unstable eigenvalues of $\mathbb{A}_{\pm}(\lam)$ can be counted as $\lam \to +\infty$ along the real axis.  
Then given bases $\{\phi_1^+,\ldots\phi_k^+\}$ and $\{ \phi_{k+1}^-,\ldots\phi_N^-\}$ of the stable manifold at $+\infty$ and the unstable manifold at $-\infty$, the idea is to define the Evans function as
$$\evans=\det(\phi_1^+,\ldots\phi_k^+,\phi_{k+1}^-,\ldots ,\phi_N^-)|_{x=0}.$$
A natural way to attempt such a procedure is to choose as bases for $\mathcal{S}^{\pm}(\lam)$ and $\mathcal{U}^{\pm}(\lam)$ the purely exponential normal modes of $W'=\mathbb{A}_{\pm}(\lam)$. Unfortunately, it is not possible to make this choice analytically with respect to $\lam$ as some eigenvalues of $\mathbb{A}_{\pm}$ may coalesce as $\lam$ varies. The utility and power of the Evans function comes from the fact that this difficulty can be surmounted, and $D(\lam)$ can be chosen to be analytic. 

The solution, due to an elegant construction of \cite{AGJ}, is to track volume forms rather than individual solutions. We associate to any collection $V_1,V_2,\ldots ,V_n$ of vectors the wedge product
$$V_1\wedge V_2\wedge\cdots \wedge V_n.$$
This determines an embedding of the manifold of $n$-dimensional bases into the manifold of $n$-forms. More precisely this determines an embedding into the submanifold of $n$-forms expressible as a single product, the \emph{pure} $n$-forms. The benefit of this approach can be seen by taking a set of $n$ solutions $W_1,W_2,\ldots ,W_n$ of the eigenvalue equation $W'=\mathbb{A}(x,\lam)W$ and noticing that the corresponding $n$-form $\zeta = W_1\wedge W_2\wedge \cdots \wedge W_n$ solves the ``lifted'' linear ODE
$$\zeta'=\mathcal{A}(x,\lam)\zeta,$$
where the operator $\mathcal{A}$ is determined by
$$\mathcal{A}=\mathbb{A}W_1\wedge W_2\wedge \cdots \wedge W_n+\cdots +W_1\wedge W_2\wedge\cdots \wedge\mathbb{A}W_n.$$
If the collection $W_1,W_2,\ldots ,W_n$ consists of eigenvectors of $\mathbb{A}_{\pm}$ with corresponding eigenvalues $\mu_1,\mu_2,\ldots ,\mu_n$, then it is immediate that the wedge $W_1\wedge W_2\wedge \cdots \wedge W_n$ is an eigenvector of $\mathcal{A}$ with corresponding \emph{simple} eigenvalue $\mu_1+\mu_2+\ldots +\mu_n$. Thus the volume form associated with any basis  of $\mathcal{S}^{\pm}(\lam)$ or $\mathcal{U}^{\pm}(\lam)$ is a simple eigenvector of $\mathcal{A}$ corresponding to a purely exponential growth or decay mode. In this lifted setting, the eigenvectors are simple, thus they depend on $\lam$ in an analytic fashion. The construction demonstrates the important fact that the eigenspaces vary analytically with respect to a parameter even when the individual eigenvectors do not.

For notational convenience, we follow the standard convention of associating the full $N$-volume forms with the complex numbers via the coordinate representation in the standard basis. That is, we write
$$V_1\wedge\cdots\wedge V_N=\det(V_1,\ldots,V_N).$$

\subsubsection{The Gap Lemma} 

The gap lemma of \cite{GZ} and \cite{KS} is the key technical result that allows Evans function techniques for the stability analysis of traveling waves to be extended to the case of viscous conservation laws. This lemma extends the analytic framework of \cite{AGJ} to cases in which the essential spectrum of the linearized operator touches the imaginary axis, and thus there is no spectral gap between the essential spectrum and the unstable half plane $\{\lambda \in \mathbf{C}\;|\;\mbox{Re}\lambda >0\}$. In the presence of such a gap, a standard argument of \cite{Co} provides a relationship between the behavior of solutions near $\pm \infty$ of a system of asymptotically constant-coefficient eigenvalue ODEs and the corresponding solutions of the limiting, constant-coefficient equations.

More precisely consider an ODE with parameter (as obtained above by rewriting the eigenvalue equation as a first-order system) 
\begin{equation} W'=\mathbb{A}(x,\lam)W, \label{EvalODE} \end{equation}
where the differentiation is with respect to $x$, and $\mathbb{A}$ is continuous in $x$ and analytic with respect to $\lambda$. Moreover, suppose also that $\mathbb{A} \rightarrow \mathbb{A}_{\pm}$ as $x \rightarrow \pm \infty$. Provided that 
\begin{equation} \int_0^{\pm \infty} |\mathbb{A}-\mathbb{A}_{\pm}|dx<+\infty ,\label{finiteint} 
\end{equation}
then there is a one-to-one correspondence between the normal modes $V_j^{\pm}e^{\mu_j^{\pm}x}$ of the constant coefficient limiting system
$$W'=\mathbb{A}_{\pm}(\lambda)W,$$
where $\mu_j^{\pm}(\lambda),V_j^{\pm}(\lambda )$ is an eigenvalue, eigenvector pair corresponding to $\mathbb{A}_{\pm}(\lambda )$, and solutions $W_j^{\pm}$ of (\ref{EvalODE}) having the same limiting behavior. That is 
$$W_j^{\pm}(\lambda ,x)= V_j^{\pm}e^{\mu_j^{\pm}x}(1+o(1))\;\; \mbox{as}\;\; x\rightarrow \pm \infty .$$
The argument in \cite{Co} uses a fixed point iteration scheme depending on the sign of differences of the real parts of the eigenvalues $\mu_j$.  In the case of strict separation of the eigenvalues, a \emph{spectral gap}, the fixed point is the uniform limit of an analytic sequence of iterates, and thus analyticity in $\lambda$ is preserved. In our case of interest there is no spectral gap; the above argument breaks down. The key observation in \cite{GZ} was that in the absence of such a gap, analyticity can be preserved provided that in lieu of (\ref{finiteint}) the stronger hypothesis
\begin{equation} |\mathbb{A}-\mathbb{A}_{\pm}|=\mathit{O}(e^{-\alpha |x|})\;\; \mbox{as}\;\; x\rightarrow \pm\infty, \label{expapproach} 
\end{equation} 
is made.

\begin{thm}[The Gap Lemma]Let $\mathbb{A}(x,\lam)$ be continuous in $x$ and analytic in $\lam$ with 
$$\mathbb{A}(x,\lam)\rightarrow\mathbb{A}_{\pm}(\lam)\;\; \mbox{as}\;\;x\rightarrow\pm\infty,$$ 
at an exponential rate $e^{-\alpha |x|},\alpha>0$, and let $\zeta^-(\lam)$ and $\eta^-(\lam)$ be analytic $n$- and $(n-k)$-forms associated to the complementary invariant subspaces of $\mathbb{A}_-(\lam)$, $C^-$ and $E^-$ with spectral gap $\beta$. Furthermore put $\tau_{C^-}$ equal to the trace of $\mathcal{A}$ restricted to $C^-$. Then there exists a solution $\mathcal{W}(x,\lam)$ of the lifted ODE $\zeta'=\mathcal{A}(x,\lam)\zeta$ of the form
$$\mathcal{W}(x,\lam)=\zeta(x,\lambda)e^{\tau_{C^-}},$$
where $\zeta$ (and thus $\mathcal{W}$) is $C^1$ in $x$ and locally analytic in $\lam$. Moreover $\eta(x,\lam)$ satisfies
$$(\frac{\partial}{\partial\lam})^j\zeta(x,\lam)=(\frac{\partial}{\partial\lam})^j\zeta^-(\lam)+\mathit{O}(e^{-\bar{\alpha}|x|}|\zeta^-(\lam)|),$$
when $x<0$ for $j=0,1,\ldots$
\label{gaplemma} \end{thm}
\noindent
See \cite{Z1,Z4,ZH,GZ} and \cite{KS} for further discussion and a proof.
Appealing to the Gap Lemma, we thus obtain bases $\{\phi_1^+(x,\lam),\ldots\phi_k^+(x,\lam)\}$ of the stable manifold at $+\infty$ and $\{\phi_{k+1}^-(x,\lam),\ldots\phi_N^-(x,\lam)\}$ of the unstable manifold at $-\infty$. Therefore we can indeed define $D(\lam)$ by
$$\evans=\det(\phi_1^+,\ldots\phi_k^+,\phi_{k+1}^-,\ldots ,\phi_N^-)|_{x=0}.$$
An important feature of the construction is that \evans\ can be chosen to be real-valued for real \lam\ .
\begin{thm}There exist bases $\phi_j^{\pm}$ such that the \evans\ satisfies
$$D(\bar{\lam})=\overline{D(\lam)}.$$
In particular \evans\ is real valued for $\lam \in \R$.
\end{thm}
\noindent
See \cite{Z1,Z4} for details. A proof involves tracing through the various steps in the construction of the Evans function, and verifying that complex symmetry is preserved at each stage.
%
%
%
\subsubsection{The Tracking Lemma}
In the calculation of the stability conditions (described below) it will be necessary to connect information about the the sign of the Evans function $D(\lambda)$ as $\lambda \longrightarrow \infty$ along the real axis to the normalizations for the bases of stable and unstable manifolds chosen at $\lambda =0$. This can be accomplished by using the tracking lemma. See \cite{GZ,Z1,Z4}. 
\begin{thm}[The Tracking Lemma]
For $\delta$ sufficiently small, solutions $w^+/w^-$ of
\begin{equation} w'=(\mathbb{A}_0(x,\delta)+\Theta(x,\delta))w,\;\;w\in\mathbf{C}^N, \label{introtrack}\end{equation}
where $\delta \rightarrow 0$ is a small parameter and 
$$|\mathbb{A}_0'|+|\Theta |\leq C\delta ,\;\; |\mathbb{A}_0|<C.$$
decaying at $+\infty /-\infty$ at rate $e^{\underline{\tilde{\tilde{\alpha}}}x}/e^{\tilde{\tilde{\overline{\alpha}}}x}$ lie always within the cones
\begin{equation}
\mathbb{K}_-=\left\{w \left|\;\;\frac{|P_{\mathcal{S}}w|}{|P_{\mathcal{U}}w|}\leq\frac{C\delta}{\eta}\right.\right\}, 
\end{equation}
and
\begin{equation}
\mathbb{K}_+=\left\{w \left|\;\;\frac{|P_{\mathcal{U}}w|}{|P_{\mathcal{S}}w|}\leq\frac{C\delta}{\eta}\right.\right\},
\end{equation}
respectively, for any $\underline{\tilde{\tilde{\alpha}}}<\liminf_{x\rightarrow +\infty}\overline{\alpha}$, $\tilde{\tilde{\overline{\alpha}}}>\limsup_{x\rightarrow -\infty}\underline{\alpha}$. Moreover there hold the following uniform growth/decay rates:
\begin{eqnarray}
\frac{|w^+(x)|}{|w^+(y)|} & \leq & Ce^{\tilde{\underline{\alpha}}|x-y|}, \\
\frac{|w^-(x)|}{|w^-(y)|} & \geq & C^{-1}e^{\tilde{\overline{\alpha}}|x-y|}, 
\end{eqnarray}
for all $x>y$, and symmetrically for $x<y$, for any $\tilde{\underline{\alpha}}>\max_x\underline{\alpha}$, $\tilde{\overline{\alpha}}<\min_x\overline{\alpha}$. (Note: Here $C$ depends in part upon the choice of $\tilde{\overline{\alpha}}$ and $\tilde{\underline{\alpha}}$.
\end{thm}
\noindent
We note that $P_{\mathcal{S}}$ and $P_{\mathcal{U}}$ are projections onto eigenspaces corresponding to two different spectral groups. Rescaling in \lam\ transforms $W'=\mathbb{A}(x,\lam)W$ to a system of the form \eqref{introtrack}.

\subsection{Discussion}
\subsubsection{Stability Conditions}
Even though the Evans function is not typically evaluable, it is possible to obtain information about its zeros in the following way. Due to a translational eigenvalue at $\lambda =0$, $D(0)=0$. One then calculates $\sgn D'(0)$ and the sign of $D(\cdot)$ as $\lambda \longrightarrow \infty$ along the real axis. (Recall: $D(\lambda)$ can be chosen to be real for real $\lambda$.) Combining this information yields a parity for the number of unstable zeros of $D(\cdot)$, hence unstable eigenvalues for the linearized operator. We call the quantity $\sgn D'(0)D(+\infty)$ the \emph{stability index}. When the signs agree, there must be an even number (possibly $0$) of real unstable eigenvalues, and when they disagree, there must be an odd number of such eigenvalues. Recall that complex eigenvalues occur in conjugate pairs, hence they do not affect the parity.
Clearly then
$$\sgn D'(0)D(+\infty)\geq 0,$$
is necessary for spectral stability. On the other hand when the index is negative, a positive growth rate is detected, and the wave under consideration is determined to be unstable. Thus the stability index is best suited as a predictor of instability. We remark that as the stability index only determines the parity of unstable eigenvalues, the condition $\sgn D'(0)D(+\infty)\geq 0$ is not sufficient on its own to conclude spectral stability. The index yields only incomplete stability information; the possibilities of complex conjugate unstable eigenvalues and/or even numbers of unstable real eigenvalues are not detected by this approach. Nonetheless, the stability index serves as a useful starting place in stability investigation.  
\subsubsection{Results}
Here we describe the two main results.
\begin{thm} The stability index for a strong detonation solution of Equations~(\ref{eq:mass})-(\ref{eq:progress}) with Lax 3-shock structure has the form
$$\tilde{\Gamma}=\sgn D'(0)D(+\infty)=\sgn\bar{\gamma}\Delta,$$
where $\bar{\gamma}$ is a constant measuring transversality of the stable/unstable manifolds of the traveling wave ODE and
\begin{equation}\Delta =\det (r_1^-,r_2^-,[U]+\vec{q}). \label{detdelta} \end{equation}
Moreover for an ideal gas, the sign of the stability index is consistent with spectral stability in the ZND limit.
\label{result1}\end{thm}
\begin{thm}Strong detonations are spectrally stable for sufficiently small $q$ provided the underlying gas-dynamical shock (of arbitrary strength) is stable.
\label{result2} \end{thm}
In these theorems $r_j^{\pm}$ are right eigenvectors of the flux Jacobian, $[U]$ is a vector of jumps in the gas-dynamical conserved quantities $\rho,m,\E$ (density, momentum, total energy), and $\vec{q}$ is the vector 
$$\left(\begin{array}{c} 0 \\ 0 \\ q\end{array}\right), $$
where $q>0$ represents the energy liberated during the exothermic chemical reaction. We note that the term $\Delta$ in (\ref{detdelta}), which appears in the stability index due to the low frequency calculation of $D'(0)$, is related the Lopatinski determinant itself a ``stability function'' for inviscid shocks. See \cite{JL} and the references therein. In Section 3 we detail the reduction of the equation $\Delta =0$ to   
\begin{equation*}
M^2[1/\rho]p_e-M-1=0,
\end{equation*} 
the well-known instability condition of Majda \cite{MBOOK} for inviscid shocks. Note that this is independent of $q$. We also remark that the finding $\tilde{\Gamma}>0$ for an ideal gas in the ZND limit has the implication that instability, if it occurs must be of ``galloping'' type, i.e. corresponding to the crossing of a complex conjugate pair of eigenvalues into the right half plane. This is consistent with both laboratory and numerical experiments.

The second theorem is related to results of Liu and Ying \cite{LYi} and Li, Liu, and Tan \cite{LLT} for various versions of the Majda model. There, full nonlinear stability is established for strong detonations in the Majda model when $q$ is sufficiently small. In \cite{LLT}, the authors prove nonlinear stability for strong detonations in a version of the Majda model with species diffusion using techniques in the spirit of the Evans function. Finally we remark again, that by the program of \cite{ZH,MZ1,MZ2}, it is expected that spectral stability should be equivalent to nonlinear stability.

\subsubsection{Extensions}
For the calculations in this paper, we have made the simplifying assumption that the equations of state are independent of the progress of the reaction. Though standard in the literature, this is clearly an idealization as the nature of the gas changes during the chemical reaction as pointed out in \cite{CHT}. One extension is to carry out the analysis in the more realistic setting of reaction-dependent equations of state as discussed in the context of the Majda model in \cite{LyZ2}. Also, we note that the $q\to 0$ argument in Section~\ref{smallq} fails since the $q=0$ gas equation is still coupled to the reaction equation through the equation of state.

Another interesting direction of future study is a more detailed examination of the effect (if any) of multiple reactants the stability index and its sign. In particular, the analysis of \cite{GS2} provides the geometric information required to evaluate (in the ZND limit) the transversality coefficient in the stability index for the interesting two-species reactions they consider. In particular, while an exothermic-exothermic two-step reaction behaves much as the one-step exothermic reaction we consider, an exothermic-endothermic two-step reaction has a richer structure \cite{FD}.

%
%
\section{Nonreacting Gas}\label{NSsection}
In this section we consider gas dynamics as modeled by the Navier-Stokes equations in one space dimension; our main focus is the calculation of the stability index. These computations will prove useful when we shift our focus to detonations in the next section. The system takes the form
\begin{align}
\rho_t& + (\rho u)_x  =  0, \label{eq:massf} \\
(\rho u)_t& + (\rho u^2+p)_x  =  (\nu u_x)_x, \label{eq:momentumf} \\
\E_t&+(u\E+up)_x  =  (\theta T_x)_x+(\nu uu_x)_x.\label{eq:energyf}
\end{align}
The system \eqref{eq:massf}-\eqref{eq:energyf} features five unknowns ($\rho,u,e,p,T)$ and three equations. The system is completed by equations of state which incorporate the physical properties of the particular gas being modeled. We thus obtain a complete description of the fluid flow by assuming that $p$ and $T$ are given functions of $\rho$ and $e$. 
The sound speed is 
$$c=\sqrt{p_{\rho}+\rho^{-2}pp_e}.$$
The three differential equations (\ref{eq:massf})-(\ref{eq:energyf}) and the equations of state give a set of five equations for the five variables $\rho$, $u$, $e$, $T$ and $p$.

For some portions of the analysis, we will further assume that the gas under consideration is ideal and polytropic so that the specific forms of the equations of state are
$$p(\rho,e)=\Gamma\rho e,\quad T(\rho,e)=c_v^{-1}e,$$
where the constants $c_v$ and $\Gamma$ are as in the previous section. We note that in this case the sound speed satisfies 
\begin{equation}
c^2=\Gamma e+\Gamma^2e=(1+\Gamma)\Gamma e =\gamma \Gamma e, \label{soundspeed}
\end{equation}
where $\gamma=1+\Gamma$. 


We rewrite the system (\ref{eq:massf})-(\ref{eq:energyf}) in terms of the conserved quantities $\rho$, $m$, and $\E$,
\begin{eqnarray}
\rho_t + m_x & = & 0, \label{massf} \\
m_t + \left(\frac{m^2}{\rho}+p\right)_x & = & \left(\nu \left(\frac{m}{\rho}\right)_x\right)_x, \label{momentumf} \\
\mathcal{E}_t+\left(\frac{m}{\rho}(\mathcal{E}+p)\right)_x & = & \left((\theta T)_x+\left(\nu \left(\frac{m}{\rho}\right)\left(\frac{m}{\rho}\right)_x\right)\right)_x.\label{energyf}
\end{eqnarray}
Rewritten once more in the form of a viscous conservation law $U_t+f(U)_x=(B(U)U_x)$, the form of the viscosity matrix $B$ becomes apparent
\begin{multline*}
\left(\begin{array}{c} \rho_t \\ m_t \\ \mathcal{E}_t \end{array}\right)+\partial_x\left(\begin{array}{c} m \\ m^2/\rho +p\\ m\rho^{-1}(\mathcal{E}+p) \end{array}\right)  = \\ \partial_x\left[\begin{pmatrix} 0 & 0 & 0 \\ -\frac{\nu m}{\rho^2} & \frac{\nu}{\rho} & 0 \\ \theta T_{\rho}+T_ee_{\rho}-\frac{\nu m^2}{\rho^3} & \theta T_ee_m+\frac{\nu m}{\rho^2} & \theta T_ee_{\E} \end{pmatrix}\left(\begin{array}{c} \rho_x \\ m_x \\ \mathcal{E}_x \end{array}\right)\right].
\end{multline*}
The Jacobian matrix $A(U)$ of the flux $f$ has the form
$$\left(\begin{array}{ccc} 0 & 1 & 0 \\ -u^2+p_{\rho}+p_ee_{\rho} & 2u-u\frac{p_{e}}{\rho} & \frac{p_{e}}{\rho} \\ u(-E-\frac{p}{\rho}+p_{\rho}+p_ee_{\rho}) & E+\frac{p}{\rho}-u^2\frac{p_e}{\rho} & u+u\frac{p_e}{\rho}
\end{array}\right).$$
To calculate the eigenvalues/eigenvectors of $A$, we use the device (see \cite{JL} and the isentropic gas section of \cite{Z2}) of conjugating by appropriate ``shift'' matrices so that the conjugated matrix has a particularly simple form. Following this procedure we obtain eigenvalues
\begin{eqnarray*}
a_1 & = & u-c, \\
a_2 & = & u, \\
a_3 & = & u+c, 
\end{eqnarray*} 
and right eigenvectors
\begin{eqnarray*}
r_1 & = & \left(\begin{array}{c}1 \\ u-c \\ \frac{u^2}{2}-cu+\frac{p}{\rho}+ e \end{array}\right), \\
r_2 & = & \left(\begin{array}{c} 1 \\ u \\ \frac{u^2}{2} \end{array}\right), \\
r_3 & = & \left(\begin{array}{c}1 \\ u+c \\ \frac{u^2}{2}+cu+\frac{p}{\rho}+ e \end{array}\right). 
\end{eqnarray*}
The left eigenvectors are 
\begin{eqnarray*}
l_1 & = &\left(p_{\rho}-\frac{p_e}{\rho}e+cu+\frac{p_eu^2}{2\rho},-c-\frac{p_e}{\rho}u,\frac{p_e}{\rho}\right), \\
l_2 & = &\left(-e-\frac{p}{\rho}+\frac{u^2}{2},-u,1\right), \\
l_3 & = &\left( p_{\rho}-\frac{p_e}{\rho}e-cu+\frac{p_eu^2}{2\rho},c-\frac{p_e}{\rho}u,\frac{p_e}{\rho}\right).
\end{eqnarray*}
We shall restrict our attention to a Lax 3-shock. That is, we suppose that the shock speed $s$ satisfies the inequalities
\begin{equation}
a_2^-<s<a_3^-,\qquad a_3^+<s \label{3shock}. 
\end{equation}
The calculations for a 1-shock follow in a similar fashion.
%
%
%
%
\subsection{Traveling Wave ODE and Linearized Equations}

The traveling wave ODE is 
\begin{eqnarray}
m' & = & 0, \label{nsemassode}  \\
\left(\frac{m^2}{\rho}+p\right)' & = & \left(\nu \left(\frac{m}{\rho}\right)'\right)', \label{nsemomentumode} \\
\left(\frac{m}{\rho}(\mathcal{E}+p)\right)' & = & \left((\theta T)'+\left(\nu \left(\frac{m}{\rho}\right)\left(\frac{m}{\rho}\right)'\right)\right)', \label{nseenergyode}
\end{eqnarray}
where without loss of generality we've taken $s=0$.
Each of the equations (\ref{nsemassode})-(\ref{nseenergyode}) may be integrated up once.
\begin{eqnarray}
m-m_- & = & 0, \label{nsemassintode}  \\
\left(\frac{m^2}{\rho}+p\right)-\left(\frac{m^2}{\rho}+p\right)_- & = & \nu \left(\frac{m}{\rho}\right)', \label{nsemomentumintode} \\
\left(\frac{m}{\rho}(\mathcal{E}+p)\right)- \left(\frac{m}{\rho}(\mathcal{E}+p)\right)_-& = & (\theta T)'+\nu \left(\frac{m}{\rho}\right)\left(\frac{m}{\rho}\right)'. \label{nseenergyintode}
\end{eqnarray}
The requirement for a connection that both endstates be rest points of the ODE, leads from (\ref{nsemassintode})-(\ref{nseenergyintode}) to the Rankine-Hugoniot conditions
\begin{eqnarray*}
[m] & = & 0, \\
m^2[1/\rho] & = & -[p], \\
\left[\frac{\mathcal{E}}{\rho}\right] & = & -\left[\frac{p}{\rho}\right]. 
\end{eqnarray*}
We suppose that $\bar{U}(x)=(\bar{\rho}(x),\bar{m}(x),\bar{\mathcal{E}}(x))^{tr}$ is a stationary profile connecting endstates $U_{\pm}=(\rho_{\pm},m_{\pm},\mathcal{E}_{\pm})^{tr}$ which satisfy the Rankine-Hugoniot conditions.  We note that such profiles, if they exist, are transverse as an appropriate Wronskian is nonvanishing. Indeed global existence for connecting profiles has been shown in \cite{Gi} for equations of state which satisfy the thermodynamic conditions of \cite{Weyl}. More precisely, existence is shown when the equation of state is assumed to be such that the isentropes are convex and do not cross in the pressure-volume plane. We note that a polytropic, ideal gas satisfies this condition. 

Then linearizing about this profile, we find equations for the evolution of small perturbations $(\rho, m,\mathcal{E})$.  These equations can be written in the general form
\begin{equation*} w_t+(Aw)_x=(Bw_x)_x, \end{equation*}
where $w=(\rho, m,\mathcal{E})^{tr}$ and the matrices $A$ and $B$ depend only on $x$. More precisely we find
\begin{align}
\rho_t&+m_x  =  0, \label{NSElinearized1}\\
m_t& +(\alpha_{21}(x)\rho+\alpha_{22}(x)m+\alpha_{23}(x)\E)_x = (b_{21}(x)\rho_x+b_{22}(x)m_x)_x, \label{NSElinearized2}
\end{align}
\begin{multline}
\E_t+(\alpha_{31}(x)\rho+\alpha_{32}(x)m+\alpha_{33}(x)\E)_x = \\  (b_{31}(x)\rho_x+b_{32}(x)m_x+b_{33}(x)\E_x)_x, \label{NSElinearized3}
\end{multline}
where the coefficient functions $\alpha_{ij}(x)$ can be expressed in terms of the entries of the flux Jacobian and derivatives of the entries of the viscosity matrix. The terms $b_{ij}(x)$ correspond to the $ij$-entries of the viscosity matrix $B$. In both cases $x$-dependence arises from evaluation along the known profile $\bar{U}(x)$. 

As our interest is in spectral stability of the profile, we focus on the eigenvalue equations corresponding to (\ref{NSElinearized1})-(\ref{NSElinearized3}). They are 
\begin{align}
\lam\rho&+m'  =  0, \\
\lam m& +(\alpha_{21}(x)\rho+\alpha_{22}(x)m+\alpha_{23}(x)\E)'  = (b_{21}(x)\rho'+b_{22}(x)m')', 
\end{align}
\begin{multline}
\lam \E+(\alpha_{31}(x)\rho+\alpha_{32}(x)m+\alpha_{33}(x)\E)'  = \\ (b_{31}(x)\rho'+b_{32}(x)m'+b_{33}(x)\E')',
\end{multline}
and the corresponding limiting system as $x\rightarrow \pm\infty$ takes the form
\begin{eqnarray*}
\lam\rho+m' & = & 0, \\
\lam m +(\alpha_{21}^{\pm}\rho+\alpha_{22}^{\pm}m+\alpha_{23}^{\pm}\E)' & = & (b_{21}^{\pm}\rho'+b_{22}^{\pm}m')', \\
\lam \E+(\alpha_{31}^{\pm}\rho+\alpha_{32}^{\pm}m+\alpha_{3 3}^{\pm}\E)' & = & (b_{31}^{\pm}\rho'+b_{32}^{\pm}m'+b_{33}^{\pm}\E')'.
\end{eqnarray*}
We make the invertible change of variables as in \cite{Z4}
\begin{equation*}
\left(\begin{array}{c} z_1 \\ z_2 \\ z_3 \end{array}\right)=\underbrace{\left(\begin{array}{ccc} 0 & 1 & 0 \\ b_{21} & b_{22} & 0 \\ b_{31} & b_{32} & b_{33} \end{array}\right)}_C\left(\begin{array}{c} \rho \\ m \\ \E \end{array}\right).
\end{equation*}
Thus 
\begin{equation*}
\left\{\begin{array}{l} z_1=m, \\ z_2 = b_{21}\rho +b_{22}m, \\ z_3=b_{31}\rho+ b_{32}m+b_{33}\E, \end{array}\right.
\end{equation*}
and
\begin{equation*}
C^{-1}=\left(\begin{array}{ccc}-b_{21}^{-1}b_{22} & b_{21}^{-1} & 0 \\ 1 & 0 & 0 \\ b_{21}^{-1}b_{33}^{-1}b_{31}b_{22}-b_{33}^{-1}b_{32} & -b_{33}^{-1}b_{21}^{-1}b_{31} & b_{33}^{-1} \end{array}\right).
\end{equation*}
We can rewrite the eigenvalue equation in $z$-coordinates as
\begin{equation*}
(B(C^{-1}z)')'=(AC^{-1}z)'+\lam C^{-1}z,
\end{equation*}
or more explicitly as
\begin{eqnarray}
0 & = & z_1'+\lam(-b_{21}^{-1}b_{22}z_1+b_{21}^{-1}z_2), \label{znseval1}\\
z_2^{''} & = & (\beta_1z_1+\beta_2z_2+\beta_3z_3)'+\lam z_1, \label{znseval2}\\
z_3^{''} & = & (\eta_1z_1+\eta_2z_2+\eta_3z_3)'+\lam g(z_1,z_2,z_3), \label{znseval3}
\end{eqnarray}
where the linear function $g$ is given by
\begin{equation}g(z_1,z_2,z_3)=(b_{21}^{-1}b_{33}^{-1}b_{31}b_{22}-b_{33}^{-1}b_{32})z_1-(b_{33}^{-1}b_{21}^{-1}b_{31})z_2+b_{33}^{-1}z_3,
\end{equation}
and the $\beta$ and $\eta$ coefficients can be calculated in terms of the $\alpha_{ij}$ and $b_{ij}$. The exact form these coefficients is not used below, so we omit the calculation.
From \eqref{znseval1}-\eqref{znseval3} it is a simple matter to recast the eigenvalue equation as a first order system of the form
\begin{equation*}
Z'=\mathbb{A}(x,\lam)Z,\quad Z=(z_1,z_2,z_3,z_2',z_3')^{tr},
\end{equation*}
with a corresponding limiting system
\begin{equation*}
Z'=\mathbb{A}_{\pm}(\lam)Z
\end{equation*}
at each of $\pm\infty$. The matrix $\mathbb{A}$ takes the form
\begin{equation*}
\mathbb{A}(x,\lam)=\left(\begin{array}{ccccc} -\lam b_{21}^{-1}b_{22} & \lam b_{21}^{-1} & 0 & 0 & 0 \\ 0 & 0 & 0 & 1 & 0 \\ 0 & 0 & 0 & 0 & 1 \\ \lam+\beta_1'-\lam b_{21}^{-1}b_{22} & \beta_2'+\lam b_{21}^{-1} & \beta_3' & \beta_2 & \beta_3 \\ \mathbb{A}_{51} & \mathbb{A}_{52} & \eta_3'+\lam b_{33}^{-1} & \eta_2 & \eta_3 
\end{array}\right),
\end{equation*}
where 
\begin{equation*}
\mathbb{A}_{51}=\eta_1'-\lam b_{21}^{-1}b_{22}+\lam (b_{21}^{-1}b_{33}^{-1}b_{31}b_{22}-b_{33}^{-1}b_{32})
\end{equation*}
and
\begin{equation*}
\mathbb{A}_{52}=\eta_2'+\lam b_{21}^{-1}-\lam b_{33}^{-1}b_{21}^{-1}b_{31}.
\end{equation*}
%
We verify the consistent splitting hypotheses, without loss of generality, in the original $w$-coordinates. The characteristic equation has the form
$$(\lam I+\mu A_{\pm}-\mu^2B_{\pm})v=0.$$ 
We obtain then, a sequence of lemmas.
\begin{lemma} For $\re\lambda >0$ the matrix $\mathbb{A}_{\pm}(\lambda)$ has eigenvalues
$$\mu_1^{\pm}(\lambda),\mu_2^{\pm}(\lam)<0<\mu_3^{\pm}(\lambda),\mu_4^{\pm}(\lambda),\mu_5^{\pm}(\lam), $$
(with ordering referring to real parts). The eigenspaces $\mathcal{S}^{\pm}(\lambda)$ and $\mathcal{U}^{\pm}(\lambda)$ associated with the eigenvalues $\mu_1^{\pm}(\lambda),\mu_2^{\pm}(\lam)$ and $\mu_3^{\pm}(\lambda),\mu_4^{\pm}(\lambda),\mu_5^{\pm}(\lam)$ respectively depend analytically on $\lambda$
\end{lemma}
%
\pf As established in the introduction, the number of positive/negative roots is constant for $\re\lam>0$, so that roots can be counted as $\lam \rightarrow +\infty$ along the real axis. There is one root with $\mu\sim \lam$ and four roots with $\mu\sim \lam^{1/2}$. See the Appendix for more details. \foorp
\noindent
Moreover in a neighborhood of $\lam=0$, a bifurcation analysis yields:
\begin{lemma}For each $j$, there are analytic extensions of $\mu_j^{\pm}(\lambda)$ to a neighborhood $N$ of $\lambda=0$. Moreover there are analytic choices of individual eigenvectors $V^{\pm}_j(\lambda)$ corresponding to $\mu_j^{\pm}(\lambda)$ in $N$. 
\end{lemma}
\pf When $\lam =0$ the characteristic equation reduces to $(\mu A_{\pm}-\mu^2B_{\pm})v=0$ which has a triple root at zero.
Nonzero roots must satisfy $\mu^{-1}\in\sigma (A^{-1}B)$ which has two nonzero eigenvalues. One of them switches signs at $\pm\infty$. The zero roots bifurcate analytically from zero, for linearizing about $(\lam, \mu)=(0,0)$, we obtain
$$(\lam I+\mu A_{\pm})v=0,$$ 
we find that  on the $-\infty$ side
$$\mu_{2,3,4}^{-}(\lam)=-\frac{\lam}{a_{1,2,3}^{-}}+\mathit{O}(\lam^2),$$
and 
$$v_{2,3,4}^-=r_{1,2,3}^-+\mathit{O}(\lam).$$ 
While on the $+\infty$ side
$$\mu_{3,4,5}^{+}(\lam)=-\frac{\lam}{a_{1,2,3}^{+}}+\mathit{O}(\lam^2),$$
and 
$$v_{3,4,5}^-=r_{1,2,3}^++\mathit{O}(\lam).$$
Then $V_j^{\pm}$ correspond to $v_j^{\pm}$ when the limiting eigenvalue equation is written as a first order system.
\foorp
\begin{lemma}
There are choices of bases
$$\mathcal{B}_{\mathcal{S}}^{\pm}(\lambda)=\{\phi_1^{\pm}(\lam),\phi_2^{\pm}(\lam)\}$$ 
and 
$$\mathcal{B}_{\mathcal{U}}^{\pm}(\lambda)=\{\phi_3^{\pm}(\lam),\phi_4^{\pm}(\lam),\phi_5^{\pm}(\lam)\}$$
of $\mathcal{S}^{\pm}(\lambda)$ and $\mathcal{U}^{\pm}(\lambda)$ which are analytic with respect to \lam\ in $N\cup \{\re \lam>0\}$. In the neighborhood $N$, they satisfy 
$$\mathcal{B}_{\mathcal{S}}^+(\lambda)=k_+(\lambda)V_1^+(\lambda)\wedge V_2^+(\lam),$$
and
$$\mathcal{B}_{\mathcal{U}}^-(\lambda)=k_-(\lam)V_3^-(\lam)\wedge V_4^-(\lam)\wedge V_5^-(\lam)$$
where $V_j^{\pm}$ are as in the previous lemma and $k_{\pm}(\lam)$ are scalar functions such that $k_{\pm}(0)=1$.
\end{lemma} 
\pf The proof follows from the previous lemma and a standard (nontrivial) result of matrix perturbation theory \cite{Ka}. \foorp
\noindent
Finally using the gap lemma (Theorem~\ref{gaplemma}), we obtain
\begin{lemma}There are bases $\mathcal{B}_{\mathcal{S}}(x,\lam)$ and $\mathcal{B}_{\mathcal{U}}(x,\lam)$ of the spaces of solutions of the eiegnvalue equations decaying at $x=\pm\infty$ which are tangent to $\mathcal{S}^+(\lam)$ as $x\rightarrow +\infty$ and $\mathcal{U}^-(\lam)$ as $x\rightarrow -\infty$. That is 
$$\mathcal{S}^+(\lam)=\lim_{x\rightarrow +\infty}\spn\mathcal{B}_{\mathcal{S}}(x,\lam),$$
and
$$\mathcal{U}^-(\lam)=\lim_{x\rightarrow -\infty}\spn\mathcal{B}_{\mathcal{U}}(x,\lam).$$
\end{lemma}
\noindent
A word on notation is in order. Working again in $z$-coordinates, we refer to the elements of the bases $\mathcal{B}_{\mathcal{S}}(x,\lam)$ and $\mathcal{B}_{\mathcal{U}}(x,\lam)$ %
by $Z_j^{\pm}(x,\lam)$, so that 
$$Z_j^{\pm}=(z_{1,j}^{\pm},z_{2,j}^{\pm},z_{3,j}^{\pm},z_{2,j}^{\pm '},z_{3,j}^{\pm '})^{tr},$$
and we denote by $z_j^{\pm}$ with a solitary subscript the first three components of $Z_j^{\pm}$. Thus
$$z_j^{\pm}=(z_{1,j}^{\pm},z_{2,j}^{\pm},z_{3,j}^{\pm})^{tr}.$$
%
\subsection{The Evans Function}
\begin{defn} The \emph{Evans function} is
$$\evans =\det(Z_1^+,Z_2^+,Z_3^-,Z_4^-,Z_5^-)|_{x=0}=\det\begin{pmatrix} z_{1,1}^+ & z_{1,2}^+ & z_{1,3}^- & z_{1,4}^- & z_{1,5}^- \\ z_{2,1}^+ & z_{2,2}^+ & z_{2,3}^- & z_{2,4}^- & z_{2,5}^- \\ z_{3,1}^+ & z_{3,2}^+ & z_{3,3}^- & z_{3,4}^- & z_{3,5}^- \\ z_{2,1}^{+'} & z_{2,2}^{+'} & z_{2,3}^{-'} & z_{2,4}^{-'} & z_{2,5}^{-'} \\ z_{3,1}^{+'} & z_{3,2}^{+'} & z_{3,3}^{-'} & z_{3,4}^{-'} & z_{3,5}^{-'} \end{pmatrix}_{|_{x=0}}.$$
\end{defn}
\noindent
As usual we are free to put at $\lam =0$
\begin{equation}
z_1^+=z_5^-=C\bar{U}_x,
\label{normalize1} \end{equation} 
and at $\lam =0$ we also set
\begin{equation}z_2^+(+\infty)=0,\quad z_3^-(-\infty)= Cr_1^-,\quad  z_4^-(-\infty)=Cr_2^-.
\label{normalize2} \end{equation}
\subsubsection{Calculation of $D'(0)$}
%
\begin{prop}The Evans function \evans\ satisfies $D(0)=0$ and 
$$\sgn D'(0)=\sgn \gamma_{\text{NS}}\det(r_1^-,r_2^-,[U]),$$ 
where 
$$\gamma_{\text{NS}}=\left(\begin{array}{cc} z_{2,1}^+ & z_{2,2}^+ \\ z_{3,1}^+ & z_{3,2}^+\end{array}\right),$$
which measures transversality of the intersection of stable/unstable manifolds in the traveling wave ODE.
\label{nsdprime}
\end{prop}
%
%
\pf That $D(0)=0$ follows immediately from (\ref{normalize1}). Using the Leibniz rule to compute $D'(0)$, we find
\begin{multline*}D'(0)=\det(\partial_{\lam}Z_1^+,Z_2^+,Z_3^-,Z_4^-,Z_5^-)|_{x=0}+\cdots\\ +\det(Z_1^+,Z_2^+,Z_3^-,Z_4^-,\partial_{\lam}Z_5^-)|_{x=0}.\end{multline*}
We combine the two nonzero determinants above
\begin{equation} D'(0) =\det(Z_1^+,Z_2^+,Z_3^-,Z_4^-,\tilde{Z})|_{x=0},
\label{NSEdprime} 
\end{equation}
where 
$$\tilde{Z}=\partial_{\lam}(Z_5^--Z_1^+).$$ 
Differentiating the eigenvalue equation with respect to \lam\ leads to the equations satisfied by $\tilde{z}$.
\begin{equation}
B(C^{-1}\tilde{z})' = AC^{-1}\tilde{z}+[U]. \label{integratedvariational}
\end{equation}
Also at $\lam =0$, the eigenvalue equations simplify considerably for $j=1,2,3,4$ to (omitting $\pm$)
\begin{eqnarray*}
0 & = & z_{1,j}', \\
z_{2,j}^{''} & = & (\beta_1z_{1,j}+\beta_2z_{2,j}+\beta_3z_{3,j})', \\
z_{3,j}^{''} & = & (\eta_1z_{1,j}1+\eta_2z_{2,j}+\eta_3z_{3,j})',
\end{eqnarray*}
which can be integrated up using the boundary conditions supplied by (\ref{normalize1}) and (\ref{normalize2}). Thus
\begin{eqnarray}
B(C^{-1}z_j^+)' & = & AC^{-1}z_j^+ \quad j=1,2, \label{integratedfastmodes}\\
B(C^{-1}z_3^-)' & = & AC^{-1}z_3^- -a_1^-r_1^- ,\label{integratedslowmode1}\\
B(C^{-1}z_4^-)' & = & AC^{-1}z_4^--a_2^-r_2^- .\label{integratedslowmode2}
\end{eqnarray}
The first equation of each of (\ref{integratedvariational})-(\ref{integratedslowmode2}) allows a simplification in the first row of the  determinant (\ref{NSEdprime}). Namely,
\begin{eqnarray*}
\tilde{z}_1 & = & -[\rho], \\
z_{1,j}^+ & = & 0 \quad j=1,2, \\
z_{1,3}^- & = & a_1^-(r_1^-)_1, \\
z_{1,4}^- & = & a_2^-(r_2^-)_1. 
\end{eqnarray*}
The second equation of each of (\ref{integratedvariational})-(\ref{integratedslowmode2}) allows a row operation to simplify the fourth row,
\begin{eqnarray*}
\tilde{z}_2' & = &\beta_1\tilde{z}_1+\beta_2\tilde{z}_2+\beta_3\tilde{z}_3 +[m], \\
z_{2,j}^{+'} & = &\beta_1 z_{1,j}^++\beta_2 z_{2,j}^++\beta_3z_{3,j}^+,\quad j=1,2, \\
 z_{2,3}^{-'} & = &\beta_1 z_{1,3}^-+\beta_2 z_{2,3}^- +\beta_3z_{3,3}^- -a_1^-(r_1^-)_2, \\
z_{2,4}^{-'} & = &\beta_1 z_{1,4}^-+\beta_2 z_{2,4}^- +\beta_3z_{3,4}^- -a_2^-(r_2^-)_2,
\end{eqnarray*}
while the third equation indicates that a row operation will simplify the fifth row of (\ref{NSEdprime}).
\begin{eqnarray*}
\tilde{z}_3' & = &\eta_1\tilde{z}_1+\eta_2\tilde{z}_2+\eta_3\tilde{z}_3 +[\E], \\
z_{3,j}^{+'} & = &\eta_1 z_{1,j}^++\eta_2 z_{2,j}^++\eta_3z_{3,j}^+,\quad j=1,2, \\
 z_{3,3}^{-'} & = &\eta_1 z_{1,3}^-+\eta_2 z_{2,3}^- +\eta_3z_{3,3}^- -a_1^-(r_1^-)_3, \\
z_{3,4}^{-'} & = &\eta_1 z_{1,4}^-+\eta_2 z_{2,4}^- +\eta_3z_{3,4}^- -a_2^-(r_2^-)_3,
\end{eqnarray*}
where we've used $(r_j^-)_i$ to denote the $i$th component of $r_j^-$.
Putting these operations together, (\ref{NSEdprime}) simplifies to 
\begin{equation}
\det\left(\begin{array}{ccccc} 0 & 0 & a_1^- & a_2^- & -[\rho] \\ z_{2,1}^+ & z_{2,2}^+ & * & * & * \\ z_{3,1}^+ & z_{3,2}^+ & * & * & * \\ 0 & 0 & -a_1^-(r_1^-)_2 & -a_2^-(r_2^-)_2 & [m] \\ 0 & 0 & -a_1^-(r_1^-)_3 & -a)2^-(r_2^-)_3 & [\E] \end{array}\right)_{|_{x=0}}.
\label{simplified}
\end{equation}
From (\ref{simplified}) it follows that 
\begin{equation*}
D'(0)=a_1^-a_2^-\det\left(\begin{array}{cc}z_{2,1}^+ & z_{2,2}^+ \\ z_{3,1}^+ & z_{3,2}^+ \end{array}\right)\det(r_1^-,r_2^-,[U]).
\end{equation*}
From the shock inequalities \eqref{3shock} it follows that $a_1^-a_2^->0$, and the proposition is proved.
\foorp
%
%
%
%
%
\subsubsection{Large \lam\ Behavior}
We appeal to the Appendix to determine $\sgn D(\lam)$ as $\lam\to +\infty$ along the real axis and complete the calculation of the stability index. In the Appendix the calculation is carried out for abstract ``real viscosity'' systems of the form:
\begin{equation*}
U_t+F(U)_x=(B(U)U_x)_x,
\end{equation*}
where 
$$U=\left(\begin{array}{c} u \\ v\end{array}\right),\quad F=\left(\begin{array}{c} f \\ g \end{array}\right),\quad B=\left(\begin{array}{cc} 0 & 0 \\ b_1 & b_2 \end{array}\right), $$ 
and
$$u,f\in \R^{n-r},\quad v,g\in\R^r,\quad b_1\in\R^{r\times (n-r)},\quad b_2\in \R^{r\times r}.$$
Noting the abuse of notation, for this section we adopt the notation of Appendix. Then from Lemma \ref{Clemma1} it follows that for \lam\ real and sufficiently large,
$$D(\lam)\ne 0.$$ 
Furthermore we can relate $\sgn D(\lam)$ for \lam\ large to our normalizations at $\lam =0$ by Lemma \ref{Dinfty}. We find
for \lam\ real and sufficiently large that 
\begin{equation}
\sgn \tilde{D}(\lam)=\sgn \mathbb{S}^+\det(\pi\mathbb{W}^+,\varepsilon\mathbb{S}^+)\det(\pi\mathbb{W}^-)|_{\lam=0},
\label{nspullback}
\end{equation}
where $\tilde{D}$ is the Evans function computed in the original $w$-coordinates. $D$ and $\tilde{D}$ differ by a nonvanishing real factor. In \eqref{nspullback}, $\mathbb{S}^+$ is a basis for the one-dimensional stable subspace of $\tilde{A}=A_{11}-A_{12}b_2^{-1}b_1$, where $A_{ij}$ are the entries of the Jacobian matrix of $F$, in particular $(A_{11},A_{12})=df(\bar{U})$. In this case $\tilde{A}$ is simply the particle velocity (the ``original'' $u$) and is thus negative for a 3-shock. $\mathbb{W}^{\pm}$ denote bases for decaying solutions at each of $\pm\infty$, and $\pi,\varepsilon$ denote projection and extension respectively.

Working now in $z$-coordinates and following the discussion of the Appendix, the form of \eqref{nspullback} simplifies considerably, and $D(\lam)$ satisfies
\begin{lemma}For \lam\ real and sufficiently large
\begin{equation}
\sgn D(\lam) =\sgn \gamma_{\text{NS}} (\mathbb{S}^+)^2\det(r_1^-,r_2^-,\bar{U}_x). \label{nslargelam}
\end{equation}
\label{nsinfty}
\end{lemma}
%
\subsection{The Stability Index}
Therefore combining Lemma \ref{nsinfty} and Proposition \ref{nsdprime} we have
\begin{prop}
The stability index for a viscous Lax 3-shock is
$$\sgn D'(0)D(+\infty)=\sgn \det(r_1^-,r_2^-,[U])\det(r_1^-,r_2^-,\bar{U}_x|_{-\infty}),$$
where $[U]=([\rho],[m],[\E])^{tr}$ is a vector of jumps and $\bar{U}_x=(\bar{\rho}_x,\bar{m}_x,\bar{\E}_x)^{tr}$.
\end{prop}
\noindent
We remark that the stability index is unaffected by change in the viscosity matrix. That is, the index agrees with that of an artificially parabolic system at least in the presence of \eqref{symmetrizability}-\eqref{blockstructure}. Further we note that the calculation of $D'(0)$ captures low-frequency information. In this setting, this corresponds to ``inviscid'' behavior or the stability of shocks solutions of the Euler Equations,  
\begin{eqnarray*}
\rho_t + (\rho u)_x & = & 0, \label{eulermass} \\
(\rho u)_t + (\rho u^2+p)_x & = & 0, \label{eulermomentum} \\
\E_t+ u\left(\E+p\right)_x & = & 0.\label{eulerenergy}
\end{eqnarray*}
and thus $\Delta$ takes the form of a Lopatinski determinant, familiar from the stability analysis of such shocks. See \cite{JL} for the calculation of this determinant for the multi-dimensional Euler equations. In the weak-shock limit we note that
$$[U]\sim \bar{U}_x\sim r_3^-,$$
which implies that the stability index satisfies
$$\sgn D'(0)D(+\infty) \sim \sgn\det(r_1^-,r_2^-,r_3^-)^2=+1$$
in the weak shock limit consistent with stability. Moreover in the ideal gas case, combining the nonvanishing of $\Delta$ which is well known for an ideal gas, for example see \cite{MBOOK} or \cite{serretransition}, the global existence result of \cite{Gi} which guarantees transversality of connections, and nonvanishing of $D(+\infty)$, we can then conclude consistency with stability for shocks of any strength in the ideal gas case.

We remark that in the more general case, some rudimentary knowledge about the connecting orbit is neccessary to evaluate
$$\sgn\det (r_1^-,r_2^-,\bar{U}_x),$$
namely the direction of $\bar{U}_x|_{-\infty}$. 
%
%
\subsection{Isentropic Gas Dynamics}
\begin{figure}[ht]
\centerline{\epsfxsize10truecm\epsfbox{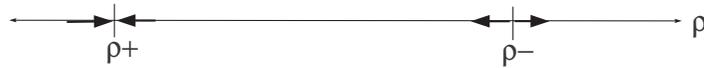}}
\caption{The Phase Portrait}
\label{2shockphase}
\end{figure}

We note that the in case of isentropic gas dynamics,
\begin{align}
\rho_t&+(\rho u)_x=0, \\
(\rho u)_t&+(\rho u^2+p)_x=(\nu u_x)_x, 
\end{align}
where the pressure satisfies
$$p=p(\rho),\quad c^2=p'(\rho)>0,$$
an analogous but simpler calculation leads to the stability index
\begin{eqnarray*}
\tilde{\Gamma}=\sgn D'(0)D(+\infty) & = & \sgn \gamma_i^2\det(\mathbb{S}^+)^2\det(r_1^-,[U])\det(r_1^-,\bar{U}_x) \\            
& = & \sgn \det(r_1^-,[U])\det(r_1^-,\bar{U}_x),
\end{eqnarray*}
for a 2-shock. Here  $[U]=([\rho],[m])^{tr}$ is a vector of jumps in the gas-gynamical conserved quantities, $\bar{U}_x=(\bar{\rho}_x,\bar{m}_x)^{tr}$, $r_1^-=(1,u_--c_-)^{tr}$ is the outgoing right eigenvector, and the term $\gamma_i$ is a transversality coefficient.
In this case, evaluation of the index is straightforward since the phase space for the traveling wave ODE is one-dimensional. The phase portrait is shown in Figure~\ref{2shockphase}.
Then $\sgn\det(r_1^-,[U])\det(r_1^-,\bar{U}_x)$ becomes simply
$$\sgn [\rho]\bar{\rho}_x(u_--c_-)^2=+1,\text{consistent with stability,}$$
since both $[\rho]$ and $\bar{\rho}_x$ are negative. Moreover, we note that for one-dimensional isentropic gas dynamics, Humpherys \cite{Hu} has shown that there are no unstable \emph{real} eigenvalues for shocks of arbitrary strength.

%
%
%
\section{Reacting Flow}\label{RNSsection}
In this section we extend the analysis of the previous section to our main interest: strong detonation waves, which are particular traveling wave solutions of equations~(\ref{eq:mass})-(\ref{eq:progress}). We assume ``ignition temperature kinetics,'' that is, $\varphi$ vanishes for temperatures below some ignition temperature and is identically 1 for some larger value of $T$. We also assume that the pressure $p$ and the temperature $T$ are given functions of $\rho$ and $e$, and thus are independent of the progress of the reaction.  At some points in the analysis we further specify that the gas is ideal and polytropic, so 
$$p=\Gamma \rho e,\quad T=c_v^{-1}e.$$
Here $\tilde{e}$ is related to $e$, the specific internal energy not due to reaction, by 
$$\tilde{e}=e+qY.$$
Using this relation, the third and fourth equations can be combined, and the energy balance equation can be rewritten as
$$ \left( \rho \left(\frac{u^2}{2}+e\right)\right)_t+\left(\rho u(\frac{u^2}{2}+e)+up\right)_x = (\theta T_x)_x+qk\rho Y\varphi(T)+(\nu uu_x)_x.$$
Using this form of the energy balance equation, we note that the system fits in the general framework (\ref{generalkinematic})-(\ref{generalreaction}) discussed below with 
$U=(\rho, m,\E)^{tr}$, $z=\rho Y$, $Q=\vec{q}=(0,0,q)^{tr}$, $\Phi(U)=k\varphi(T)$, $D^1(U)=d$, and  $D^2(U,z)=(-\frac{dz}{\rho},0,0)$. We remark at this point that since the flux $f$ in the kinematic variables $U$ is as in the Navier-Stokes case, the left and right eigenvectors $l_j$,$r_j$ of the flux Jacobian $A(U)$ and their corresponding eigenvalues $a_j$ are precisely as calculated in the previous section

In accordance with our previous analysis, we assume that the strong detonation has Lax 3-shock structure, so 
$$a_2^-<s<a_3^-,\quad s<a_3^+.$$
We will without loss of generality set $s=0$, so the reaction front is stationary. The shock inequalities imply that $u_{\pm}<0$, so that fluid particles cross the reaction front from right to left, or alternatively the front is ``moving'' to the right connecting an unburned state at $+\infty$ to a completely burned state at $-\infty$. Without loss of generality we normalize $\rho_+=1$ so that the total reactant variable $z$ satisfies
\begin{equation}
z_+=1,\quad z_-=0. \label{zendstates}
\end{equation}
We also assume that the endstates are such that the temperature on the unburned side is below ignition so $\varphi=0$ and that the temperature on the burned side is sufficiently large so that $\varphi =1$. We'll see that this first assumption is necessary to ensure that the end state at $+\infty$ is a rest point of the traveling wave ODE. Thus $\varphi$ satisfies
\begin{equation}
\varphi_+=0,\quad \varphi_-=1. \label{phiendstates}
\end{equation}

%
%
%
%
%
%
\subsection{Traveling Wave ODE and Linearized Equations}
The ($s=0$) traveling wave ODE is
\begin{align}
&\;m_x  =  0, \\
&\left(\frac{m^2}{\rho}+p\right)_x =  \left(\nu\left(\frac{m}{\rho}\right)_x\right)_x, \\
&\left(\frac{m\tilde{\E}}{\rho}+\frac{mp}{\rho}\right)_x = (\theta T_x)_x+\left(q\rho d\left(\frac{z}{\rho}\right)_x\right)_x+\left(\nu\left(\frac{m}{\rho}\right)\left(\frac{m}{\rho}\right)_x\right)_x, \\
&\left(\frac{mz}{\rho}\right)_x=(dz_x)_x+\left(-\frac{dz}{\rho}\rho_x\right)_x-k\varphi(T)z.
\end{align}
Fixing a state at $-\infty$, the first three equations above can be integrated once to 
\begin{align*}
&m-m_- = 0, \\
&\left(\frac{m^2}{\rho} +p\right)-\left(\frac{m^2}{\rho}+p\right)_- =\nu\left(\frac{m}{\rho}\right)_x, \\
&\left(\frac{m}{\rho}(\tilde{\E}+p)\right)-\left(\frac{m}{\rho}(\tilde{\E}+p)\right)_- = \theta T_x+q\rho d\left(\frac{z}{\rho}\right)_x+\nu\left(\frac{m}{\rho}\right)\left(\frac{m}{\rho}\right)_x.
\end{align*}
Requiring the state at $+\infty$ to be a rest point of the ODE yields from the first two equations above the familiar Rankine-Hugoniot conditions
%
%
%
\begin{align}
[m] & =  0, \label{RH1} \tag{RH1}\\
\left[\frac{m^2}{\rho}+p\right] & =  0,\label{RH2} \tag{RH2}
\end{align}
where the brackets as usual indicate the difference between the state at $+\infty$ and that at $-\infty$. In the case of the third equation, we find
\begin{equation*}
\left[\frac{m\tilde{\E}}{\rho}+\frac{mp}{\rho}\right]=\left[\frac{m\E}{\rho}+\frac{mp}{\rho}\right]+q\left[\frac{mz}{\rho}\right],
\end{equation*}
since $\tilde{\E}=\E+qz$. Moreover from (\ref{zendstates}), the third jump condition becomes
\begin{equation}
0=\left[\frac{m\E}{\rho}+\frac{mp}{\rho}\right]+qu_+.\label{RH3} \tag{RH3}
\end{equation}
%
%
%
%
%
%
Linearizing about a profile $(\bar{\rho},\bar{m},\bar{\E},\bar{z})$, we find 
\begin{eqnarray}
\rho_t+m_x & = & 0, \label{lindet1} \\
m_t+(\alpha_{21}\rho+\cdots+\alpha_{23}\E)_x & = & (b_{21}\rho_x+b_{22}m_x)_x, \label{lindet2} \\
\E_t+(\alpha_{31}\rho+\cdots+\alpha_{33}\E)_x & = & (b_{31}\rho_x+\cdots+b_{33}\E_x)_x+qkl, \label{lindet3} \\
z_t+(v_1\rho+v_2m+v_4z)_x & = & (dz_x)_x+(\tilde{d}\rho_x)_x-kl, \label{lindet4}
\end{eqnarray}
where $\alpha_{ij}$ and $b_{ij}$ are as in the previous section, and 
\begin{equation*}
l=\{l_zz+l_{\rho}\rho+l_mm+l_{\E}\E\},
\end{equation*} 
with
\begin{eqnarray*}
l_z(x) & = & \varphi(\bar{T}), \\
l_{\rho}(x) & = & \varphi'(\bar{T})(\bar{T}_{\rho}+\bar{T}_e\bar{e}_{\rho}), \\
l_m(x) & = & -\varphi'(\bar{T})\bar{T}_e\frac{\bar{u}}{\bar{\rho}},\\
l_{\E}(x) & = & \varphi'(\bar{T})\bar{T}_{e}/\bar{\rho},
\end{eqnarray*}
and
\begin{eqnarray*}
v_1(x) & = & -\bar{u}\bar{Y}-\frac{d\bar{z}\bar{\rho}_x}{\bar{\rho}^2}, \\
v_2(x) & = & \frac{\bar{z}}{\bar{\rho}}, \\
v_4(x) & = & \bar{u}+\frac{d\bar{\rho}_x}{\bar{\rho}}.
\end{eqnarray*}
We note that due to the structure of $\varphi$ and (\ref{phiendstates}), it follows that $l_{\rho}$, $l_m$, and $l_{\E}$ vanish at both $\pm\infty$ while $l_{z+}=0$ and $l_{z-}=1$. The equations (\ref{lindet1})-(\ref{lindet4}) can also be written in the more compact form 
\begin{eqnarray*}
(Bw')' & = & (Aw)'+w_t+\vec{q}kg(w,z), \\
(dz')'+(\tilde{d}w')' & = & (V_ww)'+(V_zz)'+z_t-kg(w,z). 
\end{eqnarray*}
In an abuse of notation we've written $\tilde{d}$ to stand for both the $1\times 3$ matrix $(-\frac{dz}{\rho},0,0)$ and the $(1,1)$-entry of that matrix. The meaning will be clear from the context.

%
%
%
%
Under the invertible change of coordinates
\begin{equation}
\left(\begin{array}{c} \z_1 \\ \z_2 \\ \z_3 \\ \z_4 \end{array}\right)=\underbrace{\left(\begin{array}{cccc} 0 & 1 & 0 & 0 \\ b_{21} & b_{22} & 0 & 0 \\ b_{31} & b_{32} & b_{33} & 0 \\ \tilde{d} & 0 & 0 & d \end{array}\right)}_C\left(\begin{array}{c} \rho \\ m \\ \E \\ z \end{array}\right),\label{zetacoord}
\end{equation}
the eigenvalue equation corresponding to (\ref{lindet1})-(\ref{lindet4}) takes the simple form
\begin{eqnarray}
0 & = & \z_1'+\lam(-b_{21}^{-1}b_{22}\z_1+b_{21}^{-1}\z_2), \label{zetaeval1}\\
\z_2^{''} & = & (\beta_1\z_1+\cdots +\beta_3\z_3)'+\lam\z_1, \label{zetaeval2}\\
\z_3^{''} & = & (\eta_1\z_1+\cdots+\eta_3\z_3)'+\lam g(\hat{\z}) +qkl\cdot \z, \label{zetaeval3}\\
\z_4^{''} & = & (\theta_1\z_1+\cdots+\theta_4\z_4)'+\lam h(\z)-kl\cdot \z .\label{zetaeval4}
\end{eqnarray}
We have used the notation $\z=(\z_1,\z_2,\z_3,\z_4)$ and $\hat{\z}=(\z_1,\z_2,\z_3)$. The linear functions $g$ and $h$ take the form
$$
g(\hat{\z})=(b_{21}^{-1}b_{33}^{-1}b_{31}b_{22}-b_{33}^{-1}b_{32})\z_1-(b_{33}^{-1}b_{21}^{-1}b_{31})\z_2+b_{33}^{-1}\z_3,
$$
and 
$$
h(\z)=d^{-1}b_{21}^{-1}b_{22}\tilde{d}\z_1-d^{-1}b_{21}^{-1}\tilde{d}\z_2+d^{-1}\z_4.
$$
The coefficients $\beta_j$ and $\eta_j$ are as in the previous chapter, while the $\theta_j$ depend on $d,\tilde{d},b_{ij}$ and $v_j$. Also $l\cdot \z=(l_1\z_1+\cdots +l_4\z_4)$ where $l_j$ depends on $l_{\rho},l_m,l_{\E},l_z,b_{ij},d$ and $\tilde{d}$.

The change of coordinates matrix $C$ has block structure respecting the division of the variables into gas dynamical variables $w=(\rho ,m, \E)$, and reaction variable $z$. Thus we have 
$$\z=\left(\begin{array}{c|c} C_{NS} & 0 \\ \hline \begin{array}{ccc} \tilde{d} & 0 & 0 \end{array} & d\end{array}\right)\left(\begin{array}{c} w \\ z\end{array}\right),$$
where $C_{NS}$ denotes the change of variables in the nonreacting gas dynamics case considered in the previous section. The inverse of the matrix $C$ also respects this structure and is given by
\begin{equation*}
C^{-1}=\left(\begin{array}{c|c} C_{NS}^{-1} & 0 \\ \hline \begin{array}{ccc} d^{-1}b_{21}^{-1}b_{22}\tilde{d}^{-1} & -d^{-1}b_{21}^{-1}\tilde{d} & 0\end{array} & d^{-1}\end{array}\right).
\end{equation*} 
We note that when $\lam =0$, by a substitution of (\ref{zetaeval4}) into (\ref{zetaeval3}) through the term $l~\cdot~\z$ which appears in both equations, we recover a third equation in which every term is differentiated; this means that we can integrate up the first three equations subject to appropriate boundary conditions. This key fact is what will allow for a simplification of the Evans function determinant via row operations in the calculation of $D'(0)$.
The eigenvalue equations (\ref{zetaeval1})-(\ref{zetaeval4}), can be written as a first order system of the form
\begin{equation}
Z'=\mathbb{A}(x,\lam)Z,\quad Z=(\z_1,\z_2,\z_3,\z_4,\z_2',\z_3',\z_4')^{tr},\label{detevalfirst}
\end{equation}
with corresponding limiting constant-coefficient system
\begin{equation}
Z'=\mathbb{A}_{\pm}(\lam)Z. \label{detevallimit}
\end{equation}
%
%
%
%
The characteristic equation in original $(w,z)$-coordinates (with $w=(w_1,w_2,w_3)^{tr}$) is 
$$ M_{\pm}\left(\begin{array}{c} w \\ z\end{array}\right)=\left(\begin{array}{c} 0 \\ 0 \end{array}\right)$$
where the characteristic matrix $M_{\pm}$ has block structure
\begin{equation*}
\left(\begin{array}{c|c} M^{ul}_{\pm} & M^{ur}_{\pm} \\ \hline M^{ll}_{\pm} & M^{lr}_{\pm} \end{array}\right),
\end{equation*}
with blocks
\begin{equation*}
 M^{ul}_{\pm}=\begin{pmatrix}\lam & \mu & 0 \\ \mu a_{21}-\mu^2b_{21} & \lam+\mu a_{22}-\mu^2b_{22} & \mu a_{23} \\ \mu a_{31}-\mu^2b_{31} & \mu a_{32}-\mu^2b_{32} & \lam +\mu a_{33}-\mu^2b_{33} \end{pmatrix},
\end{equation*}
\begin{equation*}
 M^{ur}_{\pm}=\begin{pmatrix} 0 \\ 0 \\ -qk\varphi_{\pm} \end{pmatrix},
\end{equation*}
and
\begin{align*}
M^{ll}_{\pm}&=\begin{pmatrix} \mu v_1-\mu^2\tilde{d} & \mu v_2 & 0\end{pmatrix}, \\
M^{lr}_{\pm}&=\lam -\mu^2 d+\mu v_4+k\varphi_{\pm}.
\end{align*}

As usual we take advantage of consistent splitting and count the number of stable/unstable roots in the limit $\lam\to +\infty$.
Due to the incomplete parabolicity of the viscosity matrix $B$, we expect one root which scales as
$$\mu\sim \tilde{\mu}\lam,\quad \tilde{\mu}\sim 1,$$
hence the matrix $M_{\pm}$ takes the form
\begin{equation*}
\left(\begin{array}{ccc|c}  \lam & \tilde{\mu}\lam & 0  & 0 \\ \tilde{\mu}\lam a_{21}-\tilde{\mu}^2\lam^2b_{21} & \lam+\tilde{\mu}\lam a_{22}-\tilde{\mu}^2\lam^2b_{22} & \tilde{\mu}\lam a_{23} & 0 \\ \tilde{\mu}\lam a_{31}-\tilde{\mu}^2\lam^2b_{31} & \tilde{\mu}\lam a_{32}-\tilde{\mu}^2\lam^2b_{32} & \lam +\tilde{\mu} a_{33}-\tilde{\mu}^2\lam^2b_{33} &  -qk\varphi_{\pm} \\  \hline  \tilde{\mu}\lam u_{\pm}Y_{\pm}-\tilde{\mu}^2\lam^2\tilde{d} & \tilde{\mu}\lam Y_{\pm} & 0 & N \end{array}\right)_{\pm},
\end{equation*}
where
$$ N=\lam -\tilde{\mu}^2\lam^2 d+\tilde{\mu}\lam u_{\pm}+k\varphi_{\pm}.$$
Whereupon dividing each row by the highest power of \lam, we find that in the limit $\lam\to +\infty$, the roots satisfy a block triangular system
\begin{equation}
\left(\begin{array}{ccc|c} 1 & \tilde{\mu} & 0 & 0 \\ b_{21} & b_{22} & 0 & 0 \\ b_{31} & b_{32} & b_{33} & 0 \\ \hline \tilde{d} & 0 & 0 & d \end{array}\right)\begin{pmatrix}w_1 \\ w_2 \\ w_3 \\ z \end{pmatrix}=\begin{pmatrix} 0 \\ 0 \\ 0 \\ 0 \end{pmatrix} \label{hypscaling},
\end{equation}
where the upper left-hand block is easily recognized as the gas dynamics block which appears in the Appendix. Thus there is one root 
$$-\tilde{\mu}^{-1}\in \sigma (u_{\pm}),$$ 
and therefore since we are considering detonations with 3-shock structure, there is one unstable root at each of $\pm\infty$.
The remaining roots scale as
$$\mu \sim \tilde{\mu}\lam^{1/2},\quad \tilde{\mu}\sim 1.$$
Thus, upon substituting and dividing as before, we obtain a different block triangular system 
\begin{equation}
\left(\begin{array}{ccc|c} 1 & 0 & 0 & 0 \\ -\tilde{\mu}^2b_{21} & -\tilde{\mu}^2b_{22}+1 & 0 & 0 \\ \tilde{\mu}^2b_{31} & -\tilde{\mu}^2b_{32} & -\tilde{\mu}^2b_{33}+1 & 0 \\ \hline -\tilde{\mu}^2\tilde{d} & 0 & 0 & -\tilde{\mu}^2d+1 \end{array}\right)\begin{pmatrix}w_1 \\ w_2 \\ w_3 \\ z \end{pmatrix}=\begin{pmatrix} 0 \\ 0 \\ 0 \\ 0 \end{pmatrix}. \label{parscaling}
\end{equation}
It follows immediately from \eqref{parscaling} that $w_1=0$, and thus roots must satisfy 
\begin{equation*}
\left[ \left(\begin{array}{cc|c}b_{22} & 0 & 0 \\ b_{32} & b_{33} & 0 \\ \hline 0 & 0 & d \end{array}\right)-\tilde{\mu}^{-2}I\right]\begin{pmatrix} w_2 \\ w_3 \\ z \end{pmatrix}=\begin{pmatrix} 0 \\ 0 \\ 0 \\ 0 \end{pmatrix},
\end{equation*}
that is, $\tilde{\mu}^{-2}$ must be an eigenvalue of the block diagonal matrix
\begin{equation}
\left(\begin{array}{cc|c}b_{22} & 0 & 0 \\ b_{32} & b_{33} & 0 \\ \hline 0 & 0 & d \end{array}\right).
\label{blockdiag}
\end{equation}
This yields 3 stable and 3 unstable roots at each of $\pm\infty$. We note here that due to the block diagonal structure of \eqref{blockdiag}, there is a 2-1 split at each infinity of kinematic versus reactive roots. Summarizing, we've found that there are 4 stable and 3 unstable roots at each of $\pm\infty$ as long as  $\re\lam>0$.

%
%
%
Consistent splitting breaks down at $\lam =0$. We look first at the $x=+\infty$ case. In this case due to \eqref{phiendstates}, the characteristic equation has lower block triangular structure when $\lam =0$. Thus the characteristic matrix has the form
\begin{equation*}
\left(\begin{array}{ccc|c}  0 & \mu & 0  & 0 \\ \mu a_{21}-\mu^2b_{21} & \mu a_{22}-\mu^2b_{22} & \mu a_{23} & 0 \\ \mu a_{31}-\mu^2b_{31} & \mu a_{32}-\mu^2b_{32} & \mu a_{33}-\mu^2b_{33} & 0 \\  \hline  \mu u_+-\mu^2\tilde{d} & \mu  & 0 & -\mu^2 d+\mu u_+ \end{array}\right).
\end{equation*}
From our analysis of the nonreacting case, we know that the stable roots corresponding to the upper left-hand kinematic block do not vanish at $\lam=0$, thus they correspond to \emph{fast} kinematic modes. On the other hand, looking at the reaction block,
$$-\mu^2d+\mu u_+=0,$$
it's clear that there are solutions $\mu =0$ corresponding to a slow unstable reactive mode and $\mu=\frac{u_+}{d}$, a fast stable reactive mode. A bifurcation analysis as in the previous section shows that the roots have analytic extensions in a neighborhood of \lam =0.

In the $x=-\infty$ case, we find when $\lam =0$ that, due to \eqref{zendstates}, the characteristic matrix is upper block triangular
\begin{equation*}
\left(\begin{array}{ccc|c}  0 & \mu & 0  & 0 \\ \mu a_{21}-\mu^2b_{21} & \mu a_{22}-\mu^2b_{22} & \mu a_{23} & 0 \\ \mu a_{31}-\mu^2b_{31} & \mu a_{32}-\mu^2b_{32} & \mu a_{33}-\mu^2b_{33} & -qk \\  \hline  0 & 0  & 0 & -\mu^2 d+\mu u+k \end{array}\right).
\end{equation*}
From our previous analysis of the nonreacting case, we know that two of the unstable roots from the gas block vanish at $\lam=0$. By block structure, the corresponding vectors have the expansion%
$$\begin{pmatrix} r_j^- \\ 0\end{pmatrix}+O(\lam),\quad j=1,2, $$
where $r_j^-$ is an eigenvector of the flux Jacobian $A$. The roots coming from the reaction block must satisfy
\begin{equation*}
\mu = \frac{u_-}{2d}\mp \frac{\sqrt{u_-^2+4dk}}{2d}\neq 0\;\;\text{since }d,k>0.
\end{equation*}   
This implies that all reactive modes are fast modes on the $x=-\infty$ side. We note here that this discussion shows that all the reactive modes of our interest (stable at $+\infty$/unstable at $-\infty$) are fast modes, and thus asymptotically vanish in both kinematic and reaction components. We note here that the inclusion of the differing block structure at each of $\pm\infty$ corrects an omission in the abstract development considered in \cite{Z4}.

Finally, applying the gap lemma we obtain bases 
$$\{Z_1^+,Z_2^+,Z_3^+\}\;\;\text{and}\;\; \{Z_4^-,Z_5^-,Z_6^-,Z_7^-\}$$
spanning the stable manifold at $+\infty$ and the unstable manifold at $-\infty$. We use the notation
\begin{equation*}
Z_j^{\pm}=(\z_{1,j}^{\pm},\z_{2,j}^{\pm},\z_{3,j}^{\pm},\z_{4,j}^{\pm},\z_{2,j}^{\pm '},\z_{3,j}^{\pm '},\z_{4,j}^{\pm '})^{tr},
\end{equation*}
when we need to indicate the components. Also we note here that we will follow our standard convention and identify $Z_1^+/Z_7^-$ with the derivative of the profile since it satisfies the linearized system and decays at both $\pm\infty$. We'll also let $Z_2^+$ and $Z_3^+$ correspond to fast kinematic and reactive modes respectively. On the $-\infty$ side we'll let $Z_4^-$ and $Z_5^-$ correspond to slow kinematic modes while $Z_6^-$ is a fast kinematic mode. Note that implicit in this assignment is the assumption that the profile approaches the burned endstate parallel to the reactive mode. This assumption is generic in the case where the reaction is slower and actually occurs in the ZND limit as evidenced in our discussion of existence in Section 1. 
%
%
%
%
\subsection{The Evans Function}
\begin{defn} The \emph{Evans function} is
\begin{equation}
\evans =\det(Z_1^+,Z_2^+,Z_3^+,Z_4^-,Z_5^-,Z_6^-,Z_7^-)|_{x=0}.
\end{equation}
\end{defn}
We make the usual normalizations at $\lam =0$. Namely we put 
\begin{equation}
\left(\begin{array}{c}\z_{1,1}^+ \\ \z_{2,1}^{+} \\ \z_{3,1}^{+} \\ \z_{4,1}^{+} \end{array}\right)=\left(\begin{array}{c}\z_{1,7}^- \\ \z_{2,7}^{-} \\ \z_{3,7}^{-} \\ \z_{4,7}^{-}\end{array}\right)=C\left(\begin{array}{c} \bar{\rho}_x \\ \bar{m}_x \\ \bar{\E}_x \\ \bar{z}_x \end{array}\right), \label{detnorm1}
\end{equation} 
and the other fast modes satisfy
\begin{equation}
\left\{\begin{array}{l} Z_2^+(+\infty)=Z_3^+(+\infty)= 0, \\ Z_6^-(-\infty)=0, \end{array}\right. \label{detnorm2}
\end{equation}
while the slow modes satisfy
\begin{equation}
\left(\begin{array}{c}\z_{1,4}^- \\ \z_{2,4}^{-} \\ \z_{3,4}^{-} \\ \z_{4,4}^{-} \end{array}\right)=C\left(\begin{array}{c} r_1^- \\ 0 \end{array}\right),\quad \left(\begin{array}{c}\z_{1,5}^- \\ \z_{2,5}^{-} \\ \z_{3,5}^{-} \\ \z_{4,5}^{-} \end{array}\right)=C\left(\begin{array}{c} r_2^- \\ 0 \end{array}\right). \label{detnorm3}
\end{equation}
%
\subsubsection{Calculation of $D'(0)$}

\begin{prop}The Evans function \evans\ satisfies $D(0)=0$ and $D'(0)=\gamma_d\Delta$ where
\begin{equation}
\Delta=\det (r_1^-,r_2^-,[U]+\vec{q}),
\label{reactdelta}
\end{equation}
and 
\begin{equation}
\gamma_d=\det\left(\begin{array}{cccc} \z_{2,1}^+ & \z_{2,2}^+ & \z_{2,3}^+ & \z_{2,6}^- \\  \z_{3,1}^+ & \z_{3,2}^+ & \z_{3,3}^+ & \z_{3,6}^- \\ \z_{4,1}^+ & \z_{4,2}^+ & \z_{4,3}^+ & \z_{4,6}^- \\ \z_{4,1}^{+'} & \z_{4,2}^{+'} & \z_{4,3}^{+'} & \z_{4,6}^{-'} \end{array}\right)
\label{gammad}
\end{equation}
measures transversality of the stable and unstable manifolds in the traveling wave ODE.
\label{detdprime}
\end{prop}

Before beginning the proof we remark that in contrast to the case considered earlier $\gamma_d$ is not extreme, that is, it involves fast modes from both infinities. 

\pf As usual $D(0)=0$ follows immediately from the normalization (\ref{detnorm1}) chosen for the basis elements. Applying the Leibniz rule,
\begin{equation}D'(0)=\det(\partial_{\lam}Z_1^+,Z_2^+,\ldots ,Z_7^-)|_{x=0}+\cdots +\det(Z_1^+,\ldots ,Z_6^-,\partial_{\lam}Z_7^-)|_{x=0}, 
\end{equation}
and combining the two nonzero determinants in the above equation, we obtain
\begin{equation}
D'(0)=\det(Z_1^+,Z_2^+,Z_3^+,Z_4^-,Z_5^-,Z_6^-,\tilde{Z})|_{x=0},
\end{equation}
where $\tilde{Z}=\tilde{Z^-}-\tilde{Z}^+=\partial_{\lam}Z_{7}^--\partial_{\lam}Z_1^+$. For clarity we perform the necessary manipulations in the original $(w,z)$-coordinates and then translate the results to the $\z$-coordinates. Thus we write down the general form of the eigenvalue equation 
\begin{eqnarray}
(Bw')' & = & (Aw)'+\lam w+\vec{q}k\varphi z \label{wzeval1} \\
(dz')'+(\tilde{d}w')' & = & (V_ww)'+(V_zz)'+\lam z-k\varphi z \label{wzeval2}.
\end{eqnarray}
When $\lam =0$, we can make a substitution from (\ref{wzeval2}) into (\ref{wzeval1}) so that the first equation has the form 
\begin{equation}
(Bw')' = (Aw)'+\vec{q}(-(dz')'-(\tilde{d}w')'+(V_ww)'+(V_zz)'), \label{substituted}
\end{equation}
in which every term is differentiated and hence (\ref{substituted}) can be integrated to 
\begin{equation*}
Bw'-Bw_{\pm}= Aw-Aw_{\pm}+\vec{q}(-dz'+dz_{\pm}-\tilde{d}w'+\tilde{d}w_{\pm}+V_ww-V_ww_{\pm}+V_zz-V_zz_{\pm}),
\end{equation*}
with the $\pm$ subscripts indicating boundary conditions to be supplied by the normalizations (\ref{detnorm1})-(\ref{detnorm3}). Therefore the fast modes satisfy
\begin{equation}
Bw'+\vec{q}dz'+\vec{q}\tilde{d}w' = (A+\vec{q}V_w)w+\vec{q}V_zz, \label{detfast}
\end{equation}
while the slow modes satisfy 
\begin{equation}
Bw'+\vec{q}dz'+\vec{q}\tilde{d}w' = (A+\vec{q}V_w)w+\vec{q}V_zz-a_j^-r_j^-.
\end{equation}
On the other hand $(\tilde{w}^{\pm},\tilde{z}^{\pm})$ satisfy the variational equations at $\lam =0$
\begin{eqnarray}
(B\tilde{w}^{\pm '})' & = & (A\tilde{w}^{\pm})'+\bar{U}_x+\vec{q}k\varphi \tilde{z}^{\pm} \label{detvar1} \\
(d\tilde{z}^{\pm'})'+(\tilde{d}\tilde{w}^{\pm '})' & = & (V_w\tilde{w}^{\pm})'+(V_z\tilde{z}^{\pm})'+\bar{z}_x-k\varphi \tilde{z}^{\pm} \label{detvar2}.
\end{eqnarray}
We make the same substitution from (\ref{detvar2}) into (\ref{detvar1}) so that every term is a derivative. Then we integrate $(\tilde{w}^+,\tilde{z}^+)$ from $x$ to $+\infty$,
\begin{equation}
B\tilde{w}^{+'}+\vec{q}d\tilde{z}^{+'}+\vec{q}\tilde{d}\tilde{w}^{+'} = (A+\vec{q}V_w)\tilde{w}^++\vec{q}V_z\tilde{z}^++\bar{U}-U_++\bar{z}-z_+, \label{plus}
\end{equation} 
and $(\tilde{w}^-,\tilde{z}^-)$ from $-\infty$ to $x$,
\begin{equation}
B\tilde{w}^{-'}+\vec{q}d\tilde{z}^{-'}+\vec{q}\tilde{d}\tilde{w}^{-'} = (A+\vec{q}V_w)\tilde{w}^-+\vec{q}V_z\tilde{z}^-+\bar{U}-U_-+\bar{z}-z_-. \label{minus}
\end{equation} 
It follows by subtracting (\ref{plus}) from (\ref{minus}) that $(\tilde{w},\tilde{z})$ satisfy 
\begin{equation}
B\tilde{w}^{'}+\vec{q}d\tilde{z}^{'}+\vec{q}\tilde{d}\tilde{w}^{'} = (A+\vec{q}V_w)\tilde{w}+\vec{q}V_z\tilde{z}+[U]+\vec{q}
\end{equation} 
Translating this information to $\z$-coordinates, we have for fast modes (kinematic and reactive) $j=1,2,3,6$
\begin{eqnarray*}
0 & = & \z_{1,j} \\
\z_{2,j}' & = & \beta_1\z_{1,j}+\cdots +\beta_3\z_{3,j} \\
\z_{3,j}' & = & \eta_1\z_{1,j}+\cdots +\eta_3\z_{3,j}+qk(-\z_{4,j}'+\theta_1\z_{1,j}+\cdots +\theta_4\z_{4,j}).
\end{eqnarray*}
Note that the equation for $\z_{4,j}''$ is unchanged. We don't get any simplification from that equation.
Also when $j=4$, 
\begin{eqnarray*}
0 & = & \z_{1,4}-a_1^-(r_1^-)_1 \\
\z_{2,4}' & = & \beta_1\z_{1,4}+\cdots +\beta_3\z_{3,4}-a_1^-(r_1^-)_2 \\
\z_{3,4}' & = & \eta_1\z_{1,4}+\cdots +\eta_3\z_{3,4}+qk(-\z_{4,4}'+\theta_1\z_{1,4}+\cdots +\theta_4\z_{4,4})-a_1^-(r_1^-)_3, 
\end{eqnarray*}
and similarly for $j=5$, 
\begin{eqnarray*}
0 & = & \z_{1,5}-a_2^-(r_2^-)_1 \\
\z_{2,5}' & = & \beta_1\z_{1,5}+\cdots +\beta_3\z_{3,5}-a_2^-(r_2^-)_2 \\
\z_{3,5}' & = & \eta_1\z_{1,5}+\cdots +\eta_3\z_{3,5}+qk(-\z_{4,5}'+\theta_1\z_{1,5}+\cdots +\theta_4\z_{4,5})-a_2^-(r_2^-)_3. 
\end{eqnarray*}
Finally $\tilde{\z}$ satisfies
\begin{eqnarray*}
0 & = & \tilde{\z}_1+[\rho] \\
\tilde{\z}_2' & = & \beta_1\tilde{\z}_1+\cdots +\beta_{3}\tilde{\z}_3 +[m] \\
\tilde{\z}_3' & = & \eta_1\tilde{\z}_1+\cdots + \eta_{3}\tilde{\z}_3+qk(-\tilde{\z}_4'+\theta_1\tilde{\z}_1+\cdots+\theta_4\tilde{\z}_4)+[\E]+q
\end{eqnarray*}
These equations allow us to perform row operations simplifying rows 1, 5, and 6. Then rearranging rows and columns we find that 
\begin{equation*}
D'(0)=\det\left(\begin{array}{c|c} \begin{array}{lccr} 0 & \cdots  & \cdots &0 \\  0 &\cdots & \cdots &0 \\ 0 & \cdots &  \cdots &0 \end{array} & \begin{array}{ccc} a_1^-r_1^- & a_2^-r_2^- & [U]+\vec{q} \end{array} \\ \hline  \begin{array}{cccc} \z_{2,1}^+ & \z_{2,2}^+ & \z_{2,3}^+ & \z_{2,6}^- \\  \z_{3,1}^+ & \z_{3,2}^+ & \z_{3,3}^+ & \z_{3,6}^- \\ \z_{4,1}^+ & \z_{4,2}^+ & \z_{4,3}^+ & \z_{4,6}^- \\ \z_{4,1}^{+'} & \z_{4,2}^{+'} & \z_{4,3}^{+'} & \z_{4,6}^{-'} \end{array} & \begin{array}{lcr} * & * & * \\ \vdots & \vdots & \vdots \\ \vdots & \vdots & \vdots \\ * & * & * \end{array} \end{array}\right)
\end{equation*}
from which the result follows.
\foorp
%
%
%
%
%
\subsubsection{Large \lam\ Behavior} 
To finish the calculation of the stability index, it remains to determine
$$\sgn\evans\quad\mbox{as}\quad \lam\to +\infty ,\mbox{real}.$$
Here block triangular structure plays a key role. In particular it allows for certain ``diagonal'' kinematic and reaction components of eigenvectors to be separately analytically specifiable. This means that we will be able to calculate the sign in two pieces. For the gas dynamical piece the analysis of Appendix C applies, while the reaction piece can be treated directly.

For ease of comparison and consistency with our previous analyses, we shift now to the notation of Appendix C. Thus, $(u,v)$ represent the gas dynamical variables, $\tilde{A}=A_{11}-A_{12}b_2^{-1}b_1$ where $A_{ij}$ is the $ij$-entry of the flux Jacobian and $b_1$, $b_2$ are the nonzero blocks of the matrix $B$. Recall that in this case $\tilde{A}$ is simply the particle velocity ($u$ in the usual notation). For 3-shocks, this quantity is uniformly negative.  We also use $z$ to represent the reaction variable.
\begin{prop}For real \lam\ sufficiently large,
\begin{equation}
\sgn D(\lam) =\sgn \underbrace{\det(\mathbb{S}^+)\det(\pi\mathbb{W}^+,\epsilon\mathbb{S}^+)\det(\pi\mathbb{W}^-)}_{\text{kinematic only}\;j=1,2,4,5,6}\times\underbrace{z_3^+z_7^-}_{\text{reaction}}\ne 0
\label{detlargelamnotzero}
\end{equation} 
where $\pi$ denotes projection of the kinematic variable $W=(u,v,v')$ onto $(u,v)$ components, and
$\mathbb{S}(x)$ is a real basis of the stable subspace
of $\tilde{A}$, and $\varepsilon u:=(u, -b_2^{-1}
b_1u)$ denotes extension.
\end{prop}
%
%
%
\pf For real \lam\ sufficiently large, from the block triangular structure in \eqref{hypscaling} and \eqref{parscaling}, it follows that at any (fixed) $x$ the stable/unstable manifolds of the frozen eigenvalue equation are spanned by vectors of the forms
\begin{equation}
\begin{pmatrix} u \\ -b_2^{-1}b_1u \\ * \\ * \\ * \end{pmatrix},\quad\begin{pmatrix} 0 \\ v \\ \mp\tilde{\mu}\lam^{1/2}v \\ * \\ *\end{pmatrix}\;\;\text{(kinematic)}
\end{equation} 
and 
\begin{equation}
\begin{pmatrix}0 \\ 0 \\ 0 \\ 0 \\ 1 \\ \pm\lam^{1/2}d^{-1/2} \end{pmatrix}\;\;\text{(reaction)}.
\end{equation}
Therefore the Evans function satisfies
\begin{multline*}
D(\lam)	  \sim \det\left\{ \left(\begin{array}{ccc} 0 & 0 & 0 \\ v_1 & v_2 & 0 \\ -\tilde{\mu}_1\lam^{1/2}v_1 & -\tilde{\mu}_2\lam^{1/2}v_2 & 0 \\ * & * & 1 \\ * & * & -\lam^{1/2}d^{-1/2} \end{array}\right)\underbrace{\left(\begin{array}{c|c} \beta_1^+ & 0 \\ \hline * & \alpha_z^+ \end{array}\right)}_{\bar{\alpha}^+}, \right . \\  
\left.\left(\begin{array}{cccc} u & 0 & 0 & 0 \\ -b_2^{-1}b_1u & v_1 & v_2 & 0 \\ * & \tilde{\mu}_1\lam^{1/2}v_1 & \tilde{\mu}_2\lam^{1/2}v_2 & 0 \\ * & * & * & 1 \\ * &*&*& \lam^{1/2}d^{-1/2}\end{array}\right)\underbrace{\left(\begin{array}{c|c} \begin{array}{cc} \alpha_1^- & \alpha_2^- \\ \beta_2^- & \beta_1^- \end{array} & 0 \\ \hline * & \alpha_z^- \end{array}\right)}_{\bar{\alpha}^-} \right\}
\end{multline*}
where $\alpha_j^{\pm}$ and $\beta_j^{\pm}$ are as in the gas dynamics case, so
\begin{equation}
\det\left( \beta_1^+ \right) \quad\mbox{and}\quad
\det\left(\begin{array}{cc} \alpha_1^- & \alpha_2^- \\ \beta_2^- & \beta_1^- \end{array}\right)
\end{equation}
are real nonzero quantities while $\alpha_{z}^{\pm}$ are real, nonvanishing scalar functions of $x$. Note that block triangular structure plays a key role here as it allows the matrices $\bar{\alpha}^{\pm}$ to be block triangular. The RHS is equal to the product of 
\begin{equation}
\det \left(\begin{array}{c|c} V_1 & V_2\end{array}\right) \label{bigmatrix}
\end{equation}
where
\begin{equation}
V_1=\begin{pmatrix}  0 & 0  & 0 \\  v_1 & v_2 & 0 \\ -\tilde{\mu}_1\lam^{1/2}v_1 & -\tilde{\mu}_2\lam^{1/2}v_2 & 0 \\ * & * & 1 \\ 
* & * & -\lam^{1/2}d^{-1/2} \end{pmatrix},
\end{equation} 
\begin{equation}
V_2=\begin{pmatrix}  u & 0 & 0 & 0 \\  -b_2^{-1}b_1u & v_1 & v_2 & 0 \\  * & \tilde{\mu}_1\lam^{1/2}v_1 & \tilde{\mu}_2\lam^{1/2}v_2 & 0 \\ * & * &* & \lam^{1/2}d^{-1/2} \end{pmatrix},
\end{equation}
and
\begin{equation}
\det\begin{pmatrix} \bar{\alpha}^+ & 0 \\ 0 & \bar{\alpha}^- \end{pmatrix} \label{littlematrix}
\end{equation}
Swapping columns, it's clear that the determinant \eqref{bigmatrix} simplifies to
\begin{multline}
-\det\begin{pmatrix}  0 & 0 & u & 0 & 0 \\ v_1 & v_2 & -b_2^{-1}b_1u & v_1 & v_2 \\ 
-\tilde{\mu}_1\lam^{1/2}v_1 & -\tilde{\mu}_2\lam^{1/2}v_2 & * & \tilde{\mu}_1\lam^{1/2}v_1 & \tilde{\mu}_2\lam^{1/2}v_2
\end{pmatrix}\times \\
\det\begin{pmatrix}  1 & 1 \\ -\left(\frac{\lam}{d}\right)^{1/2} & \left(\frac{\lam}{d}\right)^{1/2} \end{pmatrix}\label{swapped}
\end{multline}
which is broken into kinematic and reaction pieces. The kinematic part then reduces as in Lemma \ref{Clemma1} to 
$$\det(\mathbb{S}^+)\det(\pi\mathbb{W}^+,\varepsilon\mathbb{S}^+)\det(\pi\mathbb{W}^-),$$ while 
the sign of the reaction part is (where we have included the minus sign from \eqref{swapped}) 
$$-z_3^+z_7^-$$ whence the result follows. Nonvanishing follows from corresponding result for nonreacting gas dynamics and the fact the $\alpha_z^{\pm}$ are nonvanishing for all $x$.
\foorp
Nonvanishing of $D(+\infty)$ yields that the stability index $\tilde{\Gamma}$ satisfies either
$$\tilde{\Gamma}=\sgn \gamma_d\Delta,$$
or
$$\tilde{\Gamma}=-\sgn\gamma_d\Delta,$$
as model parameters are modeled smoothly. Thus we find that the \emph{relative stability index}, defined to be 
$$\sgn\gamma_d\Delta,$$
gives a measure of \emph{spectral flow}. That is, changes in the sign of the relative stability index indicate a change in stability. 

In this case, however, we can do more and actually evaluate the absolute stability index,
$$\sgn D'(0)D(+\infty),$$
by relating our formula above for large \lam\ to the normalizations we've chosen at $\lam=0$. 
\begin{prop} For \lam\ real and sufficiently large, 
\begin{equation}
\sgn D(\lam) =-\sgn \det(\mathbb{S}^+)\det(\pi\mathbb{W}^+,\varepsilon\mathbb{S}^+)\det(\pi\mathbb{W}^-)z_3^+z_7^-|_{\lam=0}\ne 0 \label{pullbackform}
\end{equation}
\label{dettolamzero}
\end{prop}
\pf The form of \eqref{detlargelamnotzero} allows us to connect to $\lam=0$ separately in kinematic and reaction terms. Since the pair $(A,B)$ satisfy the semidissipativity conditions, the result of the Appendix holds. As for the reaction term, pull back to $\lam =0$ is evident since ``reaction'' vectors have the form
$$ (*,*,*,*,*,\underbrace{1}_z,*),$$
hence projection onto the $z$ component never vanishes, and thus cannot change sign.
\foorp

%
%
%
\subsection{The Stability Index}
Combining Propositions \ref{detdprime} and \ref{dettolamzero} we obtain the result.
\begin{thm}The stability index for a strong detonation with Lax 3-shock structure is 
$$\tilde{\Gamma}:=-\sgn \gamma_d\Delta \gamma_{\text{NS}}z_3^+z_7^- \det(r_1^-,r_2^-,U_6^-),$$
where $\gamma_{\text{NS}}$ involves only the gas dynamics components.
\end{thm}
\pf The stability index has the form $\tilde{\Gamma}=\sgn D'(0)D(+\infty)$. From Proposition~\ref{detdprime} we find that $D'(0)=\gamma_d\Delta$. On the other hand from Proposition~\ref{dettolamzero}, the term $\det(\pi\mathbb{W}^-)$ simplifies to $\det(r_1^-,r_2^-,U_6^-)$ as in Appendix C, while the term $\det(\pi\mathbb{W}^+,\varepsilon\mathbb{S}^+)$ breaks into the product of
$$\gamma_{\text{NS}}=\det\begin{pmatrix} \z_{2,1}^+ & \z_{2,2}^+ \\ \z_{3,1}^+ & \z_{3,2}^+ \end{pmatrix} $$
and $\det(\mathbb{S}^+)$ the latter of which cancels with the factor $\det(\mathbb{S}^+)$ appearing in \eqref{pullbackform}.   \foorp
%
%
%
%
%
%
%
\subsection{Reduction of $\Delta$ and Nonvanishing of $\Delta$ for Ideal Gas}

The term $\Delta$ in the stability index has the form of a modified Lopatinski determinant
$$\Delta =\det (r_1^-,r_2^-,[U]+\vec{q}).$$
In place of the determinant above, it is more convenient for calculations to rewrite the quantity as the dot product of an appropriate left eigenvector and the jumps in the conserved quantities. Thus 
$$\Delta =l_3^-\cdot ([U]+\vec{q})=l_3^-\cdot [U] +l_3^-\cdot\vec{q}.$$
The left eigenvector $l_3^-$ is
$$l_3^-=\left(p_{\rho}-\frac{p_e}{\rho}e-cu+\frac{p_eu^2}{2\rho},c-\frac{p_e}{\rho}u,\frac{p_e}{\rho}\right),$$
while
$$[U]=([\rho],0,[\mathcal{E}])^{tr}$$ 
and $\vec{q}=(0,0,q)^{tr}$. It is then straightforward to calculate 
\begin{equation}
\Delta =\left(p_{\rho}-\frac{cm}{\rho}\right)[\rho]+\frac{p_em^2}{2\rho}\left(\frac{[\rho]}{\rho^2}+[1/\rho]\right)+\frac{p_e}{\rho}\rho_+[e]+\frac{p_e}{\rho}q
\end{equation}
Now using 
$$\frac{[\rho]}{\rho_-^2}+[1/\rho]=-\frac{1}{\rho_-}[\rho][1/\rho]$$
in our expression for $\Delta$ and simplifying, we obtain
\begin{equation}
\Delta=\left(-\frac{p_em^2}{2\rho^2}[1/\rho]+p_{\rho}-\frac{cm}{\rho}\right)[\rho]+\frac{p_e}{\rho}\rho_+[e]+\frac{p_e}{\rho}q\end{equation}
%
%
%
Similarly as in \cite{serretransition}, we have
\begin{claim}$[e]+\langle p\rangle [1/\rho]=-\frac{qu_+}{m}$ where $\langle p\rangle =\frac{1}{2}(p_++p_-)$. \label{claim2} \end{claim}
\pf From \eqref{RH3}
$$m\left[\frac{\mathcal{E}}{\rho}\right]+m\left[\frac{p}{\rho}\right] =-u_+q,$$
or 
\begin{equation} [E]+\left[\frac{p}{\rho}\right]=-\frac{u_+q}{m}. \label{solvefore} \end{equation}
Using the relationship $E=\frac{u^2}{2}+e$, we can use (\ref{solvefore}) to calculate an expression for $[e]$, the jump in the specific internal energy. Thus we find
\begin{equation} [e]=-\left(\left[\frac{p}{\rho}\right]+\frac{m^2}{2}[\rho^{-2}]\right)-\frac{u_+q}{m} \label{ejump} \end{equation}
The term in parentheses in (\ref{ejump}) can be simplified to 
\begin{equation}
\left[\frac{p}{\rho}\right]+\frac{m^2}{2}[\rho^{-2}]= \langle p\rangle [1/\rho]. 
\end{equation} 
The claim then follows by substitution into (\ref{ejump}). \foorp

From Claim \ref{claim2}, we find that 
$$\frac{p_e}{\rho}\rho_+[e]=-\frac{p_e}{\rho}\rho_+\langle p\rangle [1/\rho]-\frac{p_e}{\rho}q,$$
and upon substituting into the expression for $\Delta$, we find that the $q$ terms cancel out. Therefore
\begin{equation}
\Delta = [\rho]\left(-\frac{p_em^2}{2\rho^2}[1/\rho]+p_{\rho}-\frac{cm}{\rho}+\frac{p_e}{\rho^2}\langle p\rangle\right)
\end{equation}
Since $[\rho]\neq 0$, the condition $\Delta =0$ may be written as
\begin{equation}
p_{\rho}-\frac{cm}{\rho}+\frac{p_e}{\rho^2}\left(\frac{[p]}{2}+\langle p\rangle\right) = 0 \label{delta01}
\end{equation}
Next we note that 
$$\left(\frac{[p]}{2}+\langle p\rangle\right)=p_+$$
and thus (\ref{delta01}) simplifies to
\begin{equation}
p_{\rho}-\frac{cm}{\rho}+\frac{p_ep_+}{\rho^2} = 0. \label{delta02}
\end{equation}
Finally using the fact that the sound speed $c$ satisfies $c^2=p_{\rho}+\rho^{-2}pp_e$, we find that 
$$p_{\rho}+\frac{p_ep_+}{\rho^2}=c^2+\frac{p_e}{\rho^2}[p],$$
and thus we reduce (\ref{delta02})
\begin{eqnarray*}
c^2+\frac{p_e}{\rho^2}[p]-\frac{mc}{\rho} & = & 0 \\
c^2-\frac{p_em^2}{\rho^2}[1/\rho]-\frac{mc}{\rho} & = & 0\;\; \mbox{by \eqref{RH2}} \\
1-\frac{p_em^2}{\rho^2c^2}[1/\rho]-\frac{m}{\rho c} & = & 0 
\end{eqnarray*}
Since the Mach number $M$ satisfies $M=-\frac{m}{\rho c}$, we obtain Majda's condition for inviscid shock instability
\begin{equation}
M^2[1/\rho]p_e-M-1=0. \label{majdacond}
\end{equation} 
See \cite{MBOOK}. We remark here that the jumps in the detonation formula above refer to endstates of the whole wave and not of the Neumann shock at the leading edge.
%
%
%
%
%
%
%
%
%
Using that $\Delta$ reduces to the Majda condition for the stability of inviscid shocks, we now show that $\Delta$ does not vanish when the equation of state is assumed to be of the form
$$p(\rho,e)=\Gamma\rho e.$$
We denote the compression ratio by
$$r=\frac{\rho_+}{\rho_-},$$
and we denote by $r^*$ the compression ratio of the Neumann shock. It follows immediately from the Rankine-Hugoniot diagram, Figure \ref{hugoniot}, that
\begin{equation}1<r<r^*.\label{compratio} \end{equation} 
Also in the case of an ideal gas, we find similarly as \cite{serretransition}
\begin{equation}1<r^*<\frac{\gamma +1}{\gamma -1}=1+\frac{2}{\Gamma}\label{idealgamma} \end{equation} 
Combining (\ref{compratio}) and (\ref{idealgamma}) we find that
\begin{equation} 0<r-1<\frac{2}{\Gamma} \end{equation}
Specializing the Majda condition to the ideal gas case yields, 
\begin{eqnarray*}
\Gamma (r-1)M^2-M-1 & < & 2M^2-M-1 \\
	& = & (2M+1)(M-1). 
\end{eqnarray*}
The quantity $(2M+1)(M-1)$ is nonpositive for $-\frac{1}{2}\leq M\leq 1$ and for strong detonations the Mach number $M$ satisfies $0<M < 1$. Finally since $\Delta=[\rho](\Gamma (r-1)M^2-M-1)$ by our reduction above and since $[\rho]<0$ for 3-shock detonations, we conclude that $\Delta>0$ for all strong detonation waves under the ideal gas assumption.


\subsection{Evaluation of $\tilde{\Gamma}$ in the Ideal Gas ZND limit}
Here we use the structure of the singular manifolds contructed in the existence argument of \cite{GS1} to determine the sign of the transversality coefficient $\gamma_d$ and the other terms in the stability index. Recall that this coefficient is defined by the determinant

\begin{equation*}\gamma_d=\det\left(\begin{array}{cccc} \z_{2,1}^+ & \z_{2,2}^+ & \z_{2,3}^+ & \z_{2,6}^- \\  \z_{3,1}^+ & \z_{3,2}^+ & \z_{3,3}^+ & \z_{3,6}^- \\ \z_{4,1}^+ & \z_{4,2}^+ & \z_{4,3}^+ & \z_{4,6}^- \\ \z_{4,1}^{+'} & \z_{4,2}^{+'} & \z_{4,3}^{+'} & \z_{4,6}^{-'} \end{array}\right).\end{equation*}
In order to use the structure of the singular manifolds, we must translate from the $\zeta$-coordinates of our stability analysis to the $uTYZ$-coordinates of the existence argument. This is accomplished in two steps. First, using the original  $\z$-coordinate transformation \eqref{zetacoord} and the $\lam =0$ eigenvalue equation, we can connect the $\z$-coordinates to a set of intermediate coordinates: $\rho_j,\E_j,z_j,z_j'$ via the linear transformation (dropping $\pm$) 

\begin{equation}\begin{pmatrix} \z_{2,j} \\ \z_{3,j} \\ \z_{4,j} \\ \z_{4,j}' \end{pmatrix} =\begin{pmatrix} b_{21} & 0 & 0 & 0 \\ b_{31} & b_{33} & 0 & 0 \\ \tilde{d} & 0 & d & 0 \\ \tilde{d}'+b_{21}^{-1}\alpha_{23}\tilde{d} & b_{21}^{-1}\alpha_{23}\tilde{d}& 0 & d \end{pmatrix}\begin{pmatrix} \rho_j \\ \E_j \\ z_j \\ z_j' \end{pmatrix}. \label{trany1}\end{equation}

Note that the determinant of the linear transformation above is simply

$$b_{21}b_{33}d^2=-\frac{\nu m}{\rho^3}\theta c_v^{-1}d^2>0.$$

To obtain \eqref{trany1} note that from \eqref{zetacoord} we know
\begin{align}
\z_{2,j}& = b_{21}\rho_j \label{zeta2} \\
\z_{3,j}& = b_{31}\rho_j +b_{33}\E_j  \label{zeta3} \\
\z_{4,j}& = \tilde{d}\rho_j + dz_j. \label{zeta4}
\end{align}
Also from \eqref{zeta4} it follows that 
$$\z_{4,j}'=\tilde{d}'\rho_j +\tilde{d}\rho_j'+dz_j'.$$
We take advantage of the second eigenvalue equation at $\lam =0$,
$$(\alpha_{21}\rho_j+\alpha_{23}\E_j)' = (b_{21}\rho_j')',$$
and the fact that we are interested in the behavior of fast modes, so we can integrate up once to obtain
$$\rho_j' = b_{21}^{-1}\alpha_{21}\rho_j+b_{21}^{-1}\alpha_{23}\E_j,$$
which we can use to write $\rho_j'$ in terms of $\rho_j$ and $\E_j$. Substituting into the equation for $\z_{4,j}'$, we get
\begin{equation}
\z_{4,j}'= (\tilde{d}'+b_{21}^{-1}\alpha_{21}\tilde{d})\rho_j+b_{21}^{-1}\alpha_{23}\tilde{d}\E_j +dz_j'. \label{zeta4prime}
\end{equation}
Combining equations \eqref{zeta2}-\eqref{zeta4} and \eqref{zeta4prime} yields \eqref{trany1}.

The second step is to connect $\rho_j,\E_j, z_j, z_j'$, the variations in the conserved quantities, to $u_j,T_j,Y_j,Z_j$, the coordinates in which the construction of the singular manifolds has been accomplished. We use the relationships
\begin{align}
\rho & u = m \\
\rho & (u^2/2+c_vT) = \E \\ 
\rho & Y = z, 
\end{align}
and we linearize about the profile to obtain (note that subscript $j$ indicates a variation while a bar indicates that a quantity is evaluated on the profile)
\begin{align}
\rho_j & = -\bar{\rho}\bar{u}^{-1}u_j  \label{rho} \\
\E_j & = mu_j+c_v\bar{\rho}T_j-\bar{\rho}\bar{u}^{-1}u_j(\bar{u}^2/2+c_v\bar{T}) \label{E} \\
z_j & = -\bar{\rho}\bar{u}^{-1}\bar{Y}u_j+\bar{\rho}Y_j. \label{z}
\end{align} 
Also since $z'=\rho'Y+\rho Y'$, we linearize  and use known relationships among the variations to obtain 
\begin{equation}
z_j'=\bar{\rho}'Y_j+\bar{Y}(b_{21}^{-1}\alpha_{21}\rho_j+b_{21}^{-1}\alpha_{23}\E_j)+\bar{Y}'\rho_j+\bar{\rho}(-\frac{m}{\bar{\rho}d}Z_j+\frac{m}{\bar{\rho}d}Y_j), \label{zjprime}
\end{equation}
where $Z$ is defined as in \cite{GS1} by 
$$Z=Y-\frac{\rho d}{m}Y_x.$$
Finally we can use \eqref{rho}-\eqref{z} and \eqref{zjprime} to write down the coordinate change as
\begin{equation}\begin{pmatrix} \rho_j \\ \E_j \\ z_j \\ z_j' \end{pmatrix} =\begin{pmatrix} -\rho u^{-1} & 0 & 0 & 0 \\ m-\rho u^{-1}(u^2/2+c_vT) & c_v\rho & 0 & 0 \\ -\rho u^{-1}Y & 0 & \rho & 0 \\ M_1 & M_2 & \rho'+\frac{m}{d} & -\frac{m}{d} \end{pmatrix}\begin{pmatrix} u_j \\ T_j \\ Y_j \\ Z_j \end{pmatrix}, \label{trany2}\end{equation}
where
\begin{align}
M_1& =   -\rho u^{-1}(b_{21}^{-1}\alpha_{21}Y+Y')+(m-\rho u^{-1}(u^2/2+c_vT))b_{21}^{-1}\alpha_{23}Y \\
M_2& =  b_{21}^{-1}\alpha_{23}Yc_v\rho
\end{align}
and we have dropped the bars on all terms in the matrix. Moreover, we note that the determinant of the matrix of this second coordinate change is
$$c_v\rho^3u^{-1}\frac{m}{d}>0.$$
Therefore we find, by virtue of the fact that the determinants of the two change of coordinates matrices are positive, that 
\begin{equation}\sgn\gamma_d=\sgn\det\begin{pmatrix} u_1^+ & \cdots & u_6^- \\ T_1^+ & \cdots & T_6^- \\ Y_1^+ & \cdots & Y_6^- \\ Z_1^+ & \cdots & Z_6^- \end{pmatrix},\end{equation}
By a simpler calculation proceeding as above, we also find that 
\begin{equation}
\sgn\gamma_{\text{NS}}=\sgn\det\begin{pmatrix} \z_{2,1}^+ & \z_{2,2}^+ \\ \z_{3,1}^+ & \z_{3,2}^+ \end{pmatrix}=\sgn\det\begin{pmatrix}u_1^+ & u_2^+ \\ T_1^+ & T_2^+ \end{pmatrix}.
\end{equation}
Finally the sign of $\tilde{\Gamma}$ can be determined by a careful examination of the structure of the singular manifolds from which the solution is constructed.

Recall from Chapter~1 that the structure of the singular manifolds is as in Figure~\ref{singularflow}. In this figure, the $T$-axis is perpendicular to the $uZ$-plane on the page, while the fourth missing direction in the phase space is the $Y$ direction. Also the parabolic curve, $\mathcal{C}$, on which the slow (reactive) flow takes place is not in a $T=\text{constant}$ plane, see Figure~\ref{manifold}. The diagonal dotted line indicates the intersection of the $T=T_i$ plane and $\mathcal{C}$. Thus the bold segments of the branches of $\mathcal{C}$ are below ignition temperature and there are no slow dynamics on those portions of the curve. Furthermore, we note that the unburned state $(u_+,T_+,Y_+,Z_+)$ is a degenerate rest point due to the ignition temperature assumption with a 3-dimensional stable manifold. On the other hand the burned state, $(u_-,T_-,0,0)$, has a 2-dimensional unstable manifold (featuring a reactive and a kinematic direction) and a 2-dimensional stable manifold. Note that the $T$ and $Y$ directions are both \emph{stable} directions, see the discussion of existence in Chapter~1 particularly equations~\eqref{redmom}-\eqref{redprog} and equations~\eqref{laysingmom}-\eqref{laysingextra}.

\begin{figure}[ht]
\centerline{\epsfxsize6truecm\epsfbox{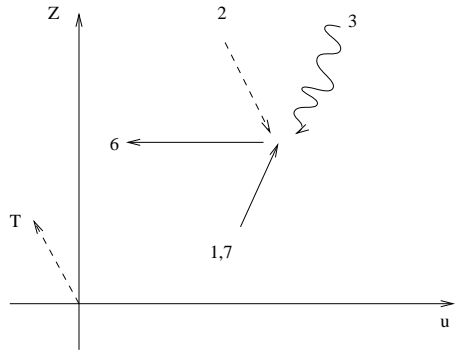}}
\caption{The Point of Evaluation}
\label{eval}
\end{figure}

We evaluate 
\begin{equation}
-\sgn\gamma_d\gamma_{\text{NS}}\det(r_1^-,r_2^-,U_6^-)z_3^+\bar{z}_x|_{-\infty}
\label{beast}
\end{equation}
at the ``corner'' where the fast manifold which approaches the burned state intersects with the opposite branch of $\mathcal{C}$. In figure~\ref{eval} we see a schematic indicating the relevant vectors in the calculation. The $1,7$ arrow corresponds to the profile, while the $6$ arrow corresponds to the kinematic direction at $-\infty$. The dashed $2$ arrow represents the stable manifold in the $T$-direction, while the curvy $3$ line corresponds to the missing stable $Y$-direction.  
Then we find that $\sgn\gamma_d\gamma_{\text{NS}}$ can be computed as
\begin{equation}
\sgn\det\begin{pmatrix} + & 0 & 0 & - \\ - & - & 0 & - \\ * & 0 & - & 0 \\ + & 0 & 0 & 0 \end{pmatrix}\det\begin{pmatrix} + & 0 \\ - & - \end{pmatrix}=-1.
\label{sgngdgns}
\end{equation}
These choices force 
\begin{equation}
\sgn\det(r_1^-,r_2^-,U_6^-)=-1,
\label{sgnnsdet}
\end{equation}
To see \eqref{sgnnsdet}, note that 
$\det(r_1^-,r_2^-,U_6^-|_{-\infty})=l_3^-\cdot U_6^-$ where $l_3^-$ is the appropriate left-eigenvect
or, (simplified due to ideal gas assumption)
\begin{equation}
l_3^-=\left(-cu+\frac{\Gamma u^2}{2},c-\Gamma u,\Gamma\right)^{tr}.
\label{idealleft}
\end{equation}
Then since 
\begin{equation}
U_6^-=(\rho_6^-,m_6^-,\E_6^-)^{tr}=(\rho_6^-,0,\E_6^-)^{tr},
\label{nsprofile}
\end{equation}
we obtain by combining \eqref{idealleft} and \eqref{nsprofile} that 
\begin{equation}
l_3^-\cdot U_6^-=(-cu+\frac{\Gamma u^2}{2})\rho_6^-+\Gamma\E_6^-.
\label{firstdot}
\end{equation}
To take advantage of the signs we know, we translate into $u,T$ coordinates. 
We use \eqref{trany1} and \eqref{trany2} in \eqref{firstdot}, and then some elementary simplification yields
\begin{equation}
l_3^-\cdot U_6^-=\left(c\rho-\frac{c_vp}{u}\right)u_6^-+\Gamma c_v\rho T_6^-.
\label{finaldot}
\end{equation}
To evaluate the sign of \eqref{finaldot}, we note that because we consider a 3-shock, $u<0$, and thus the coefficient
$$\left(c\rho-\frac{c_vp}{u}\right)$$
is positive. Finally since both $\sgn u_6^-$ and $\sgn T_6^-$ are $-1$, \eqref{sgnnsdet} follows.

Also we find then that 
\begin{equation}
\sgn z_3^+=\sgn \bar{\rho}Y_3^+=-1.
\label{plusreactive}
\end{equation}
Lastly, we find that 
\begin{equation} 
\sgn z_7^-=\sgn\bar{z}_x|_{-\infty}=\sgn (\bar{\rho}\bar{Y}_x|_{-\infty}-\bar{\rho}\bar{u}^{-1}\bar{Y}\bar{u}_x|_{-\infty})=+1.
\label{profile}
\end{equation} 
Tracking the signs computed in \eqref{sgngdgns}-\eqref{profile} and combining with \eqref{beast} and the fact that we computed 
$$\Delta >0,$$ 
for all ideal gas strong detonations, we discover that in the ideal gas ZND limit, the stability index satisfies
\begin{eqnarray*}
\tilde{\Gamma} &= & -\sgn\gamma_d\gamma_{\text{NS}}\Delta\det(r_1^-,r_2^-,U_6^-)z_3^+\bar{z}_x|_{-\infty}\\
& =&(-1)(-1)(+1)(-1)(-1)(+1)=+1
\label{beast2}
\end{eqnarray*}
which is consistent with stability. This completes Theorem~\ref{result1}
	
Finally, we also remark here that under the ideal gas assumption we found $\sgn\Delta$ by comparison to the Neumann Shock. Restating this, we found that if the Neumann Shock is ``index stable,'' \emph{i.e.} satisfies the stability index necessary criterion, then the corresponding strong detonation is also ``index stable.'' A natural question then is: under what conditions on the equation of state does this remain true? Or perhaps more importantly, what is the actual stability relationship, not just the relationship between stability indexes, between the shock and the strong detonation? The partial information gathered from the stability index approach definitely motivates further investigation into this question.

\subsection{Stability Index for Multiple Reactants}

In actual combustion, the chemical reactions involved are typically more complicated than a single one-step reaction. It's natural to model reactions in an $s+1$ component gas by a system of the general form
\begin{eqnarray}
U_t+f(U)_x & = & (B(U)U_x)_x+Q\Phi (U)z \label{generalkinematic} \\
z_t+(v(U)z)_x & = & (D^1(U)z_x)_x+(D^2(U,z)U_x)_x-\Phi (U)z, \label{generalreaction}
\end{eqnarray}
where $U=(\rho,m,\E)^{tr}$ and $f$ and $B$ are as in the Navier-Stokes model for gas dynamics considered above. The vector $z=\rho Y\in\R^s$ measures the quantities of each of the reactants; the constant matrix $Q\in \R^{3\times s}$ records the heat released in each reaction, thus 
$$ Q=\left(\begin{array}{ccc} 0 & \cdots & 0 \\ 0 & \cdots & 0 \\ q_1 & \cdots & q_s \end{array}\right); $$
The positive definite matrix $\Phi(U)\in\R^{s\times s}$ incorporates the reaction rates; the diffusion matrices satisfy $D^1\in\R^{s\times s}$ and $D^2\in\R^{s\times 3}$; and the scalar function $v(U)$ is simply $\frac{m}{\rho}=u$, the velocity. 

It turns out that the number of equations required to model the chemistry can be reduced through the use of progress variables. See \cite{FD} or \cite{Lyng} for further discussion. Our analysis largely applies to the more general multiple reactant case with the following important exception below.

In the more complicated multi-species case the large \lam\ calculation for the reaction block is not so straightforward. In particular the connection to $\lam=0$ is not clear; we indicate here some partial results along this line. Examining the reaction portion of the characteristic equation on the $+\infty$ side, we find since the reaction function satisfies $\Phi =0$,
\begin{equation}
(\mu^2D-\mu V+\lam I)z= 0.
\end{equation}
Since the convection in the reaction equation is just the background velocity, $V$ is scalar, hence it commutes with the species diffusion. This implies that there is a symmetrizer $V^0$ such that $V^0V$ and $V^0D$ are symmetric, and $\re V^0D>0$. Thus arguments using Theorem \ref{serrelemma} apply.  On the other hand, things are trickier on the $-\infty$ side the characteristic equation reads
\begin{equation}
(\mu^2D-\mu V+\lam I-\Phi)z=0.
\end{equation}
Here, we point out that the argument of \cite{BSZ}, Lemma 7.2 goes through word for word for constant-coefficient eqautions of the form $\lam z+Az'+Bz''+Cz=0$ provided that $A$ is symmetric and $B$ and $C$ are positive definite, i.e. $\re B,\re C>0$, (more generally there is a coordinate transformation for which this is true) with the slightly modified Lyapumov function $H(v,w)=(1/2)\ip{Aw}{w}+\re\ip{x}{(\lam + C)v}$. Indeed the conclusions become somewhat stronger extending to all $\re\lam\ge 0$ and not only $\re\lam >0$ as in \cite{BSZ}. In our case $A=v$ is scalar, $B=D$ and $C=\Phi$, so we find that the machinery of \cite{BSZ} may be applied provided that there exist coordinates in which $D$ and $\Phi$ are both positive: in particular if $D$ is scalar.
%
%
%
%
%
\subsection{Small $q$ Stability}\label{smallq} Here we examine the $q\to 0$ limit, and take advantage of the simplification in the equations when $q=0$. Using a continuity argument, we prove \emph{strong spectral stability}, that is, nonexistence of eigenvalues with $\re\lam\ge 0$ and $\lam\ne 0$, as well as transversality of the connecting profile $\gamma_d\ne 0$, and low-frequency stability $\Delta\ne 0$. It is expected that by following the program of \cite{ZH} and \cite{MZ1,MZ2,MZ3}, these three properties  should be sufficient to conclude full nonlinear orbital stability. Such a result would be an extension to the reacting Navier-Stokes model of the results of  \cite{LYi},\cite{LLT} in which the authors show nonlinear stability of strong detonations in the Majda model as $q\to 0$. We note that in both \cite{LYi} and \cite{LLT}, as here, working with the integrated equations is a key ingredient in the analysis.

We rewrite the system \eqref{eq:mass}-\eqref{eq:progress} using $U$ to represent the vector of gas dynamical variables and $z=\rho Y$ in the first three equations, 
\begin{align}
&(U+\vec{q}z)_t+(f(U)+\vec{q}uz)_x  =  (B(U)U_x)_x, \label{gas} \\
&(\rho Y)_t+(\rho uY)_x  =  (\rho dY_x)_x-k\varphi(T)\rho Y. \label{react}
\end{align}
Now we define $\tilde U$ by $\tilde{U}=\int_{-\infty}^{x}(U+\vec{q}z)$, and consider the linearized eigenvalue equation  from \eqref{gas} in the integrated variable $\tilde{U}$,
\begin{equation}
\lam \tilde{U}+A\tilde{U}'-(A-u)\vec{q}z=B\tilde{U}'',\label{intgaseval}
\end{equation}
and from \eqref{react}
\begin{multline}
\lam(\bar{\rho}Y+\rho\bar{Y})+(\bar{\rho}\bar{u}Y+\bar{\rho}u\bar{Y}+\rho\bar{u}\bar{Y})'=(\bar{\rho}dY'+\rho d\bar{Y}_x)' \\-k\bar{\varphi}\rho\bar{Y}-k\bar{\varphi}\bar{\rho}Y-k\bar{\varphi}'(\bar{T})T\bar{\rho}\bar{Y}. \label{ylinreact}
\end{multline}
Note that when $q=0$, \eqref{intgaseval} reduces to the integrated eigenvalue equation for gas dynamics \emph{about a gas dynamical profile}.

\begin{lemma}
The limiting system (in integrated form) has no eigenvalues with $\re\lam\ge 0$ provided the limiting shock is stable. (In particular, the limiting shock is stable if the amplitude is sufficiently small.)
\end{lemma}
\pf Since the limiting shock is spectrally stable, we have by standard considerations \cite{ZH} that the integrated eigenvalue equation for gas dynamics supports no eigenvalues on $\re\lam\ge 0$, and so $(\rho,u,\E)$ identically vanish for any eigenfunction of the limiting eigenvalue equations as $q\to 0$. Thus \eqref{ylinreact} has the simpler form 
$$\lam(\bar{\rho}Y)+(\bar{\rho}\bar{u}Y)'=(\bar{\rho}dY')'-k\bar{\varphi}\bar{\rho}Y.$$
We note that $\bar{\rho}\bar{u}=m$ is real and constant on the profile. Then taking the standard complex $L^2$ innerproduct of the above equation with $Y$, we have 
$$\ip{\lam(\bar{\rho}Y)}{Y}+\ip{(mY)'}{Y}=\ip{(\bar{\rho}dY')'}{Y}-\ip{k\bar{\varphi}\bar{\rho}Y}{Y},$$
so that integrating by parts, and taking real parts yields
\begin{equation}
\re\lam\ip{\bar{\rho}Y}{Y}+0+\re\ip{k\bar{\varphi}\bar{\rho}Y}{Y}+\re\ip{Y'}{\bar{\rho}dY'}=0.
\label{important}\end{equation}
since $\ip{mY'}{Y}=\ip{(mY)'}{Y}=-\ip{mY}{Y'}=\overline{\ip{mY}{Y'}}$ so that $\ip{(mY)'}{Y}$ is purely imaginary. But for \eqref{important} to hold, we must have either $\re\lam<0$ or $\re\lam=0$ and also $Y'\equiv 0$ which implies that $Y$ is constant. But then we need also $\re\ip{k\bar{\varphi}\bar{\rho}Y}{Y}=0$, so that the constant value for $Y$ must be $0$. We conclude that the only nontrivial solutions must correspond to $\re\lam<0$. 

In the case that species diffusion is neglected things are even simpler; we write the reaction equation as
\begin{equation}
\lam z+(vz)'=-\bar{\varphi} z, \label{qzeval2}
\end{equation}
where $z=\rho Y$ as usual. Applying the Gap Lemma to \eqref{qzeval2}, we find that behavior at $+\infty$ is governed by the limiting constant-coefficient equations. These are easily seen to support no stable modes. Indeed, at $+\infty$, we have $\bar{\varphi}\equiv 0$, so there is no reaction at all. Then the equation
$$\lam z=-(vz)',$$
which can be rewritten in terms of the new variable $w=vz$ as
\begin{equation}
w'=-\left(\frac{\lam}{v}\right)w,
\end{equation} 
is obviously blowing up for $\re\lam>0$.
\foorp

\begin{cor}
For sufficiently small $q$, the integrated eigenvalue equations associated with \eqref{gas} and \eqref{react} have no $\re\lam\ge 0$ eigenvalues provided that the limiting shock is stable.
\end{cor}
\pf The respective Evans functions vary continuously, and the limiting Evans function is nonvanishing on $\re\lam\ge 0$.
\foorp
\begin{prop} If the limiting shock is stable, then then small-$q$ detonations are strongly spectrally stable. That is, they have no eigenvalues for $\re\lam\ge 0$ and $\lam\ne 0$. Moreover, there holds $D'(0)=\gamma_d\Delta\ne 0$, where $\gamma_d$ and $\Delta$ are as defined in \eqref{gammad} and \eqref{reactdelta}
\end{prop}
\pf
When $\lam\ne 0$, we may integrate the divergence form gas eigenvalue equation to deduce that $(U+\vec{q}z)$ has zero integral, and thus $\tilde U$, defined as $\int_{-\infty}^{x}(U+\vec{qq}z)$ lies in $L^2$ if $U,z$ do since by the gap lemma functions decay exponentially if at all as do their integrals \cite{ZH}. Thus existence of an eigenfunction for the linearized eigenvalue equation and its integrated version are equivalent. Existence of a transverse connection $\gamma_d\ne 0$ follows likewise by continuity from the result for the limiting equations, provided there exists a transverse connection for the limiting gas-dynamical shock. Finally $\Delta\ne 0$ follows by inspection from the corresponding property for the limiting gas-dynamical shock, since $\Delta$ for $q=0$ reduces to this case, which property is a necessary condition for stability \cite{ZH}. 
\foorp
The results of this section clearly extend to the case of multiple reactants. Indeed for the zero species diffusion case, the number and type of reactants plays no role. For the $D\ne 0$ case, the arguments above carry through if $D$ and $\Phi$ are simultaneously positive, or if there exists a constant coordinate change making them both positive. In particular the argument applies if the diffusion $D$ is scalar and $\Phi=\varphi K$ for $K$ constant and $\varphi$ scalar.

\appendix
\section{Real Viscosity}\label{appendix}
For completeness and convenience, we provide here a revised version of an appendix of \cite{Z4} using \cite{Zapp} and including a discussion of the extreme shock case in Section~\ref{evaluation} below. Systems modeling gas dynamics have the general form
\begin{equation}
U_t+F(U)_x=(B(U)U_x)_x,\label{realvisc}
\end{equation}
where 
$$U=\left(\begin{array}{c} u \\ v\end{array}\right),\quad F=\left(\begin{array}{c} f \\ g \end{array}\right),\quad B=\left(\begin{array}{cc} 0 & 0 \\ b_1 & b_2 \end{array}\right), $$ 
and
$$u,f\in \R^{n-r},\quad v,g,\in\R^r,\quad b_1\in\R^{r\times (n-r)},\quad b_2\in \R^{r\times r}.$$
We note that in the isentropic gas dynamics case $n=2$ and $r=1$ while in the Navier-Stokes case $n=3$ and $r=2$.
Our interest is in traveling wave solutions of the form
\begin{equation}
U=\bar{U}(x),\quad \lim_{x\rightarrow \pm\infty}\bar{U}(x)=U_{\pm}=(u_{\pm},v_{\pm}).\label{travelwave}
\end{equation}
Standard assumptions for equations in this generality are:
\begin{description}
\item[$\mathrm{H0}$] $F, B\in C^2. $
\item[$\mathrm{H1}$]
$\left\{\begin{array}{l} 
 \text{(i)} \re \sigma(b_2)>0 \text{ and} \\ 
 \text{(ii)}\binom{df}{b} \text{ full rank on }\{\bar{U}(\cdot)\}, \text{ moreover} \\
 \text{(iii)}(df_u -b_1(b_2)^{-1}df_v) \text{ is real}.
\end{array}\right.$
\item[$\mathrm{H2}$]
$
\sigma\left( dF(U_\pm)\xi\right)\ \mbox{real for}\;
\xi\in\R,\ 0\not\in \sigma(dF(U_\pm)). 
$
\item[$\mathrm{H3}$]
$
\re\sigma \left( \xi dF(U_\pm) -
\xi^2B(U_\pm)\right) \leq 0 
$
for $\xi \in\R$.
\item[$\mathrm{H4}$]
$
\hbox{Solutions of (\ref{realvisc})--(\ref{travelwave}) form a 
smooth manifold $\{
\bar{u}^{\delta}\},\ \delta\in \mathcal{U}\subset \R^\ell$.}
$
\end{description}
These hypotheses are analogous to those of the strictly
parabolic case considered in e.g \cite{Z4}, with (H1)(i) and (H1)(iii) ensuring local well-posedness.
Indeed, they are the standard set of conditions identified
by Kawashima \cite{Kaw}; for further discussion, see \cite{SeZ}.
The condition (H1)(ii), may
be motivated by consideration of the traveling wave ODE
\begin{equation}
f(u,v) \equiv f(u_-, v_-),\label{ODEmanifold} 
\end{equation}
\begin{equation}
b_1 u' + b_2 v' = g(u,v) - g(u_-,v_-). \label{ode}
\end{equation}
For, ($\mathrm{H1}$)(ii) is readily seen to be the condition
that (\ref{ode}) describes a nondegenerate ODE on the $r$-dimensional
manifold described by (\ref{ODEmanifold}); thus, this is a reasonable
nondegeneracy condition to impose in the study of viscous profiles.
Condition ($\mathrm{H1}$)(iii) also arises in the
analysis of the eigenvalue ODE; see the discussion of consistent splitting in Appendix A2, \cite{Z4}. In the 
symmetrizable case it holds automatically.

We remark, finally, that ($\mathrm{H1}$)(ii) (indeed, all
of hypothesis ($\mathrm{H1}$)) is satisfied for gas and plasma
dynamics precisely when particle and shock velocities are distinct,
which is always the case along a shock; for a study of viscous
profiles in these contexts, see \cite{Gi}, \cite{Gel}, \cite{FreS}.
 
Let $i_+$ denote the
dimension of the stable subspace of $df^1(u_+)$,
$i_-$ denote the dimension of the unstable subspace of $df^1(u_-)$,
and $i:=i_+ + i_-$.  
Let $d_+$ denote the dimension within the submanifold $f\equiv \text{\rm 
constant}$ of the stable manifold at $(u_+,v_+)$ of traveling
wave ODE \eqref{ode}, and $d_-$ the dimension
of the unstable manifold at $(u_-,v_-)$, and $d:=d_-+d_+$.
Then, we have the following result 
analogous to that of Majda and Pego [MP] in the strictly
parabolic case.

\begin{lemma}
Under assumptions (H0)--(H3), 
$(u_\pm,v_\pm)$
are hyperbolic rest points of the reduced traveling
wave ODE \eqref{ode}.
In particular, traveling wave solutions satisfy
\begin{equation}
|(d/dx)^k \big((\bar u(x),\bar v(x))-(u_\pm,v_\pm)\big)|
\le Ce^{-\theta |x|}, \quad  k=0,\dots,4,
\label{realexpdecay}
\end{equation}
as $x\to \pm \infty$.
Moreover, the type of the connection agrees with
the (hyperbolic) type of the shock, in the sense that 
\begin{equation}
d-r=i-n.
\label{ireal}
\end{equation}
\end{lemma}
\pf 
Integrating \eqref{ode} from $-\infty$ to $x$ and
rearranging, we may write \eqref{ODEmanifold}--\eqref{ode}
in the alternative form
\begin{equation}
\begin{pmatrix} u \\v\end{pmatrix}'=
\begin{pmatrix} f_u & f_v \\
b_1 & b_2 \end{pmatrix}^{-1}
\begin{pmatrix} 0 \\ g- g_- \end{pmatrix}.
\label{uvode}
\end{equation}
Linearizing \eqref{uvode} about $U_\pm$, we obtain
\begin{equation}
\begin{pmatrix} u \\v\end{pmatrix}'=
\begin{pmatrix} f_u & f_v \\
b_1 & b_2 \end{pmatrix}^{-1}
\begin{pmatrix} 0 & 0 \\ g_u & g_v \end{pmatrix}
_{|(U_\pm)}
\begin{pmatrix} u \\v\end{pmatrix},
\label{luvode}
\end{equation}
or, setting 
$$
\begin{pmatrix} z_1 \\ z_2 \end{pmatrix} :=
\begin{pmatrix} f_u & f_v \\
b_1 & b_2 \end{pmatrix}_{|(U_\pm)}
\begin{pmatrix} u \\ v\end{pmatrix},
$$
the pair of equations 
$$
z_1'=0,
$$
and
\begin{equation}
z_2'=
\begin{pmatrix} g_u & g_v \end{pmatrix}
\begin{pmatrix} f_u & f_v \\
b_1 & b_2 \end{pmatrix}^{-1}
\begin{pmatrix} 0 \\ I_r \end{pmatrix}
_{|(U_\pm)}
z_2,
\label{linman}
\end{equation}
the latter of which evidently describes the linearized ODE
on manifold \eqref{ODEmanifold}.
Observing that
$$
\det
\begin{pmatrix} g_u & g_v \end{pmatrix}
\begin{pmatrix} f_u & f_v \\
b_1 & b_2 \end{pmatrix}^{-1}
\begin{pmatrix} 0 \\ I_r \end{pmatrix}_{|(U_\pm)}
=
\det
\begin{pmatrix} f_u & f_v\\ g_u & g_v \end{pmatrix}
\begin{pmatrix} f_u & f_v \\
b_1 & b_2 \end{pmatrix}^{-1}_{|(U_\pm)}
\ne 0
$$
by (H2) and (H1)(i), we find that the coefficient matrix
of \eqref{linman} has no zero eigenvalues.
On the other hand, it can have no
nonzero purely imaginarly eigenvalues $i\xi$,
since otherwise
$$
\begin{pmatrix} f_u & f_v\\ g_u & g_v \end{pmatrix}
\begin{pmatrix} f_u & f_v \\ b_1 & b_2\end{pmatrix}^{-1}
 \begin{pmatrix} 0\\v\end{pmatrix}=
i\xi
\begin{pmatrix} 0\\v\end{pmatrix},
$$
and thus
$$
\Big[ -i\xi \begin{pmatrix} f_u & f_v\\ g_u & g_v \end{pmatrix}
-\xi^2 \begin{pmatrix} 0 & 0 \\ b_1 & b_2\end{pmatrix} \Big]
\Big[\begin{pmatrix} f_u & f_v \\ b_1 & b_2\end{pmatrix}^{-1}\begin{pmatrix} 0\\v\end{pmatrix}
\Big]
=\begin{pmatrix} 0\\0\end{pmatrix}
$$
for $\xi\ne 0 \in \R$,
in violation of (H3).
Thus, we find that $U_\pm$ are hyperbolic
rest points, from which \eqref{realexpdecay} follows.
Relation \eqref{ireal} now follows from Lemma \ref{dimensionlemma}, below.
\foorp
\noindent
The linearized eigenvalue equations about $\bar{U}(\cdot)$ are:
\begin{equation}
(A_{11}u + A_{12}v)' = -\lambda u,
\label{appeval1}
\end{equation}
and
\begin{equation}
(b_1 u' + b_2 v')' = (A_{21}u + A_{22}v)', \label{appeval2}
\end{equation}
where
$$
\begin{pmatrix}
0 & 0 \\
b_1 & b_2
\end{pmatrix} := B(\bar{U}),
$$
$$
\begin{pmatrix}
A_{11} & A_{12} \\
A_{21} & A_{22}
\end{pmatrix} U:= dF(\bar{U})U - dB(\bar{U}) (U,\bar{U}'_{x}),
$$
and `$'$' denotes $\partial/\partial x$; in particular, note
that 
$$
(A_{11}, A_{12}) = df(\bar{U}).
$$
Utilizing the change of variables
\begin{equation}
\left(\begin{array}{c}z_1 \\ z_2 \end{array}\right) =\left(\begin{array}{cc} A_{11} & A_{12} \\ b_1 & b_2 \end{array}\right)\left(\begin{array}{c} u \\ v \end{array}\right),
\end{equation}
we can write the eigenvalue equation as a first order system
\begin{equation}
Z'=\mathbb{A}(x,\lam)Z,\quad Z=(z_1,z_2,z_2')^{tr}.\label{realviscfirstsystem}
\end{equation}
We note here that $z_2$ has $r$ components. 

The consistent splitting 
hypothesis can be verified by a limiting analysis as
$\lambda\to+\infty$, carried out without loss of generality in original
coordinates $W$, for which the asymptotic characteristic equations
become: 
\begin{equation}
\det
\begin{pmatrix}
\mu A_{11} + \lambda & \mu A_{12} \\
\mu A_{12} - \mu^2b_1 & \mu A_{22} - \mu^2b_2 + \lambda
\end{pmatrix}_\pm  
\begin{pmatrix} u \\ v \end{pmatrix} = \begin{pmatrix} 0 \\ 0 \end{pmatrix}.
\label{realviscchar}
\end{equation}
This yields $n-r$ roots $\mu\sim \tilde{\mu}\lambda$, 
$\tilde{\mu}\sim 1$, satisfying
\begin{equation}
\begin{pmatrix}
\tilde{\mu}A_{11} + I & A_{12}\\
b_1 & b_2 \end{pmatrix}_\pm 
\begin{pmatrix} u \\ v\end{pmatrix} = \begin{pmatrix} 0  \\ 0\end{pmatrix}, 
\label{roots1}
\end{equation}
or
\begin{equation}
-\tilde{\mu}^{-1}\in \sigma (A_{11} - A_{12} b_2^{-1}
b_1)_\pm , \label{growth1}
\end{equation}
and $2r$ roots $\mu \sim \tilde{\mu}\lambda^{1/2},\ \tilde{\mu}\sim
1$, satisfying
\begin{equation}
\begin{pmatrix}
I & 0 \\
-\tilde{\mu}^2 b_1 & - \tilde{\mu}^2 b_2 + I
\end{pmatrix}_\pm \begin{pmatrix} u  \\ v\end{pmatrix} = \begin{pmatrix} 0  \\0\end{pmatrix},
\label{roots2}
\end{equation}
or
\begin{equation}
\tilde{\mu}^{-2}\in\sigma(b_2). \label{growth2}
\end{equation}
By assumption ($\tilde{\mathrm{H1}}$)(iii), (\ref{growth1}) yields a fixed
number $k$/$(n-r-k)$ of stable/unstable roots, independent of $x$, and
thus of $\pm$.  Likewise, ($\tilde{\mathrm{H1}})$(i) implies that
(\ref{growth2}) yields $r$ stable, $r$ unstable roots.
Combining, we find the desired consistent splitting, with $(k+r)$/ $(n-k)$
stable/unstable roots at both $\pm \infty$.  We can thus define an
Evans function as usual as
\begin{equation}
D(\lam)= \det (Z^+_1, \dots Z^+_{k+r}, Z^-_{k+r+1}, \dots
Z^-_{n+r})_{|x=0} \label{realviscevans}
\end{equation}
where $\{Z^+_1, \dots Z^+_{k+r}\}$, $\{Z^-_{k+r+1}, \dots Z^-_{n+r}\}$
span the stable manifold at $+\infty$, unstable manifold at $-\infty$
of (\ref{realviscfirstsystem}). Notice that the Evans function in $Z$ coordinates is just a
constant multiple of the corresponding Evans function defined in 
$W=(u,v,v')^t$ coordinates. 

\subsection{Stability Index}
The stability index is defined to be 
\begin{equation}
\tilde{\Gamma}:=\sgn (\partial_{\lam})^{\ell}D(0)D(+\infty). \label{index}
\end{equation}
The low frequency calculations of $D'(0)$ are detailed for the Navier-Stokes model in Section 3. Note that $\ell =1$ in this case. Here we evaluate the sign of \evans\ as $\lam\to +\infty$ along the real axis.
\begin{lemma} \label{Clemma1} Let $\tilde D(\cdot)$ denote the alternative Evans function
\begin{equation}
\tilde D(\lambda):= \det (W^+_1, \dots W^+_{k+r}, W^-_{k+r+1}, \dots
W^-_{n+r})_{|x=0} 
\label{Dtilde} 
\end{equation}
computed in the original coordinates $W$.
Then, for real $\lam$ sufficiently large, there holds 
\begin{equation}
\sgn \tilde D(\lambda ) = 
\sgn\det (\mathbb{S}^+, \mathbb{U}^+)  \det(\pi\mathbb{W}^+, \varepsilon \mathbb{S}^+)
\det (\varepsilon \mathbb{U}^-, \pi \mathbb{W}^-) 
\ne 0,
\label{largelam1}
\end{equation}
where $\pi$ denotes projection of $W=(u,v,v')$ onto $(u,v)$ components, and
$\mathbb{S}(x),\ \mathbb{U}(x)$ are real bases of the stable/unstable subspaces
of $(A_{11}-A_{12}b_2^{-1}b_1)$ (note: $(n-r)$
dimensional), with $\varepsilon u:=(u, -b_2^{-1}
b_1u)$ denoting extension. \label{deltanotvanish}
\end{lemma}
%
\begin{proof}
Recalling ($\tilde{\mathrm{H1}}$)(ii), we have that
the coordinate change$(u,v)\to(z_1,z_2)$ is invertible, 
and so we may work equivalently in $(u,v,z_3)$ coordinates,
where $z_3:= b_1 u' + b_2 v'$.  
Then, we find from
(\ref{roots1}), (\ref{roots2}) that the 
stable/unstable manifolds of the frozen
eigenvalue equation at any (fixed) $x$ are spanned by vectors of form 
$$
\begin{pmatrix}
u \\
-b_2^{-1}b_1 u \\
*
\end{pmatrix},
$$
with $u$ an unstable/stable eigenvector of $(A_{11} - A_{12}
b_2^{-1}b_1)$,
$-\tilde{\mu}^{-1}$ the corresponding eigenvalue; and vectors
$$
\begin{pmatrix}
0 \\
v \\
\mp \tilde{\mu}\lambda^{1/2}v
\end{pmatrix}
$$
with $v$ an eigenvector of $b_2$, $-\tilde{\mu}^{-2}$ the
corresponding eigenvalue. 
Rescaling and applying the tracking lemma, we thus obtain
\begin{multline}
D(\lam)\sim \det\left( \underbrace{\begin{array}{cc} u & \cdots \\ -b_2^{-1}b_1u & \cdots \\ * & \cdots \end{array}}_{n-r} \underbrace{\begin{array}{ccc} 0 & \cdots & 0  \\ v & \cdots & v \\ -\tilde{\mu}\lam^{1/2}v & \cdots & \tilde{\mu}\lam^{1/2}v \end{array}}_{2r}\right) \\
\times\det \left(\begin{array}{cc|cc} \alpha_{+}^1 & 0 & \alpha_+^2 & 0 \\ 0 & \alpha_-^1 & 0 & \alpha_-^2 \\ \hline \beta_+^2 & 0 & \beta_1^+ & 0 \\ 0 & \beta_-^2 & 0 & \beta_-^1 \end{array}\right) \label{largelamdelta}
\end{multline}
where both
$$\det\left(\begin{array}{cc} \alpha_+^1 & \alpha_+^2 \\ \beta_+^2 & \beta_+^1 \end{array}\right) \quad\mbox{and}\quad
\det\left(\begin{array}{cc} \alpha_-^1 & \alpha_-^2 \\ \beta_-^2 & \beta_-^1 \end{array}\right) $$
are real, nonzero quantities. The right hand side of \eqref{largelamdelta} can be rewritten as 
\begin{equation}
\det(\mathbb{S}^+,\mathbb{U}^+)\det(\mathbb{V})^2\det\begin{pmatrix}  \alpha_+^1 & \alpha_+^2 \\ \beta_+^2 & \beta_+^1 \end{pmatrix}\det\begin{pmatrix} \alpha_-^1 & \alpha_-^2 \\ \beta_-^2 & \beta_-^1\end{pmatrix},
\end{equation}
where by $\det(\mathbb{V})$ we refer to the $r\times r$ determinant coming from the $v$ component of \eqref{largelamdelta}. On the other hand, the term 
$$\det(\pi\mathbb{W}^+,\varepsilon\mathbb{S}^+)$$
in \eqref{largelam1} can be simplified to 
\begin{equation}
\det(\pi\widetilde{\mathbb{W}}^+,\pi\mathbb{U}^+,\varepsilon\mathbb{S}^+)\det\left(\begin{array}{cc|c} \alpha_+^1 & \alpha_+^2 & 0 \\ \beta_+^2 & \beta_+^1 & 0 \\ \hline 0 & 0 & I \end{array}\right), \label{simplifiedterm}
\end{equation} 
where $(\widetilde{\mathbb{W}}^+,\mathbb{U}^+)=\mathbb{W}^+$, but $\det(\pi\widetilde{\mathbb{W}}^+,\pi\mathbb{U}^+,\varepsilon\mathbb{S}^+)$ has the form
$$
\det\left(\begin{array}{c|cc}0 & \mathbb{U}^+ & \mathbb{S}^+ \\ \mathbb{V} & * & * \end{array}\right).
$$
Finally we see that \eqref{simplifiedterm} is simply 
\begin{equation}
\det(\mathbb{V})\det(\mathbb{S}^+,\mathbb{U}^+)\det\begin{pmatrix}  \alpha_+^1 & \alpha_+^2 \\ \beta_+^2 & \beta_+^1 \end{pmatrix}. \label{firsthalf}
\end{equation}
Similarly, the term 
$$\det(\varepsilon\mathbb{U}^-,\pi\mathbb{W}^-)$$ simplifies to 
\begin{equation}
\det(\mathbb{V})\det(\mathbb{S}^-,\mathbb{U}^-)\det\begin{pmatrix} \alpha_-^1 & \alpha_-^2 \\ \beta_-^2 & \beta_-^1\end{pmatrix} \label{secondhalf}
\end{equation}
combining \eqref{firsthalf} and \eqref{secondhalf} we find since $\sgn \det(\mathbb{S}^+,\mathbb{U}^+)=\sgn\det(\mathbb{S}^-,\mathbb{U}^-)$, that the expressions \eqref{largelam1} and \eqref{largelamdelta} agree modulo the real, positive factor $\det(\mathbb{S}^+,\mathbb{U}^+)^2$.  
\end{proof}
%
%
\noindent
We further make the assumptions of \emph{semidissipativity}

(+) \quad  There exist symmetrizers $A^0_\pm$ such that
$A^0_\pm A_\pm$ are symmetric and $\re A^0_\pm B_\pm \leq 0$, 

and \emph{block structure}

(++) \quad $\displaystyle{ (A^0_\pm)^{1/2} B_\pm (A^0_\pm)^{-1/2}
= \begin{pmatrix}
0 & 0 \\
0 & \tilde{b}_2
\end{pmatrix}_\pm }.
$

Both of these assumptions hold for the compressible Navier-Stokes equations. When they hold, more can be said about the sign of $D(\lam)$ for large real \lam\ .
\begin{lemma}
Assuming (+)--(++), there holds, for
sufficiently large, real $\lambda$:
\begin{equation}
\sgn \tilde D(\lambda) =
\sgn\det (\mathbb{S}^+, \mathbb{U}^+) \det (\pi\mathbb{W}^+,
\varepsilon\mathbb{S}^+) 
\det( \varepsilon\mathbb{U}^-, \pi \mathbb{W}^-)_{|_{\lam=0}} \ne 0,
\label{tolamzero}\end{equation}

where $\tilde D(\cdot)$ as in (\ref{Dtilde}) denotes the Evans function
computed in original coordinates $W=(u,v,v')^t$. \label{Dinfty}
\end{lemma}

\pf Without loss of generality, we may take $A$
symmetric, $B = \left( \begin{matrix}
0 & 0 \\
0 & b_2
\end{matrix}\right)$, and $\re\,~b_2 >0$, by the transformation 
$$
A \to (A^0)^{\frac{1}{2}} A (A^0)^{-\frac{1}{2}}, 
\quad B \to (A^0)^{\frac{1}{2}} B(A^0)^{-\frac{1}{2}}.  
$$
It is sufficient to show that quantity (\ref{tolamzero}) does not
vanish in the class (+)--(++).  For, since $D(\lambda)$ does not
vanish either, for real $\lambda$ sufficiently large,
we can then establish the result by homotopy of the symmetric matrix $ A_\pm$ to an invertible diagonal matrix
(straightforward, using the unitary decomposition $A=UDU^*$, 
$U^*U=I$, and the fact that unitary matrices are homotopic
either to $I$ or $-I$)
and of 
$ B_\pm$ to
$\begin{pmatrix} 
0 & 0 \\
0 & I_r
\end{pmatrix}$
(e.g., by linear interpolation of the positive definite $b^{11}$ to $I_r$),
in which case it can be seen by explicit computation
that \eqref{tolamzero} is independent of $\lambda\in [0,+\infty]$.
%
We note that the endpoint of this homotopy
is on the boundary of but not in
the Kawashima class, since eigenvectors of $A$ are
in kernel of $B$; indeed, our definition of semidissipativity is
is {\it not} the ``strict'' dissipativity condition of Kawashima,
but a nonstrict version.
However, it suffices for the present, purely linear-algebraic purpose.

We begin by examining $\det (\pi\mathbb{W}^+,\varepsilon\mathbb{S}^+)$. When $\lambda=0$, a bifurcation analysis as in \cite{Z4} of the limiting constant-coefficient equations at $\pm \infty$
shows that the projections $\pi$ of 
slow modes of $\mathbb{W}^+$ may be chosen as the unstable
eigenvectors $r^+_j$ of $A$, corresponding to outgoing characteristic modes,
and the projections of fast modes as the stable (i.e. $\re \mu<0$)
solutions of
$$
(A - \mu B)_\pm \binom{u}{v} = \binom{0}{0}, 
$$
or without loss of generality
$$
\begin{pmatrix}
A_{11} & A_{12} \\
A_{21} & A_{22}-\mu b_2
\end{pmatrix}_\pm \begin{pmatrix} u \\ v \end{pmatrix} = \begin{pmatrix} 0 \\ 0\end{pmatrix},
$$
and thus of form
$$
\begin{pmatrix}
- (A_{11})^{-1}A_{12}v \\
v
\end{pmatrix}, 
$$
where
$$
\left( b^{-1}_2(A_{22} - A_{21}(A_{11})^{-1}A_{12})-\mu
I\right)v=0. 
$$
Likewise, using $b_1=0$, we find from 
the definitions of $\mathbb{S}$, $\varepsilon$ in the statement of
Lemma \ref{deltanotvanish} that
stable solutions $\mathbb{S}^+$ are in the stable subspace of
$A_{11}$, with $\varepsilon\, \mathbb{S}^+=\binom{\mathbb{S}^+}{0}$, hence 
vectors $\varepsilon\mathbb{S}^+$ lie
in the intersection of the stable subspace of $A$ and 
the kernel of $B$.  
Our claim is that these three subspaces are independent,
spanning $\C^n$.  Rewording this assumption, we are claiming that
the stable subspace of $(A)^{-1}B$, the center subspace 
$\ker B$ intersected with the stable subspace of $A$, and the unstable
subspace of $A$ are mutually independent.
(Note: that dimensions are correct follows by consistent splitting). 
But, this follows by Lemma \ref{serrelemma} below. Similar considerations apply to $\det( \varepsilon\mathbb{U}^-, \pi \mathbb{W}^-)$.
\foorp

A key to the calculations above is the following modified lemma of Serre from \cite{Z4}. We begin by fixing notation. For a matrix $M$, we denote by $\mathcal{S}(M)$ and $\mathcal{U}(M)$ the stable and unstable subspaces of $M$. Also we we have
$$\mathcal{N}_{V}(M)=\{v\in V\; | \; \langle v,M\rangle <0\}$$
and
$$\mathcal{P}(M)=\{v\; | \; \langle v,Mv\rangle >0\}.$$

\begin{lemma}[Modified Serre's Lemma \cite{Z4}]Let $A$ be a symmetric, invertible matrix and let $B$ be a positive semidefinite matrix, $\re (B)\geq 0$. Then
\begin{enumerate}
\item the subspaces $\mathcal{S}(A^{-1}B)\oplus \mathcal{N}_{\mbox{ker}B}(A)$ and $\mathcal{U}(A)$ are transverse.
\item the subspaces $\mathcal{S}(BA^{-1})\oplus \mathcal{N}_{\mbox{ker}B}(A)$ and $\mathcal{U}(A)$ are transverse.
\end{enumerate} \label{serrelemma}
\end{lemma}
\pf It suffices to prove the first claim, as the second follows immediately from the first due to the similarity transform $B\rightarrow ABA^{-1}$. Since $A$ is symmetric, the unstable subspace $\mathcal{U}(A)$ is equal to the subspace $\mathcal{P}(A)$, which for notational convenience we denote by $\mathcal{P}$. Suppose, in order to obtain a contradiction, that $x_0\neq 0$ lies both in the subspace  $\mathcal{S}(A^{-1}B)\oplus \mathcal{N}_{\mbox{ker}B}(A)$ and in $\mathcal{P}$. Thus we suppose that 
$$x_0=x_1+x_2$$
where $x_1\in\mathcal{S}(A^{-1}B)$ and $x_2\in \mathcal{N}_{\mbox{ker}B}(A)$. Define $x(t)$ by the ordinary differential equation
\begin{equation}x'=A^{-1}Bx,\quad x(0)=x_0. \end{equation}
It follows then that $x(t)\rightarrow x_2$ as $t\rightarrow +\infty$ and thus 
\begin{equation}\lim_{t\rightarrow +\infty}\langle x(t),Ax(t)\rangle \leq 0. \end{equation}
However,
\begin{eqnarray*}
\langle x, Ax\rangle' & = & 2\langle A^{-1}Bx,Ax\rangle \\
		         & = & 2\langle Bx,x\rangle  \\
                                  & \geq & 0.
\end{eqnarray*}
This implies that $\langle x_0,Ax_0\rangle \leq 0$ which is a contradiction to our assumption that $x_0$ belongs to $\mathcal{P}$. \foorp
%
%
%
\subsection{Evaluation of $\tilde{\Gamma}$}\label{evaluation}

Using the $z$-coordinates we may regard the traveling wave ODE as an $r$ dimensional first order dynamical system. We denote by $d_{\pm}$ the dimensions of the stable manifold at $z_{2+}$ and the unstable manifold at $z_{2-}$. We also define $d$ to be the sum $d=d_{+}+d_{-}$. It follows then that $1\leq d_{\pm}\leq r$. Also we denote by $i_{\pm}$ the number of characteristics entering the shock from the left $(-)$ and the right $(+)$. We put $i=i_++i_-$. Corresponding to \cite{MP}, we have:
\begin{lemma} \label{dimensionlemma}
The following relations hold
\begin{enumerate}
\item $n-i_+=r-d_++\dim \mathcal{U}(\tilde{A}_+)$
\item $n-i_-= r-d_-+\dim \mathcal{S}(\tilde{A}_-)$
\end{enumerate}
where $\tilde{A}=(A_{11}-A_{12}b_2^{-1}b_1)$. Moreover $n-i=r-d$. \label{majdapego}
\end{lemma}
\pf Equating the dimensions of $\mathbb{Z}^+$ at $\lam =0$ and as $\lam\rightarrow\infty$ we find
$$\dim \mathcal{U}(A_+)+d_+=\dim \mathcal{U}(\tilde{A}_+)+r,$$
or
$$(n-i_+)+d_+=\dim \mathcal{U}(\tilde{A}_+)+r.$$
Similarly as $x\rightarrow -\infty$, we find
$$(n-i_-)+d_-=\dim \mathcal{S}(\tilde{A}_-)+r.$$
That $n-i=r-d$ follows by adding the two equations and noting that $\dim \mathcal{U}(\tilde{A}_+)$ and $\dim\mathcal{S}(\tilde{A}_-)$ are constant and sum to $n-r$. \foorp

\begin{cor}\label{extreme} For (right) extreme shocks, $i_+=n$, there holds also $d_+=r$. Thus the connection is also extreme and $\dim \mathcal{U}(\tilde{A})\equiv 0$.
\end{cor}
\pf This follows at once from Lemma~\ref{majdapego} due to the fact that $d_+\leq r$ and $\dim\mathcal{U}(\tilde{A})\geq 0$. \foorp
\noindent
The import of Lemma~\ref{majdapego} is that the ``parabolic'' and ``hyperbolic'' types of connections agree.
From Corollary~\ref{extreme}, we may deduce that $\gamma$ for an extreme right (i.e. $n$-shock) Lax profile consists of a Wronskian involving only modes from th $+\infty$ side, and is therefore explicitly evaluable.
For, working now in $(z_1,z_2,z_2')$ coordinates, we obtain $\gamma$ as a determinant of $z_2$ components only.

Moreover, the expression \eqref{tolamzero} simplifies greatly. In the $(z_1,z_2)$ coordinates we have that $\mathbb{U}=\emptyset$, $\mathbb{S}$
is full dimension $n-r$, and $\epsilon \mathbb{S}$ consists
of vectors of the simple form $(z,0)$.
This means that $\det (\mathbb{S},\mathbb{U})$ simplifies to
just $\det \mathbb{S}$, while $\det(\pi\mathbb{Z}^+, \epsilon \mathbb{S}^+)$
simplifies to  the product of $\det \mathbb{S}^+$ and
$\tilde{\gamma}$.  {\it Therefore, this term, similarly as
in the strictly parabolic case, cancels with term $\tilde{\gamma}$
in the computation of the stability index.}

Finally, $\det (\epsilon \mathbb{U}^-, \pi \mathbb{Z}^-)$
simplifies to $\det(r_1^-,\dots, r_{n-1}^-, \bar{u}')$,
times the determinant of the coordinate transformation
from $W$ to $Z$ coordinates, the latter determinant cancelling
with a like factor appearing above. 
We are left in the end with the following very simple formula.
\begin{prop}
In the case of an extreme right ($n$-shock) Lax profile,
\begin{equation}
\tilde{\Gamma}=\sgn\det(r_1^-,\ldots,r_{n-1}^-,[u])\det(r_1^-,\ldots,r_{n-1}^-,\bar{u}'/|\bar{u}'|(-\infty))
\end{equation}
\end{prop}
We emphasize that this is \emph{identical with the stability index in the strictly parabolic case}. The only very weak information required from the connection problem is the orientation of $\bar{u}'$ as $x\to -\infty$, i.e. the direction in which the profile leaves along the one-dimensional unstable manifold. We remark that in the case of isentropic gas dynamics, the traveling wave ODE is scalar, and thus the orientation of $\bar{u}'$ is determined by the direction of the connection. See Section~2.4 or \cite{Lyng} for further details.

%
%
%
%
%
\nocite{BarEgo,MenPlohr}

\bibliographystyle{abbrv}

\bibliography{refs}

\def\cprime{$'$}
\begin{thebibliography}{10}

\bibitem{AGJ}
J.~Alexander, R.~Gardner, and C.~Jones.
\newblock A topological invariant arising in the stability analysis of
  travelling waves.
\newblock {\em J. Reine Agnew Math.}, 410:167--212, 1990.

\bibitem{BarEgo}
A.~Barmin and S.~Egorushkin.
\newblock Stability of shock waves.
\newblock {\em Adv. Mech.}, 15(1-2):3--37, 1992.

\bibitem{BSZ}
S.~Benzoni-Gavage, D.~Serre, and K.~Zumbrun.
\newblock Alternate {E}vans functions and viscous shock waves.
\newblock {\em SIAM J. Math. Anal.}, 32(5):929--962 (electronic), 2001.

\bibitem{BMR}
A.~Burlioux, A.~J. Majda, and V.~Roytburd.
\newblock Theoretical and numerical structure for unstable one-dimensional
  detonations.
\newblock {\em SIAM Journal on Applied Mathematics}, 51(2):303--343, 1991.

\bibitem{CampWood1}
C.~Campbell and D.~Woodhead.
\newblock The ignition of gases by an explosion wave. {I}.carbon monoxide and
  carbon mixtures.
\newblock {\em J. Chem. Soc.}, 129:3010--3021, 1926.

\bibitem{CampWood2}
C.~Campbell and D.~Woodhead.
\newblock Striated photographic records of explosion waves.
\newblock {\em J. Chem. Soc.}, 130:1572--1578, 1927.

\bibitem{CHT}
G.-Q. Chen, D.~Hoff, and K.~Trivisa.
\newblock On the {N}avier-{S}tokes equations for exothermically reacting
  compressible fluid.
\newblock {\em Acta Mathematicae Applicatae Sinica, English Series}, 18(1),
  2002.

\bibitem{Co}
W.~A. Coppel.
\newblock {\em Stability and Asymptotic Behavior of Differential Equations}.
\newblock D. C. Heath and Co., Boston, Mass., 1965.

\bibitem{Ev1}
J.~Evans.
\newblock Nerve axon equations: {I}. {L}inear approximations.
\newblock {\em Indiana University Math Journal}, 21:877--855, 1972.

\bibitem{Ev2}
J.~Evans.
\newblock Nerve axon equations: {II}. {S}tability at rest.
\newblock {\em Indiana University Math Journal}, 22:75--90, 1972.

\bibitem{Ev3}
J.~Evans.
\newblock Nerve axon equations: {III}. {S}tability of the nerve impulse.
\newblock {\em Indiana University Math Journal}, 22:577--593, 1972.

\bibitem{Ev4}
J.~Evans.
\newblock Nerve axon equations: {IV}. {T}he stable and unstable impulse.
\newblock {\em Indiana University Math Journal}, 24:1169--1190, 1975.

\bibitem{Fen}
N.~Fenichel.
\newblock Geometric singular perturbation theory for ordinary differential
  equations.
\newblock {\em J. Differential Equations}, 31(1):53--98, 1979.

\bibitem{FD}
W.~Fickett and W.~Davis.
\newblock {\em Detonation: Theory and Experiment}.
\newblock U of California Press, Berkeley, 1979.

\bibitem{FickWood}
W.~Fickett and W.~Wood.
\newblock Flow calculations for pulsating one-dimensional detonations.
\newblock {\em Physics of Fluids}, 9:903--916, 1966.

\bibitem{FreS}
H.~Freist{\"u}hler and P.~Szmolyan.
\newblock Existence and bifurcation of viscous profiles for all intermediate
  magnetohydrodynamic shock waves.
\newblock {\em SIAM J. Math. Anal.}, 26(1):112--128, 1995.

\bibitem{GZ}
R.~Gardner and K.~Zumbrun.
\newblock The gap lemma and geometric criteria for the instability of viscous
  shocks.
\newblock {\em Communications in Pure and Applied Mathematics}, 51(7):797--855,
  1998.

\bibitem{GS1}
I.~Gasser and P.~Szmolyan.
\newblock A geometric singular perturbation analysis of detonation and
  deflagration waves.
\newblock {\em SIAM Journal of Math Analysis}, 24(4):968--986, 1993.

\bibitem{GS2}
I.~Gasser and P.~Szmolyan.
\newblock Detonation and deflagration waves with multistep reaction schemes.
\newblock {\em SIAM Journal of Applied Math}, 55(1):175--191, 1995.

\bibitem{Gel}
I.~M. Gel{\cprime}fand.
\newblock Some problems in the theory of quasilinear equations.
\newblock {\em Amer. Math. Soc. Transl. (2)}, 29:295--381, 1963.

\bibitem{Gi}
D.~Gilbarg.
\newblock The existence and limit behavior of the one-dimensional shock layer.
\newblock {\em Amer. J. Math.}, 73:256--274, 1951.

\bibitem{GordMooHarp}
W.~Gordon, A.~Mooradian, and S.~Harper.
\newblock Limit and spine effects in hydrogen-oxygen detonations.
\newblock In {\em Seventh Symposium (International) on Combustion}, pages
  752--759. Academic Press, 1959.

\bibitem{He}
D.~Henry.
\newblock {\em Geometric theory of semilinear parabolic equations}.
\newblock Springer-Verlag, Berlin, 1981.

\bibitem{HR}
M.~Hesaaraki and A.~Razani.
\newblock Detonative travelling waves for combustions.
\newblock {\em Applicable Analysis}, 77(3-4):405--418, 2001.

\bibitem{Hu}
J.~Humpherys.
\newblock On spectral stability of strong shocks for isentropic gas dynamics.
\newblock in preparation.

\bibitem{JL}
H.~K. Jenssen and G.~Lyng.
\newblock Evaluation of the {L}opatinski determinant for multi-dimensional
  {E}uler equations, 2002.
\newblock appendix to \cite{Zhandbook}.

\bibitem{KS}
T.~Kapitula and B.~Sandstede.
\newblock Stability of bright solitary-wave solutions to perturbed nonlinear
  {S}chr\"odinger equations.
\newblock {\em Phys. D}, 124(1-3):58--103, 1998.

\bibitem{KasStew}
A.~Kasimov and D.~S. Stewart.
\newblock Spinning instability of gaseous detonations.
\newblock {\em Journal of Fluid Mechanics}, 466:179--203, 2002.

\bibitem{Ka}
T.~Kato.
\newblock {\em Perturbation Theory for Linear Operators}.
\newblock Springer-Verlag, Berlin, 1995.
\newblock Reprint of the 1980 edition.

\bibitem{Kaw}
S.~Kawashima.
\newblock {\em Systems of a hyperbolic-parabolic type with applications to the
  equations of magnetohydrodynamics}.
\newblock PhD thesis, Kyoto University, 1983.

\bibitem{LeeStew}
H.~Lee and D.~Stewart.
\newblock Calculation of linear detonation instability: One-dimensional
  instability of plane detonation.
\newblock {\em J. Fluid Mech.}, 216:103--132, 1990.

\bibitem{LLT}
D.~Li, T.-P. Liu, and D.~Tan.
\newblock Stability of strong detonation waves to combustion model.
\newblock {\em Journal of Mathematical Analysis and Applications},
  201:516--531, 1996.

\bibitem{Li2}
T.~Li.
\newblock On the {R}iemann problem for a combustion model.
\newblock {\em SIAM Journal of Mathematical Analysis}, 24(1):59--75, 1993.

\bibitem{Li4}
T.~Li.
\newblock On the initiation problem for a combustion model.
\newblock {\em Journal of Differential Equations}, 112:351--373, 1994.

\bibitem{Li1}
T.~Li.
\newblock Rigorous asymptotic stability of a {C}hapman-{J}ouget detonation wave
  in the limit of small resolved heat release.
\newblock {\em Combustion Theory and Modeling}, 1(3):259--270, 1997.

\bibitem{Li3}
T.~Li.
\newblock Stability of strong detonation waves and rates of convergence.
\newblock {\em Electronic Journal of Differential Equations}, 1998(9):1--77,
  1998.

\bibitem{Li5}
T.~Li.
\newblock Stability and instability of detonation waves.
\newblock In {\em Hyperbolic Problems: Theory, Applications \& Numerics;
  Seventh International Conference in Z\"{u}rich}, 1999.

\bibitem{LYi}
T.-P. Liu and L.~Ying.
\newblock Nonlinear stability of strong detonations for a viscous combustion
  model.
\newblock {\em SIAM Journal of Math Analysis}, 26(3):519--528, 1995.

\bibitem{LY}
T.-P. Liu and S.~Yu.
\newblock Nonlinear stability of weak detonation waves for a combustion model.
\newblock {\em Communications in Mathematical Physics}, 204:551--586, 1999.

\bibitem{Lyng}
G.~Lyng.
\newblock {\em One Dimensional Stability of Detonation Waves}.
\newblock PhD thesis, Indiana University, 2002.

\bibitem{LyZ2}
G.~Lyng and K.~Zumbrun.
\newblock A stability index for detonation waves in {M}ajda's model for
  reacting flow, 2003.
\newblock preprint.

\bibitem{M}
A.~Majda.
\newblock A qualitative model for dynamic combustion.
\newblock {\em SIAM Journal of Applied Math}, 41(1):70--93, 1981.

\bibitem{MBOOK}
A.~Majda.
\newblock {\em Compressible Fluid Flows and Systems of Conservation Laws}.
\newblock Springer-Verlag, New York, 1983.

\bibitem{MP}
A.~Majda and R.~L. Pego.
\newblock Stable viscosity matrices for systems of conservation laws.
\newblock {\em J. Differential Equations}, 56(2):229--262, 1985.

\bibitem{Manson}
N.~Manson, C.~Brochet, J.~Brossard, and Y.~Pujol.
\newblock Vibratory phenomena and instability of self-sustained detonations in
  gases.
\newblock In {\em Ninth Symposium (International) on Combustion}, pages
  461--469. Academic Press, 1963.

\bibitem{MZ1}
C.~Mascia and K.~Zumbrun.
\newblock Pointwise {G}reen's function bounds and stability of relaxation
  shocks.
\newblock Preprint, 2001.

\bibitem{MZ2}
C.~Mascia and K.~Zumbrun.
\newblock Stability of viscous shock profiles for dissipative symmetric
  hyperbolic-parabolic systems.
\newblock Preprint, 2001.

\bibitem{MZ3}
C.~Mascia and K.~Zumbrun.
\newblock Pointwise {G}reen's function bounds for shock profiles with
  degenerate viscosity.
\newblock {\em Archive for Rational Mechanics and Analysis}, 2003 (to appear).

\bibitem{MenPlohr}
R.~Menikoff and B.~J. Plohr.
\newblock The {R}iemann problem for fluid flow of real materials.
\newblock {\em Reviews of Modern Physics}, 61(1):75--130, 1999.

\bibitem{Mundy}
G.~Mundy, F.~Ubbelhode, and I.~Wood.
\newblock Fluctuating detonations in gases.
\newblock {\em Proc. Roy. Soc. A}, 306:171--178, 1968.

\bibitem{RV}
J.~Roquejoffre and J.~Vila.
\newblock Stability of {ZND} detonation waves in the {M}ajda combustion model.
\newblock {\em Asymptotic Analysis}, 18:329--348, 1998.

\bibitem{serretransition}
D.~Serre.
\newblock La transition vers l'instabilit\'e pour les ondes de choc
  multi-dimensionnelles.
\newblock {\em Trans. Amer. Math. Soc.}, 353(12):5071--5093 (electronic), 2001.

\bibitem{SeZ}
D.~Serre and K.~Zumbrun.
\newblock Boundary layer stability in real vanishing viscosity limit.
\newblock {\em Comm. Math. Phys.}, 202:547--569, 2001.

\bibitem{SK}
Y.~Shizuta and S.~Kawashima.
\newblock Systems of hyperbolic-parabolic type with applications to the
  discrete {B}oltzmann equation.
\newblock {\em Hokkaido Math. J.}, 14:435--457, 1984.

\bibitem{SS1}
M.~Short and D.~S. Stewart.
\newblock Low-frequency two-dimensional linear instability of plane detonation.
\newblock {\em J. Fluid Mech.}, 340:249--295, 1997.

\bibitem{SS2}
M.~Short and D.~S. Stewart.
\newblock Cellular detonation stability. {I}. {A} normal-mode linear analysis.
\newblock {\em J. Fluid Mech.}, 368:229--262, 1998.

\bibitem{SS3}
M.~Short and D.~S. Stewart.
\newblock The multi-dimensional stability of weak-heat-release detonations.
\newblock {\em J. Fluid Mech.}, 382:109--135, 1999.

\bibitem{Sz}
A.~Szepessy.
\newblock Dynamics and stability of a weak detonation wave.
\newblock {\em Communications in Mathematical Physics}, 202:547--569, 1999.

\bibitem{Szm}
P.~Szmolyan.
\newblock Transversal heteroclinic and homoclinic orbits in singular
  perturbation problems.
\newblock {\em J. Differential Equations}, 92(2):252--281, 1991.

\bibitem{Weyl}
H.~Weyl.
\newblock Shock waves in arbitrary fluids.
\newblock {\em Communications in Pure and Applied Mathematics}, 2, 1949.

\bibitem{Wi}
F.~Williams.
\newblock {\em Combustion Theory}.
\newblock Benjamin/Cummings, Menlo Park, 1985.

\bibitem{Z5}
K.~Zumbrun.
\newblock Stability of large-amplitude shock profiles for compressible
  {N}avier-{S}tokes and {MHD} equations.
\newblock in preparation.

\bibitem{Z1}
K.~Zumbrun.
\newblock Stability of viscous shock waves.
\newblock Lecture Notes, Indiana University, 1998.

\bibitem{Z2}
K.~Zumbrun.
\newblock Multidimensional stability of shock waves.
\newblock Lecture Notes, Indiana University, 2000.

\bibitem{Z4}
K.~Zumbrun.
\newblock Multidimensional stability of planar viscous shock waves.
\newblock In {\em Advances in the Theory of Shock Waves}, number~47 in Progress
  in Nonlinear Differential Equations and Applications, pages 307--516.
  Birkhauser, 2001.

\bibitem{Zapp}
K.~Zumbrun.
\newblock Stability index for relaxation and real viscosity systems, 2002.
\newblock available at math.indiana.edu/home/kzumbrun (corrected appendix of
  \cite{Z4}).

\bibitem{Zhandbook}
K.~Zumbrun.
\newblock Multidimensional stability of shock fronts of compressible
  {N}avier-{S}tokes equations for gas- and magnetohydrodynamics.
\newblock In {\em Handbook of Mathematical Fluid Dynamics IV}. Elsevier, in
  preparation.

\bibitem{ZH}
K.~Zumbrun and P.~Howard.
\newblock Pointwise semigroup methods and the stability of viscous shocks.
\newblock {\em Indiana University Math Journal}, 47(4):741--871, 1998.

\bibitem{ZS}
K.~Zumbrun and D.~Serre.
\newblock Viscous and inviscid stability of multidimensional planar shock
  fronts.
\newblock {\em Indiana University Math Journal}, 48(3):937--999, 1999.

\end{thebibliography}
%
\end{document}